\documentclass[Afour,sageh,times]{sagej}

\usepackage{amsthm}
\usepackage{mathtools}
\usepackage{enumitem}
\usepackage{algorithm}
\usepackage{algpseudocode}
\usepackage{graphicx}
\algnewcommand\algorithmicinput{\textbf{INPUT:}}
\algnewcommand\INPUT{\item[\algorithmicinput]}
\algnewcommand\algorithmicoutput{\textbf{OUTPUT:}}
\algnewcommand\OUTPUT{\item[\algorithmicoutput]}

\usepackage{moreverb,url}

\usepackage[colorlinks,bookmarksopen,bookmarksnumbered,citecolor=red,urlcolor=red]{hyperref}

\newcommand\BibTeX{{\rmfamily B\kern-.05em \textsc{i\kern-.025em b}\kern-.08em
T\kern-.1667em\lower.7ex\hbox{E}\kern-.125emX}}

\newtheorem{theorem}{Theorem}
\newtheorem{lemma}{Lemma}
\newtheorem{prop}[lemma]{Proposition}
\newtheorem{cor}[lemma]{Corollary}
\newtheorem{assum}{Assumption}
\newtheorem{remark}{Remark}
\newtheorem{defn}{Definition}

\def\volumeyear{2019}

\begin{document}

\runninghead{Nishimura and Schwager}

\title{SACBP: Belief Space Planning for Continuous-Time Dynamical Systems via Stochastic Sequential Action Control}

\author{Haruki Nishimura\affilnum{1} and Mac Schwager\affilnum{1}}

\affiliation{\affilnum{1} Stanford University, USA}

\corrauth{Haruki Nishimura, 
          Department of Aeronautics and Astronautics, Stanford University, Stanford, CA 94305, USA.}

\email{hnishimura@stanford.edu}

\begin{abstract}
We propose a novel belief space planning technique for continuous dynamics by viewing the belief system as a hybrid dynamical system with time-driven switching. Our approach is based on the perturbation theory of differential equations and extends Sequential Action Control \cite{ansari2016sac} to stochastic dynamics. The resulting algorithm, which we name SACBP, does not require discretization of spaces or time and synthesizes control signals in near real-time. SACBP is an anytime algorithm that can handle general parametric Bayesian filters under certain assumptions. We demonstrate the effectiveness of our approach in an active sensing scenario and a model-based Bayesian reinforcement learning problem. In these challenging problems, we show that the algorithm significantly outperforms other existing solution techniques including approximate dynamic programming and local trajectory optimization.
\end{abstract}


\maketitle

\section{Introduction}
\label{sec: intro}
Planning under uncertainty still remains as a challenge for robotic systems. Various types of uncertainty, including unmodeled dynamics, stochastic disturbances, and imperfect sensing, significantly complicate problems that are otherwise easy.
For example, suppose that a robot needs to manipulate an object from some initial state to a desired goal. If the mass properties of the object are not known beforehand, the robot needs to simultaneously estimate these parameters and perform control, while taking into account the effects of their uncertainty; the exploration and exploitation trade-off needs to be resolved \cite{slade2017manip}.
On the other hand, uncertainty is quite fundamental in motivating some problems. For instance, a noisy sensor may encourage the robot to carefully plan a trajectory so the observations taken along it are sufficiently informative. This type of problem concerns pure information gathering and is often referred to as active sensing \cite{mihaylova2002activesensing}, active perception \cite{bajcsy1988activeperception}, or informative motion planning \cite{hollinger2014immp}.

A principled approach to address all those problems is to form plans in the belief space, where the planner chooses sequential control inputs based on the evolution of the belief state. This approach enables the robot to appropriately execute controls under stochasticity and partial observability since they are both incorporated into the belief state. Belief space planning is also well suited for generating information gathering actions \cite{platt2010bsp}.

This paper proposes a novel online belief space planning algorithm. It does not require discretization of the state space or the action space, and can directly handle continuous-time system dynamics. The algorithm optimizes the expected value of a first-order cost reduction with respect to a nominal control policy at every re-planning time, proceeding in a receding-horizon fashion. We are inspired by the Sequential Action Control (SAC) algorithm recently proposed in \cite{ansari2016sac} for model-based deterministic optimal control problems.  SAC is an online method to synthesize control signals in real time for challenging (but deterministic) physical systems such as a cart pendulum and a spring-loaded inverted pendulum. Based on the concept of SAC, this paper develops an algorithmic framework to control stochastic belief systems whose dynamics are governed by parametric Bayesian filters.

\subsection{Related Work in Belief Space Planning}
\label{sec: related work}
There is a large body of literature in belief space planning. Below we briefly review relevant work in three major types of solution methods: greedy strategies, trajectory optimization methods, and belief MDP and POMDP approaches. We then change our perspective and discuss advantages and drawbacks of closed-loop and open-loop planning schemes, while also introducing other relevant work that do not necessarily fall into the three categories mentioned above.

\paragraph{Greedy Strategies}
Belief space planning is known to be challenging for a couple of reasons. First, the belief state is continuous and can be high-dimensional even if the underlying state space is small or discrete. Second, the dynamics that govern the belief state transitions are stochastic due to unknown future observations. Greedy approaches alleviate the complexity by ignoring long-term effects and solve single-shot decision making problems sequentially. Despite their suboptimality for long-term planning, these methods are often employed to find computationally tractable solutions and achieve reasonable performance in different problems \cite{bourgault2002information,seekircher2011entropy,schwager2017info}, especially in the active sensing domain.

\paragraph{Trajectory Optimization Methods}
In contrast to the greedy approaches, trajectory optimization methods take into account multiple timesteps at once and find non-myopic solutions. In doing so, it is often assumed that the maximum likelihood observation (MLO) will always occur during planning \cite{platt2010bsp,erez2010continuouspomdps,patil2014gaussianbsp}. This heuristic assumption yields a deterministic optimal control problem, to which various nonlinear trajectory optimization algorithms are applicable. However, ignoring the effects of stochastic future observations can degrade the performance \cite{vandenberg2012beliefilqg}. Other methods \cite{vandenberg2012beliefilqg, rafieisakhaei2017tlqg} that do not rely on the MLO assumption are advantageous in that regard. In particular, belief iLQG \cite{vandenberg2012beliefilqg} performs iterative local optimization in a Gaussian belief space by quadratically approximating the value function and linearizing the dynamics to obtain a time-varying affine feedback policy. However, this method as well as many other solution techniques in this category result in multiple iterations of intensive computation and can require a significant amount of time until convergence.

\paragraph{Belief MDP and POMDP Approaches}
Belief space planning can be modeled as a Markov decision process (MDP) in the belief space, given that the belief state transition is Markovian. If the reward (or cost) is defined as an explicit function of the state and the control, the problem is equivalent to a partially observable Markov decision process (POMDP) \cite{kaelbling1998planning}.
A key challenge in POMDPs and belief MDPs has been to address problems with large state spaces. This is particularly important in belief MDPs since the state space for a belief MDP is a continuous belief space. To handle continuous spaces, \cite{couetoux2011cuct} introduce double progressive widening (DPW) for Monte Carlo Tree Search (MCTS) \cite{browne2012mcts}. In \cite{slade2017manip}, this MCTS-DPW algorithm is run in a Gaussian belief space to solve the object manipulation problem mentioned earlier. We have also presented a motion-based communication algorithm in our prior work, which uses MCTS-DPW for active intent inference with monocular vision \cite{nishimura2018activembc}.

While MCTS-DPW as well as other general purpose POMDP methods \cite{somani2013despot,sunberg2017pomcpow} are capable of handling continuous state spaces, their algorithmic concepts are rooted in dynamic programming and tree search, requiring a sufficient amount of exploration in the tree. The tree search technique also implicitly assumes discrete-time transition models. In fact, most prior works reviewed in this section are intended for discrete-time systems. A notable exception is \cite{chaudhari2013sampling}, in which a sequence of discrete-time POMDP approximations is constructed such that it converges to a true continuous-time model. However, this approach still relies on an existing POMDP solver \cite{kurniawati2008sarsop} that is designed for discrete systems. There still remains a need for an efficient and high-performance belief space planning algorithm that is capable of directly handling systems with inherently continuous-space, continuous-time dynamics, such as maneuvering micro-aerial vehicles, or autonomous cars at freeway speeds.

\paragraph{Closed-loop and Open-loop Planning}
An important aspect of belief space planning is the stochastic belief dynamics, which essentially demand that one plan over closed-loop (i.e. feedback) policies. This is a primary assumption in belief MDP and POMDP frameworks where an optimal mapping from beliefs to control actions are sought \cite{huber2009probabilistic, kochenderfer2015decision}. However, computing exact optimal policies in general is intractable due to curse of dimensionality and curse of history \cite{papadimitriou1987complexity, madani1999undecidability, somani2013despot}. Therefore, in practice we need to make certain approximations to achieve a tractable solution within a reasonable computation time, whether the solution method is online or offline. For POMDP methods, those approximations are often done by discretization and/or sampling \cite{kurniawati2008sarsop, chaudhari2013sampling, somani2013despot, sunberg2017pomcpow}. Belief iLQG \cite{vandenberg2012beliefilqg} makes an approximation by restricting the policy to a time-varying affine feedback controller and performing local optimization. Other methods that provide closed-loop policies include Feedback-based Information RoadMap (FIRM) \cite{agha2014firm}, its extension with online local re-planning \cite{agha2018slap}, and T-LQG \cite{rafieisakhaei2017tlqg}. FIRM-based methods, which are developed primarily for motion planning problems, construct a graph in the belief space whose edges correspond to local LQG controllers. T-LQG approximately decouples belief space planning into local trajectory optimization and LQG tracking via a use of the separation principle. The local trajectory optimization part achieves one of the lowest asymptotic computational complexities among existing trajectory optimization methods, and thus can be solved efficiently with a commercial nonlinear programming (NLP) solver \cite{rafieisakhaei2017tlqg}.

Even though closed-loop methods are appealing, their computation cost can be prohibitive in some challenging planning and control problems, especially when the state space is large and yet real-time performance is required. In such cases, a practical strategy to alleviate the computational burden is open-loop planning, wherein one seeks to optimize a static sequence of control actions. Often times feedback is provided through re-planning of a fixed horizon, open-loop control sequence after making a new observation. This combination of open-loop planning and feedback through re-planning is called receding horizon control (RHC) or model predictive control (MPC) \cite{huber2009probabilistic}. Although they do not account for feedback at each planning time, RHC methods in general have successfully solved challenging planning and control problems with high efficiency where fast closed-loop policy computation can be impractical. Examples of such real-time applications include trajectory tracking with quadrotors \cite{bangura2014real} and agile autonomous driving on dirt \cite{williams2018information}, although they are both in the optimal control literature. Within belief space planning, some methods based on the MLO assumption \cite{platt2013convex, patil2014gaussianbsp} are run as RHC.

\subsection{Contributions}
Our approach presented in this paper takes the form of RHC, which is reasonable for high-dimensional and continuous belief space planning problems that also demand highly-dynamic, real-time maneuvers with high frequency observation updates. However, the proposed method is significantly different than any of the previous approaches discussed in Section \ref{sec: related work}. We view the stochastic belief dynamics as a hybrid system with time-driven switching \cite{heemels2009hybrid}, where the controls are applied in continuous time and the observations are made in discrete time. A discrete-time observation creates a jump discontinuity in the belief state trajectory due to a sudden Bayesian update of the belief state. This view of belief space planning yields a continuous-time optimal control problem of a high-dimensional hybrid system. We then propose a model-based control algorithm to efficiently compute control signals in a receding-horizon fashion. The algorithm is based on Sequential Action Control (SAC) \cite{ansari2016sac}. SAC in its original form is a deterministic, model-based hybrid control algorithm, which ``perturbs'' a nominal control trajectory in a structured way so that the cost functional is optimally reduced up to the first order. The key to this approach is a careful use of the perturbation theory of differential equations that is often discussed in the mode scheduling literature \cite{egerstedt2006transition,Wardi2012switched}.
As a result, SAC derives the optimal perturbation in closed form and synthesizes control signals at a high frequency to achieve a significant improvement over other optimal control methods that are based on local trajectory optimization \cite{ansari2016sac}.


We apply the perturbation theory to parametric Bayesian filters and derive the optimal control perturbation using the framework of SAC. Even though each control perturbation is small, high-frequency control synthesis and online re-planning yield a series of control actions that is significantly different than the nominal control, reacting to stochastic observations collected during execution or online changing conditions. Furthermore, we extend the original SAC algorithm to also account for stochasticity in the future observations during planning, by incorporating Monte Carlo sampling of nominal belief trajectories. Our key contribution is the resulting continuous belief space planning algorithm, which we name SACBP.
The algorithm has the following desirable properties:
\begin{enumerate}
	\item Although the form of control perturbation is open-loop with on-line re-planning, the perturbation computed by SACBP in near real-time is optimized for better average performance over the planning horizon than a given nominal control, whether it is open-loop or closed-loop. 

	\item SACBP does not require discretization of the state space, the observation space, or the control space. It also does not require discretization of time other than for numerical integration purposes.
	\item General nonlinear parametric Bayesian filters can be used for state estimation as long as the system is control-affine and the control cost is quadratic.
	\item Stochasticity in the future observations are fully considered.
	\item SACBP is an anytime algorithm.  Furthermore, the Monte Carlo sampling part of the algorithm is naturally parallelizable.
	\item Even though SACBP is inherently suboptimal for the original stochastic optimal control problem, empirical results suggest that it is highly sample-efficient and outperforms other open-loop and closed-loop methods when near real-time performance is required.
\end{enumerate}

There exists prior work \cite{mavrommati2018ergodic} that uses SAC for active sensing, but its problem formulation relies on the ergodic control framework, which is significantly different from the belief space planning framework we propose here. We show that our SACBP outperforms projection-based ergodic trajectory optimization, MCTS-DPW, T-LQG, and a greedy method on an active multi-target tracking example. We also show that SACBP outperforms belief iLQG, MCTS-DPW, and T-LQG on a manipulation scenario.

This paper is an extension of the theory and results previously presented by the authors in \cite{nishimura2018sacbp}. Compared to the conference version, we provide a more detailed derivation of the algorithm (Section \ref{sec: SACBP}) as well as a thorough mathematical analysis of the open-loop control perturbation for stochastic hybrid systems (Section \ref{sec: analysis} and Appendix \ref{sec: appendix}). This analysis leads to a guarantee for SACBP that, with an appropriate choice of the perturbation duration, the algorithm is expected to perform no worse than the nominal control. Since the nominal control can be arbitrary, one could even provide a discrete POMDP policy derived offline as a nominal policy to ``warm-start" the planning. We also report new simulation results that involve a comparison of SACBP with T-LQG (Section \ref{sec: results}), where we observe superior performance of SACBP for real-time applications.

In the next section we derive relevant equations and present the SACBP algorithm along with a discussion on computational complexity. Section \ref{sec: analysis} provides the key results of the mathematical analysis. Section \ref{sec: results} summarizes the simulation results. Conclusions and future work are presented in Section \ref{sec: conclusions}.

\section{SACBP Algorithm}
\label{sec: SACBP}
We first consider the case where some components of the state are fully observable. We begin with this mixed observability case as it is simpler to explain, yet still practically relevant. For example, this is a common assumption in various active sensing problems \cite{schwager2017info, leny2009active, popovic2017uav} where the state of the robot is perfectly known, but some external variable of interest (e.g. a target's location) is stochastic.  In addition, deterministic state transitions are often assumed for the robot. Therefore, in Section \ref{subsec: mixed observability} we derive the SACBP control update formulae for this case. The general belief space planning where none of the state is fully observable or deterministically controlled is discussed in Section \ref{subsec: general belief space planning}. An extension to use a closed-loop policy as the nominal control is presented in Section \ref{subsec: closed-loop policy}. The computation time complexity is discussed in Section \ref{subsec: runtime analysis}.

\subsection{Problems with Mixed Observability}
\label{subsec: mixed observability}
Suppose that a robot can fully observe and deterministically control its own state $p(t) \in \mathbb{R}^{n_p}$. Other external states over which the robot does not have direct control are not known and are estimated with the belief vector $b(t) \in \mathbb{R}^{n_b}$. This belief vector characterizes a probability distribution that the robot uses for state estimation. If the belief is Gaussian, for example, the covariance matrix can be vectorized column-wise and stacked all together with the mean to form the belief vector. 
We define the augmented state as $s \triangleq (p^\mathrm{T},b^\mathrm{T})^\mathrm{T} \in \mathbb{R}^{n_s}$.

\subsubsection{Dynamics Model}
The physical state $p$ is described by the following ODE:
\begin{align}
\dot{p}(t) = f\left(p(t),u(t)\right), \label{eq: physDynamics}
\end{align}
where $u(t) \in \mathbb{R}^m$ is the control signal.
On the other hand, suppose that the belief state only changes in discrete time upon arrival of a new observation from the sensors. In the case of target tracking, for example, this means that the change of the location of the target is described by a discrete time model from the robot's perspective. The behavior of such a system can be estimated by a discrete-time Bayesian filter. We will discuss the more general continuous-discrete time filtering case in Section \ref{subsec: general belief space planning}. Let $t_k$ be the time when the
$k$-th observation becomes available to the robot. The belief state transition is given by
\begin{align}
\begin{cases}
b(t_k) = g(p(t^-_k), b(t^-_k),y_k)  \\ 
b(t) = b(t_k) &\forall t \in [t_k, t_{k+1}), \label{eq: belDynamics}
\end{cases}
\end{align}
where $t^-_k$ is infinitesimally smaller than $t_k$.
Nonlinear function $g$ corresponds to a discrete-time, parametric Bayesian filter (e.g., Kalman filter, extended Kalman filter, discrete Bayesian filter, etc.) that forward-propagates the belief for prediction, takes the new observation $y_k \in \mathbb{R}^q$, and returns the updated belief state.
The concrete choice of the filter depends on the instance of the problem.
Note that the belief state stays constant for the most of the time in \eqref{eq: belDynamics}. This is because we are viewing a discrete time model in a continuous time framework.

Equations \eqref{eq: physDynamics} and \eqref{eq: belDynamics} constitute a hybrid system with time-driven switching \cite{heemels2009hybrid}. This hybrid system representation is practical since it captures the fact that the observation updates occur less frequently than the control actuation in general, due to expensive information processing of sensor readings. Furthermore, with this representation one can naturally handle agile robot dynamics as they are without coarse discretization in time.

Given the initial state $s_0 \triangleq (p(t_0)^\mathrm{T}, b(t_0)^\mathrm{T})^\mathrm{T}$ and some control $u(t)$ for $t \in [t_0, t_f]$, the system evolves stochastically according to the hybrid dynamics equations. The stochasticity is due to a sequence of stochastic future observations that will be taken by the end of the planning horizon $t_f$. In this paper we assume that the observation interval $t_{k+1} - t_k \triangleq \Delta t_o$ is fixed, and the control signals are recomputed when a new observation is incorporated in the belief, although one can also use a variable observation interval.

\subsubsection{Perturbed Dynamics}
The control synthesis of SACBP begins with a given nominal control schedule $u(t)$ for $t \in [t_0, t_f]$. For simplicity we assume here that the nominal control schedule is open-loop, but we remind the reader that SACBP can also handle closed-loop nominal policies, which we discuss in Section \ref{subsec: closed-loop policy}. Suppose that the nominal control is applied to the system and a sequence of $T$ observations $(y_{1},\dots,y_{T})$ is obtained. Conditioned on the observation sequence, the augmented state evolves deterministically. Let $s = (p^\mathrm{T},b^\mathrm{T})^\mathrm{T}$ be the nominal trajectory of the augmented state induced by $(y_{1},\dots,y_{T})$. 

Now let us consider perturbing the nominal trajectory at a fixed time $\tau < t_1$ for a short duration $\epsilon$. The perturbed control $u^\epsilon$ is defined as
\begin{align}
u^\epsilon(t) \triangleq \begin{cases}
v &\text{if $t \in (\tau - \epsilon, \tau]$} \\
u(t) & \text{otherwise.}
\end{cases}
\label{eq: perturbed control trajectory}
\end{align}
Therefore, the control perturbation is determined by the nominal control $u(t)$, the tuple $(\tau, v)$, and $\epsilon$.
Given $(\tau, v)$, the resulting perturbed system trajectory can be written as
\begin{align}
\begin{cases}
p^{\epsilon}(t) \triangleq p(t) + \epsilon\Psi_p(t) + o(\epsilon) \\
b^\epsilon(t) \triangleq b(t) + \epsilon\Psi_b(t) + o(\epsilon),
\end{cases}
\end{align}
where $\Psi_p(t)$ and $\Psi_b(t)$ are the state variations that are linear in the perturbation duration $\epsilon$:
\begin{align}
\Psi_p(t) = \frac{\partial_+}{\partial\epsilon} p^{\epsilon}(t)\bigg|_{\epsilon=0} \triangleq \lim_{\epsilon\rightarrow 0^+} \frac{p^\epsilon(t) - p(t)}{\epsilon}\\
\Psi_b(t) = \frac{\partial_+}{\partial\epsilon} b^{\epsilon}(t)\bigg|_{\epsilon=0} \triangleq \lim_{\epsilon\rightarrow 0^+} \frac{b^\epsilon(t) - b(t)}{\epsilon}.
\end{align}
The notation $\frac{\partial_+}{\partial\epsilon}$ represents the right derivative with respect to $\epsilon$.
The state variations at perturbation time $\tau$ satisfy
\begin{align}
\begin{cases}
\Psi_p(\tau) = f(p(\tau),v) - f(p(\tau),u(\tau)) \\
\Psi_b(\tau) = 0.
\end{cases}
\end{align}
The initial belief state variation $\Psi_b(\tau)$ is zero because the control perturbation $u^\epsilon$ has no effect on the belief state until the first Bayesian update is performed at time $t_1$, according to the hybrid system model \eqref{eq: physDynamics}\eqref{eq: belDynamics}.
For $t \geq \tau$, the physical state variation $\Psi_p$ evolves according to the following first-order ODE:
\begin{align}
\dot{\Psi}_p(t) &= \frac{d}{dt}\left(\frac{\partial_+}{\partial\epsilon} p^{\epsilon}(t)\bigg|_{\epsilon=0}\right) \\
&= \frac{\partial_+}{\partial\epsilon} \dot{p}^{\epsilon}(t)\bigg|_{\epsilon=0} \\
&= \frac{\partial_+}{\partial\epsilon} f(p^{\epsilon}(t),u(t))\bigg|_{\epsilon=0} \\
&= \frac{\partial}{\partial p} f\left(p(t),u(t)\right) \Psi_p(t), \label{eq: variational equation p}
\end{align}
where the chain rule of differentiation and $p^{\epsilon}(t)\vert_{\epsilon=0} = p(t)$ are used in \eqref{eq: variational equation p}. For a more rigorous analysis, see Appendix \ref{sec: appendix}.
The dynamics of the belief state variation $\Psi_b$ in the continuous region $t \in [t_k, t_{k+1})$ satisfy $\dot{\Psi}_b(t) = 0$ since the belief vector $b(t)$ is constant according to \eqref{eq: belDynamics}.
However, across the jumps the belief state variation $\Psi_b$ changes discontinuously and satisfies
\begin{align}
\Psi_b(t_{k}) &= \frac{\partial_+}{\partial\epsilon}b^{\epsilon}(t_{k})\bigg|_{\epsilon=0} \\
&= \frac{\partial_+}{\partial\epsilon}g\left(p^{\epsilon}(t_k^-),b^{\epsilon}(t_k^-), y_k\right)\bigg|_{\epsilon=0} \\
&\begin{multlined}= \frac{\partial}{\partial p} g\left(p(t_k^-),b(t_k^-),y_k\right) \Psi_p(t_k^-) \\ + \frac{\partial}{\partial b} g\left(p(t_k^-),b(t_k^-),y_k\right) \Psi_b(t_k^-).
\end{multlined}
\label{eq: variational equation b}
\end{align}

\subsubsection{Perturbed Cost Functional}
Let us consider a total cost of the form
\begin{align}
\int_{t_0}^{t_f} c\left(p(t),b(t),u(t)\right) dt + h(p(t_f),b(t_f)),
\end{align}
where $c$ is the running cost and $h$ is the terminal cost.
Following the discussion above on the perturbed dynamics, let $J$ denote the total cost of the nominal trajectory conditioned on the given observation sequence $(y_{1},\dots,y_{T})$. Under the fixed $(\tau, v)$, we can represent the perturbed cost $J^{\epsilon}$ in terms of $J$ as
\begin{align}
J^{\epsilon} \triangleq J + \epsilon \nu(t_f) + o(\epsilon),
\label{eq: Cost variation def}
\end{align}
where $\nu(t_f) \triangleq \frac{\partial_+}{\partial\epsilon} J^{\epsilon}|_{\epsilon=0}$ is the variation of the total cost with respect to the perturbation.
For further analysis it is convenient to express the running cost in the Mayer form \cite{liberzon2012calculus}. Let $\hat{s}(t)$ be a new state variable defined by $\dot{\hat{s}}(t) = c\left(p(t),b(t),u(t)\right)$ and $\hat{s}(t_0) = 0$. Then the total cost is a function of the appended augmented state $\bar{s} \triangleq (\hat{s},s^\mathrm{T})^\mathrm{T}\in \mathbb{R}^{1+n_s}$ at time $t_f$, which is given by
\begin{align}
J = \hat{s}(t_f) + h\left(s(t_f)\right).
\end{align}
Using this form of the total cost $J$, the perturbed cost \eqref{eq: Cost variation def} becomes
\begin{align}
J^{\epsilon} = J + \epsilon \begin{bmatrix} 1 \\ \frac{\partial}{\partial p}h\left(p(t_f),b(t_f)\right) \\ \frac{\partial}{\partial b}h\left(p(t_f),b(t_f)\right)\end{bmatrix}^\mathrm{T} \overline{\Psi}(t_f) + o(\epsilon), \label{eq: Cost variation 2}
\end{align}
where $\overline{\Psi}(t_f) \triangleq \left(\hat{\Psi}(t_f), \Psi_p(t_f)^\mathrm{T}, \Psi_b(t_f)^\mathrm{T}\right)^\mathrm{T}$ and $\hat{\Psi}$ is the variation of $\hat{s}$. Note that the dot product in \eqref{eq: Cost variation 2} corresponds to $\nu(t_f)$ in \eqref{eq: Cost variation def}. The variation $\hat{\Psi}$ follows the variational equation for $t \geq \tau$:
\begin{align}
\dot{\hat{\Psi}}(t) &= \frac{d}{dt}\left(\frac{\partial_+}{\partial\epsilon} \hat{s}^\epsilon(t)\bigg|_{\epsilon=0}\right)\\
&= \frac{\partial_+}{\partial\epsilon} \dot{\hat{s}}^\epsilon(t)\bigg|_{\epsilon=0} \\
&= \frac{\partial_+}{\partial\epsilon} c(p^{\epsilon}(t),b^{\epsilon}(t),u(t))\bigg|_{\epsilon=0} \\
&\begin{multlined}
    = \frac{\partial}{\partial p}c(p(t),b(t),u(t))^\mathrm{T} \Psi_p(t) \\ + \frac{\partial}{\partial b} c(p(t),b(t),u(t))^\mathrm{T} \Psi_b(t), \label{eq: variational equation 0}
\end{multlined}
\end{align}
where the initial condition is given by $\hat{\Psi}(\tau) = c(p(\tau), b(\tau), v) - c(p(\tau), b(\tau), u(\tau))$.

The perturbed cost equation \eqref{eq: Cost variation 2}, especially the dot product expressing $\nu(t_f)$, is consequential; it tells us how the total cost changes due to the perturbation applied at some time $\tau$, up to the first order with respect to the perturbation duration $\epsilon$.
At this point, one could compute the value of $\nu(t_f)$ for a control perturbation with a specific value of $(\tau,v)$ by simulating the nominal dynamics and integrating the variational equations \eqref{eq: variational equation p}\eqref{eq: variational equation b}\eqref{eq: variational equation 0} from $\tau$ up to $t_f$.
\subsubsection{Adjoint Equations}
Unfortunately, this forward integration of $\nu(t_f)$ is not so useful by itself since we are interested in finding the value of $(\tau, v)$ that achieves the largest possible negative $\nu(t_f)$, if it exists; it would be computationally intensive to apply control perturbation at different application times $\tau$ with different values of $v$ and re-simulate the state variation $\overline{\Psi}$. To avoid this computationally expensive search, we mirror the approach presented in \cite{ansari2016sac} and introduce the adjoint system $\overline{\rho} \triangleq (\hat{\rho}, \rho_p^{\mathrm{T}}, \rho_b^{\mathrm{T}})^{\mathrm{T}}$ with which the dot product remains invariant:
\begin{align}
    \frac{d}{dt}\left(\overline{\rho}(t)^\mathrm{T}\overline{\Psi}(t)\right) = 0 \;\;\forall t \in [t_0,t_f].
\end{align}
If we let 
\begin{multline}
    \label{eq: adjoint_boundary}
    \overline{\rho}(t_f) \triangleq  \bigg(1, \frac{\partial}{\partial p}h\left(p(t_f),b(t_f)\right)^\mathrm{T}, \\ \frac{\partial}{\partial b}h\left(p(t_f),b(t_f)\right)^\mathrm{T}\bigg)^\mathrm{T}
\end{multline}
so that its dot product with $\overline{\Psi}(t_f)$ equals $\nu(t_f)$ as in \eqref{eq: Cost variation 2}, the time invariance gives
\begin{align}
\nu(t_f) &= \overline{\rho}(t_f)^\mathrm{T} \overline{\Psi}(t_f) \\
& = \overline{\rho}(\tau)^\mathrm{T} \overline{\Psi}(\tau) \\
& = \overline{\rho}(\tau)^\mathrm{T} \begin{bmatrix} c\left(p(\tau),b(\tau),v\right) - c\left(p(\tau),b(\tau),u(\tau)\right) \\ f(p(\tau),v) - f(p(\tau),u(\tau)) \\ 0 \end{bmatrix}.
\end{align}
Therefore, we can compute the first-order cost change $\nu(t_f)$ for different values of $\tau$ once the adjoint trajectory is derived.
For $t \in[t_k, t_{k+1})$ the time derivative of $\overline{\Psi}$ exists, and the invariance property leads to the following equation:
\begin{align}
    \dot{\overline{\rho}}(t)^\mathrm{T}\overline{\Psi}(t) + \overline{\rho}(t)^\mathrm{T}\dot{\overline{\Psi}}(t) = 0. \label{eq: Continuous adjoint condition}
\end{align}
It can be verified that the following system satisfies \eqref{eq: Continuous adjoint condition} with $\overline{\rho}(t) = \left(\hat{\rho}(t), \rho_p(t)^\mathrm{T}, \rho_b(t)^\mathrm{T}\right)^\mathrm{T}$:
\begin{align}
\begin{cases}
\dot{\hat{\rho}}(t) = 0 \\
\dot{\rho}_p(t) = -\frac{\partial}{\partial p} c(p(t),b(t),u(t)) - \frac{\partial}{\partial p} f(p(t),u(t))^\mathrm{T} \rho_p(t) \\
\dot{\rho}_b(t) = -\frac{\partial}{\partial b} c(p(t),b(t),u(t)).
\end{cases} \label{eq: Continuous adjoint system}
\end{align}
Analogously, across discrete jumps we can still enforce the invariance by setting
$\overline{\rho}(t_k)^\mathrm{T}\overline{\Psi}(t_k) = \overline{\rho}(t_k^-)^\mathrm{T}\overline{\Psi}(t_k^-)$, 
which holds for the following adjoint equations:
\begin{align}
\begin{cases}
\hat{\rho}(t_k^-) = \hat{\rho}(t_k) \\
\rho_p(t_k^-) = \rho_p(t_k) + \frac{\partial}{\partial p}g\left(p(t_k^-),b(t_k^-),y_k\right)^\mathrm{T}\rho_b(t_k)\\
\rho_b(t_k^-) = \frac{\partial}{\partial b}g\left(p(t_k^-),b(t_k^-),y_k\right)^\mathrm{T}\rho_b(t_k).
\end{cases} \label{eq: Discrete adjoint system}
\end{align}
Note that the adjoint system integrates backward in time as it has the boundary condition \eqref{eq: adjoint_boundary} defined at $t_f$ . More importantly, the adjoint dynamics \eqref{eq: Continuous adjoint system}\eqref{eq: Discrete adjoint system} only depend on the nominal trajectory of the system $(p,b)$ and the observation sequence $(y_{1},\dots,y_{T})$.
Considering that $\hat{\rho}(t) = 1$ at all times, the cost variation term $\nu(t_f)$ is finally given by
\begin{multline}
\nu(t_f) = c(p(\tau),b(\tau),v) - c(p(\tau),b(\tau),u(\tau)) + \\ \rho_p(\tau)^\mathrm{T}\left\{f(p(\tau),v) - f(p(\tau),u(\tau))\right\}. \label{eq: general mixed observability linear variational cost}
\end{multline}

\subsubsection{Control Optimization}
In order to efficiently optimize \eqref{eq: general mixed observability linear variational cost} with respect to $(\tau, v)$, we assume hereafter that the control cost is additive quadratic $\frac{1}{2}u^\mathrm{T} C_u u$ and the dynamics model $f(p,u)$ is control-affine with linear term $H(p)u$, where $H\colon \mathbb{R}^{n_p} \rightarrow \mathbb{R}^m$ can be any nonlinear mapping. Although the control-affine assumption may appear restrictive, many physical systems possess this property in engineering practice. As a result of these assumptions, \eqref{eq: general mixed observability linear variational cost} becomes
\begin{multline}
\nu(t_f) = \frac{1}{2} v^\mathrm{T} C_u v + \rho_p(\tau)^\mathrm{T}H(p(\tau))(v - u(\tau)) \\- \frac{1}{2} u(\tau)^\mathrm{T} C_u u(\tau).
\label{eq: specific linear variational cost}
\end{multline}

So far we have treated the observation sequence $(y_{1},\dots,y_{T})$ as given and fixed. However, in practice it is a random process that we have to take into account. Fortunately, our control optimization is all based on the nominal control schedule $u(t)$, with which we can both simulate the augmented dynamics and sample the observations. To see this, let us consider the observations as a sequence of random vectors $(Y_1,\dots,Y_T)$ and rewrite $\nu(t_f)$ in \eqref{eq: specific linear variational cost} as $\nu(t_f, Y_{1},\dots,Y_{T})$ to clarify the dependence on it. The expected value of the first order cost variation is given by
\begin{align}
\mathbb{E}[\nu(t_f)] = \int \nu(t_f, Y_{1},\dots,Y_{T}) d\mathbb{P},
\end{align}
where $\mathbb{P}$ is the probability measure associated with these random vectors. 
Although we do not know the specific values of $\mathbb{P}$, we have the generative model; we can simulate the augmented state trajectory using the nominal control and sequentially sample the stochastic observations from the belief states along the trajectory.
 
Using the linearity of expectation for \eqref{eq: specific linear variational cost}, we have
\begin{multline}
\mathbb{E}[\nu(t_f)] = \frac{1}{2}v^\mathrm{T}C_u v + \mathbb{E}[\rho_p(\tau)]^\mathrm{T}H(p(\tau))(v - u(\tau)) \\ - \frac{1}{2} u(\tau)^\mathrm{T} C_u u(\tau). \label{eq: expected linear variational cost}
\end{multline}  
Notice that only the adjoint trajectory is stochastic. We can employ Monte Carlo sampling to sample a sufficient number of observation sequences to approximate the expected adjoint trajectory.
Now \eqref{eq: expected linear variational cost} becomes a convex quadratic in $v$ for a positive definite $C_u$. Assuming that $C_u$ is also diagonal, analytical solutions are available to the following convex optimization problem with an input saturation constraint.
\begin{equation}
\begin{aligned}
&\underset{v}{\text{minimize}} &\mathbb{E}[\nu(t_f)] \\
&\text{subject to} &a \preceq v \preceq b,
\end{aligned}
\label{eq: convex program}
\end{equation}
where $a, b \in \mathbb{R}^m$ are some saturation vectors and $\preceq$ is an elementwise inequality.
This optimization is solved for different values of $\tau \in (t_0 + t_\mathrm{calc} + \epsilon, t_0 + t_\mathrm{calc} + \Delta t_o)$, where $t_\mathrm{calc}$ is a pre-allocated computation time budget and $\Delta t_o$ is the time interval between two successive observations as well as control updates. We then search over $(\tau,v^*(\tau))$ for the optimal perturbation time $\tau^*$ to globally minimize $\mathbb{E}[\nu(t_f)]$.  There is only a finite number of such $\tau$ to consider since in practice we use numerical integration such as the Euler scheme with some step size $\Delta t_c$ to compute the trajectories.
In \cite{ansari2016sac} the finite perturbation duration $\epsilon$ is also optimized using line search, but in this work we set $\epsilon$ as a tunable parameter to reduce the computation time. The complete algorithm is summarized in Algorithm \ref{algo: sacbp_algo_1}. The call to the algorithm occurs every $\Delta t_o$[s] in a receding-horizon fashion, after the new observation is incorporated in the belief.

\begin{algorithm}[t]
	\caption{SACBP Control Update for Problems with Mixed Observability}\label{algo: sacbp_algo_1}
	\begin{algorithmic}[1]
		\INPUT Current augmented state $s_0 = (p(t_0)^\mathrm{T},b(t_0)^\mathrm{T})^\mathrm{T}$, nominal control schedule $u(t)$ for $t \in [t_0, t_f]$, perturbation duration $\epsilon$
		\OUTPUT Optimally perturbed control schedule $u^{\epsilon}(t)$ for $t \in [t_0, t_f]$
		\For{$i$ = 1:$N$}
		\State Forward-simulate nominal augmented state trajectory \eqref{eq: physDynamics}\eqref{eq: belDynamics} and sample observation sequence $(y^{i}_{1},\dots,y^{i}_{T})$ along the augmented state trajectory.
		\State Backward-simulate nominal adjoint trajectory $\rho_p^{i}$, $\rho_b^{i}$ \eqref{eq: Continuous adjoint system}\eqref{eq: Discrete adjoint system} with sampled observations.
		\EndFor
		\State Compute Monte Carlo estimate: $\mathbb{E}[\rho_p] \approx \frac{1}{N}\sum_{i=1}^N \rho^{i}_p$.
		\For{($\tau = t_0 + t_\mathrm{calc} + \epsilon;~\tau \leq t_0 + t_\mathrm{calc} + \Delta t_o;~\tau \leftarrow \tau + \Delta t_c$)}
		\State Solve quadratic minimization \eqref{eq: convex program} with \eqref{eq: expected linear variational cost}. Store optimal value $\nu^*(\tau)$ and optimizer $v^*(\tau)$. 
		\EndFor
		\State $\tau^* \leftarrow\arg\min \nu^*(\tau), \;v^* \leftarrow v^*(\tau^*)$
		\State $u^{\epsilon} \leftarrow PerturbControlTrajectory(u,v^*,\tau^*,\epsilon)$ \eqref{eq: perturbed control trajectory}
		\State \textbf{return} $u^{\epsilon}$
	\end{algorithmic}
\end{algorithm}

\subsection{General Belief Space Planning Problems}
\label{subsec: general belief space planning}
If none of the state is fully observable, the same stochastic SAC framework still applies almost as is to the belief sate $b$. In this case we consider a continuous-discrete filter \cite{xie2007estimation} where the prediction step follows an ODE and the update step provides an instantaneous discrete jump. The hybrid dynamics for the belief vector are given by
\begin{align}
\begin{cases}
b(t_{k}) = g(b(t^-_k),y_k)  \\ 
\dot{b}(t) = f(b(t),u(t)) &\forall t \in [t_k, t_{k+1}). \label{eq: continuous discrete belDynamics}
\end{cases}
\end{align}
Letting $\Psi(t) = \frac{\partial_+}{\partial\epsilon} b^{\epsilon}(t)\big\vert_{\epsilon=0}$, the variational equation yields
\begin{align}
    \begin{cases}
        \Psi(t_{k}) = \frac{\partial}{\partial b} g(b(t_k^-), y_k) \Psi(t_k^-) \\
        \dot{\Psi}(t) = \frac{\partial}{\partial b} f(b(t), u(t)) \Psi(t) &\forall t \in [t_k, t_{k+1})
    \end{cases}
\end{align}
with initial condition $\Psi(\tau) = f(b(\tau), v) - f(b(\tau), u(\tau))$. 

Let the total cost be of the form:
\begin{align}
    \int_{t_0}^{t_f} c(b(t), u(t)) dt + h(b(t_f)).
\end{align}
Under the given $(\tau, v)$ and $(y_1,\dots,y_T)$, the variation $\nu(t_f)$ of the total cost can be computed as
\begin{multline}
    \label{eq: total_cost_variation_forward_integration}
    \nu(t_f) = c(b(\tau), v) - c(b(\tau), u(\tau)) \\ + \int_{\tau}^{t_f} \frac{\partial}{\partial b} c(b(t), u(t))^{\mathrm{T}} \Psi(t) dt \\
    + \frac{\partial}{\partial b} h(b(t_f))^\mathrm{T}\Psi(t_f).
\end{multline}
This is equivalent to 
\begin{multline}
     \nu(t_f) = c(b(\tau), v) - c(b(\tau), u(\tau)) \\ + \rho(\tau)^\mathrm{T} \left\{f(b(\tau), v) - f(b(\tau), u(\tau)\right\},
     \label{eq: general belief linear variational cost}
\end{multline}
where $\rho$ is the adjoint system that follows the dynamics:
\begin{align}
    \begin{cases}
        \rho(t_k^-) = \frac{\partial}{\partial b}g(b(t_k^-),y_k)^\mathrm{T} \rho(t_k)\\
        \dot{\rho}(t) = -\frac{\partial}{\partial b}c(b(t),u(t)) - \frac{\partial}{\partial b} f(b(t),u(t))^\mathrm{T}\rho(t)
    \end{cases}
    \label{eq: continuous discrete adjoint dynamics}
\end{align}
with the boundary condition $\rho(t_f) = \frac{\partial}{\partial b}h(b(t_f))$.
Under the control-affine assumption for $f$ and the additive quadratic control cost,
the expected first order cost variation \eqref{eq: general belief linear variational cost} yields
\begin{multline}
\mathbb{E}[\nu(t_f)] = \frac{1}{2}v^\mathrm{T}C_u v + \mathbb{E}[\rho(\tau)]^\mathrm{T}H(b(\tau))(v - u(\tau)) \\ - \frac{1}{2}u(\tau)^\mathrm{T}C_u u(\tau),
\label{eq: expected control-affine belief linear variational cost}
\end{multline}
where $H(b(\tau))$ is the control coefficient term in $f$. 

Although it is difficult to state the general conditions under which this control-affine assumption holds, for instance one can verify that the continuous-discrete EKF \cite{xie2007estimation} satisfies this property if the underlying system dynamics $f_{\mathrm{sys}}$ is control-affine.
\begin{align}
    \begin{cases}
        \dot{\mu}(t) = f_{\mathrm{sys}}(\mu(t),u(t)) \\
        \dot{\Sigma}(t) = A\Sigma + \Sigma A^\mathrm{T} + Q
    \end{cases}
\end{align}
In the above continuous-time prediction equations, $A$ is the Jacobian of the dynamics function $f_{\mathrm{sys}}(x(t),u(t))$ evaluated at the mean $\mu(t)$ and $Q$ is the process noise covariance. If $f_{\mathrm{sys}}$ is control-affine, so is $A$ and therefore so is $\dot{\Sigma}$. Obviously $\dot{\mu}$ is control affine as well. As a result the dynamics for the belief vector $b = (\mu^\mathrm{T}, \text{vec}(\Sigma)^\mathrm{T})^\mathrm{T}$ satisfy the control-affine assumption.

Mirroring the approach in Section \ref{subsec: mixed observability}, we can use Monte Carlo sampling to estimate the expected value in \eqref{eq: expected control-affine belief linear variational cost}. The resulting algorithm is presented in Algorithm \ref{algo: sacbp_algo_2}.
\begin{algorithm}[t]
	\caption{SACBP Control Update for General Belief Space Planning Problems}\label{algo: sacbp_algo_2}
	\begin{algorithmic}[1]
		\INPUT Current belief state $b_0 = b(t_0)$, nominal control schedule $u(t)$ for $t \in [t_0, t_f]$, perturbation duration $\epsilon$
		\OUTPUT Optimally perturbed control schedule $u^{\epsilon}(t)$ for $t \in [t_0, t_f]$.
		\For{$i$ = 1:$N$}
		\State Forward-simulate nominal belief state trajectory \eqref{eq: continuous discrete belDynamics} and sample observation sequence $(y^{i}_{1},\dots,y^{i}_{T})$ along the belief trajectory.
		\State Backward-simulate nominal adjoint trajectory $\rho^{i}$ \eqref{eq: continuous discrete adjoint dynamics} with sampled observations.
		\EndFor
		\State Compute Monte Carlo estimate: $\mathbb{E}[\rho] \approx \frac{1}{N}\sum_{i=1}^N \rho^{i}$.
		\For{($\tau = t_0 + t_\mathrm{calc} + \epsilon;~\tau \leq t_0 + t_\mathrm{calc} + \Delta t_o;~\tau \leftarrow \tau + \Delta t_c$)}
		\State Solve quadratic minimization \eqref{eq: convex program} with \eqref{eq: expected control-affine belief linear variational cost}. Store optimal value $\nu^*(\tau)$ and optimizer $v^*(\tau)$. 
		\EndFor
		\State $\tau^* \leftarrow \arg\min \nu^*(\tau)$, $v^* \leftarrow v^*(\tau^*)$
		\State $u^{\epsilon} \leftarrow PerturbControlTrajectory(u,v^*,\tau^*,\epsilon)$ \eqref{eq: perturbed control trajectory}
		\State \textbf{return} $u^{\epsilon}$
	\end{algorithmic}
\end{algorithm}

\subsection{Closed-loop Nominal Policy}
\label{subsec: closed-loop policy}
In Sections \ref{subsec: mixed observability} and \ref{subsec: general belief space planning} we assumed that the nominal control was an open-loop control schedule. However, one can think of a scenario where a nominal control is a closed-loop policy computed offline, such as a discrete POMDP policy that maps beliefs to actions \cite{kurniawati2008sarsop}. Indeed, SACBP can also handle closed-loop nominal policies. Let $\pi$ be a closed-loop nominal policy, which is a mapping from either an augmented state $s(t)$ or a belief state $b(t)$ to a control value $u(t)$. Due to the stochastic belief dynamics, the control values returned by $\pi$ in the future is also stochastic for $t \geq t_1$. This is reflected when we forward-propagate the nominal dynamics. Specifically, each sampled trajectory has a different control trajectory in addition to a different observation sequence. However, the equations are still convex quadratic in $v$ as shown below.
For problems with mixed observability, we have
\begin{multline}
\mathbb{E}[\nu(t_f)] = \frac{1}{2}v^\mathrm{T}C_u v + \mathbb{E}[\rho_p(\tau)]^\mathrm{T}H(p(\tau))\left\{v - \pi(s(\tau))\right\} \\- \frac{1}{2} \pi(s(\tau))^\mathrm{T} C_u \pi(s(\tau)). \label{eq: expected closed-loop linear variational cost}
\end{multline}
The general belief space planning case also yields a similar equation:
\begin{multline}
\mathbb{E}[\nu(t_f)] = \frac{1}{2}v^\mathrm{T}C_u v + \mathbb{E}[\rho(\tau)]^\mathrm{T}H(b(\tau))\left\{v - \pi(b(\tau))\right\} \\- \frac{1}{2} \pi(b(\tau))^\mathrm{T} C_u \pi(b(\tau)). \label{eq: expected belief closed-loop linear variational cost}
\end{multline}
Note that $s(\tau)$ and $b(\tau)$ are both deterministic since the first observation $y_1$ is not yet taken at $\tau < t_1$.
The expectations in \eqref{eq: expected closed-loop linear variational cost} and \eqref{eq: expected belief closed-loop linear variational cost} can be estimated using Monte Carlo sampling. The forward-simulation of the nominal trajectory in Line 2 of Algorithms \ref{algo: sacbp_algo_1} and \ref{algo: sacbp_algo_2} is now with the closed loop policy $\pi$, and the equations in Line 7 need to be replaced with \eqref{eq: expected closed-loop linear variational cost} and \eqref{eq: expected belief closed-loop linear variational cost}, respectively. However, the rest remains unchanged.

\subsection{Computation Time Complexity}
\label{subsec: runtime analysis}
Let us analyze the time complexity of the SACBP algorithm. The bottleneck of the computation is when the forward-backward simulation is performed multiple times (lines 1--4 of Algorithms \ref{algo: sacbp_algo_1} and \ref{algo: sacbp_algo_2}). The asymptotic complexity of this part is given by $O(N(\frac{t_f - t_0}{\Delta t_o})(M_{\text{forward}}+M_{\text{backward}}))$, where $M_{\text{forward}}$ and $M_{\text{backward}}$ are the times to respectively integrate the forward and backward dynamics between two successive observations. For a more concrete analysis let us use the Gaussian belief dynamics given by EKF as an example. For simplicity we assume the same dimension $n$ for the state, the control, and the observation. The belief state has dimension $O(n^2)$. Using the Euler scheme, the forward integration takes $M_\text{forward} = O((\frac{\Delta t_o}{\Delta t_c} + 1) n^3)$ since evaluating continuous and discrete EKF equations are both $O(n^3)$. Computation of the continuous part of the adjoint dynamics \eqref{eq: continuous discrete adjoint dynamics} is dominated by the evaluation of the Jacobian $\frac{\partial f}{\partial b}$, which is $O(n^5)$ because $O(n^3)$ operations to evaluate $f$ are carried out $O(n^2)$ times. The discrete part is also $O(n^5)$. Therefore, $M_\text{backward} = O((\frac{\Delta t_o}{\Delta t_c} + 1) n^5)$. Overall, the time complexity is $O(N(\frac{t_f - t_0}{\Delta t_o})(\frac{\Delta t_o}{\Delta t_c} + 1)n^5)$. This is asymptotically smaller in $n$ than belief iLQG, which is $O(n^6)$. See \cite{rafieisakhaei2017tlqg} for a comparison of time complexity among different belief space planning algorithms. We also remind the readers that SACBP is an online method and a naive implementation already achieves near real-time performance, computing control over a 2[s] horizon in about 0.1[s] to 0.4[s]. By near real-time we mean that a naive implementation of SACBP requires approximately $0.7\times t_\text{calc}$ to $3\times t_\text{calc}$ time to compute an action that must be applied $t_\text{calc}$[s] in the future. We expect that parallelization in a GPU and a more efficient implementation will result in real-time computation for SACBP.

\section{Analysis of Mode Insertion Gradient for Stochastic Hybrid Systems}
\label{sec: analysis}
The SACBP algorithm presented in Section \ref{sec: SACBP} as well as the original SAC algorithm \cite{ansari2016sac} both rely on the local sensitivity analysis of the cost functional with respect to the control perturbation. This first-order sensitivity term (i.e. $\nu(t_f)$ in our notation) is known as the mode insertion gradient in the mode scheduling literature \cite{egerstedt2006transition, Wardi2012switched}. In \cite{ansari2016sac} the notion of the mode insertion gradient has been generalized to handle a broader class of hybrid systems than discussed before. What remains to be seen is a further generalization of the mode insertion gradient to stochastic hybrid systems, such as the belief dynamics discussed in this paper. Indeed, the quantity we can optimize in \eqref{eq: convex program} is essentially the expected value of the first-order sensitivity of the total cost. This is not to be confused with the first-order sensitivity of the expected total cost, which would be the natural generalization of the mode insertion gradient to stochastic systems. In general, those two quantities can be different, since the order of expectation and differentiation may not be swapped arbitrarily. In this section, we provide a set of sufficient conditions under which the order can be exchanged. By doing so we show that 1) the notion of mode insertion gradient can be generalized to stochastic hybrid systems, and 2) the SACBP algorithm optimizes this generalized mode insertion gradient. Through this analysis we will see that the SACBP algorithm has a guarantee that, in expectation it performs at least as good as the nominal policy with an appropriate choice of $\epsilon$.

\subsection{Assumptions}
Let us begin with a set of underlying assumptions for the system dynamics, the control, and the cost functions. Without loss of generality, we assume that the system starts at time $t = 0$ and ends at $t = T$, with a sequence of $T$ observations $(y_1, \dots, y_T)$ made every unit time. For generality, we use notation $x$ to represent the state variable of the system in this section, in place of $b$ or $s$ that respectively represented the belief state or the augmented state in Section \ref{sec: SACBP}. This means that the analysis presented here is not restricted to belief systems where the dynamics are governed by Bayesian filters, but rather applies to a broader class of stochastic systems.
\begin{assum}[Control Model]
    \label{assum: control}
    The controls are in $\widetilde{C}^{0,m}[0, T]$, the space of piecewise continuous functions from $[0, T]$ into $\mathbb{R}^m$. We further assume that there exists some $\rho_{\max} < \infty$ such that for all $t \in [0, T]$, we have $u(t) \in B(0, \rho_{\max})$ where $B(0, \rho_{\max})$ is the closed Euclidean ball of radius $\rho_{\max}$ centered at $0$, i.e. $\Vert u(t) \Vert_2 \leq \rho_{\max}$. We denote this admissible control set by $U \triangleq \{u \in \widetilde{C}^{0,m}[0, T] \mid \forall t \in [0, T]\; u(t) \in B(0, \rho_{\max})\}$.
\end{assum}

\begin{remark}
    For the sake of analysis, the control model described above takes the form of an open-loop control schedule, where time $t$ determines the control signal $u(t)$. However the analysis can be extended to closed-loop nominal policies in a straightforward manner, which is discussed in Appendix \ref{sec: appendix}. (See Remark \ref{remark: closed_loop_policy}.)
\end{remark}

\begin{assum}[Dynamics Model]
    \label{assum: dynamics}
    Let $x_0 \in \mathbb{R}^{n_x}$ be the given initial state value at $t = 0$. Given a control $u \in U$ and a sequence of observations $(y_1, \dots, y_T) \in \mathbb{R}^{n_y}\times\cdots\times\mathbb{R}^{n_y}$, the dynamics model is the following hybrid system with time-driven switching:
    \begin{align}
        x(t) \triangleq x_i(t)\;  \forall t \in [i-1, i)\; \forall i \in \{1,2,\dots, T\},
    \end{align}
    where $x_i$ is the $i$-th ``mode" of the system state defined on $[i-1, i]$ as:
    \begin{align}
        x_{i}(i-1) &= g\left(x_{i-1}(i-1), y_{i-1}\right) \label{eq: x_jump_definition}\\
        \dot{x_i}(t) &= f\left(x_{i}(t), u(t)\right) &\forall t \in [i-1, i], \label{eq: xdot_definition}
    \end{align}
    with $x(0) = x_1(0) = x_0$. We also define the final state as $x(T) \triangleq g(x_{T}(T), y_T)$.
    
    For the transition functions $f$ and $g$ we assume the following:
    \begin{enumerate}[label=(\arabic{assum}\alph*)]
        \item \label{assum: function_f} the function $f\colon \mathbb{R}^{n_x} \times \mathbb{R}^m \rightarrow \mathbb{R}^{n_x}$ is continuously differentiable;
        \item \label{assum: function_g} the function $g\colon \mathbb{R}^{n_x} \times \mathbb{R}^{n_y} \rightarrow \mathbb{R}^{n_x}$ is continuous. It is also differentiable in $x$;
        \item \label{assum: function_f_bound} for function $f$, there exists a constant $K_1 \in [1, \infty)$ such that $\forall x^{\prime}, x^{\prime\prime} \in \mathbb{R}^{n_x}$ and $\forall u^{\prime}, u^{\prime\prime} \in B(0,\rho_{\max})$, the following condition holds:
        \begin{multline}
            \Vert f(x^{\prime}, u^{\prime}) - f(x^{\prime\prime}, u^{\prime\prime}) \Vert_2 \\ \leq K_1 \left(\Vert x^{\prime} - x^{\prime\prime} \Vert_2 + \Vert u^{\prime} - u^{\prime\prime}\Vert_2\right)
        \end{multline}

        \item \label{assum: function_g_bound} for function g, there exist finite non-negative constants $K_2, K_3, K_4, K_5$ and positive integers $L_1, L_2$ such that $\forall x \in \mathbb{R}^{n_x}$ and $\forall y \in \mathbb{R}^{n_y}$, the following relations hold:
        \begin{multline}
            \label{eq: g_norm_bound}
            \Vert g(x, y) \Vert_2 \leq K_2 + K_3\Vert x \Vert_2^{L_1} + K_4\Vert y \Vert_2^{L_2} \\
            + K_5\Vert x \Vert_2^{L_1} \Vert y \Vert_2^{L_2} 
        \end{multline}
        \begin{multline}
            \label{eq: dg_norm_bound}
            \left\Vert \frac{\partial}{\partial x}g(x, y) \right\Vert_2 \leq K_2 + K_3\Vert x \Vert_2^{L_1} + K_4\Vert y \Vert_2^{L_2} \\
            + K_5\Vert x \Vert_2^{L_1} \Vert y \Vert_2^{L_2} 
        \end{multline}
    \end{enumerate}
\end{assum}

\begin{remark}
    Assumptions \ref{assum: function_f} and \ref{assum: function_f_bound} are used to show existence and uniqueness of the solution to the differential equation \eqref{eq: xdot_definition} as well as the variational equation under control perturbation. (See Propositions \ref{prop: existence} and \ref{prop: state_variation} in Appendix \ref{sec: appendix}.) Note that Assumption \ref{assum: function_f_bound} is essentially a Lipschitz continuity condition, which is a quite common assumption in the analysis of nonlinear ODEs, as can be seen in Assumption 5.6.2 in \cite{elijah1997optimization}  and Theorem 2.3 in \cite{khalil2002nonlinear}. Assumptions \ref{assum: function_g} and \ref{assum: function_g_bound} are the growth conditions on $x$ across adjacent modes. Recall that in belief space planning where the system state $x$ is the belief state $b$, the jump function $g$ corresponds to the observation update of the Bayesian filter. The form of the bound in \eqref{eq: g_norm_bound} and \eqref{eq: dg_norm_bound} allows a broad class of continuous functions to be considered as $g$, and is inspired by a few examples of the Bayesian update equations as presented below.
    \begin{prop}[Bounded Jump for Univariate Gaussian Distribution] 
    Let $b = (\mu, s)^\mathrm{T} \in \mathbb{R}^2$ be the belief state, where $\mu$ is the mean parameter and $s > 0$ is the variance. Suppose that the observation $y$ is the underlying state $x \in \mathbb{R}$ corrupted by additive Gaussian white noise $v \sim \mathcal{N}(0, 1)$. Then, the Bayesian update function $g$ for this belief system satisfies Assumption \ref{assum: function_g_bound}.
    \begin{proof}
    The Bayesian update formula for this system is given by $g(b,y) = \hat{b} \triangleq (\hat{\mu}, \hat{s})^\mathrm{T}$, where
        \begin{align}
            \hat{\mu} &= \mu + \frac{s}{s + 1}(y - \mu) \\
            \hat{s} &= s - \frac{s^2}{s + 1}
        \end{align}
        is the update step of the Kalman filter. Rearranging the terms, we have
        \begin{align}
            g(b,y) = \frac{1}{s + 1}
            \begin{pmatrix} 
                \mu + sy \\
                s          
            \end{pmatrix}
        \end{align}
        and consequently,
        \begin{align}
            \frac{\partial}{\partial b}g(b, y) = \frac{1}{(s + 1)^2}
            \begin{pmatrix}
                s + 1 & y \\
                0 & 1
            \end{pmatrix}.
        \end{align}
        We will show that the function $g$ satisfies Assumption \ref{assum: function_g_bound}. For the bound on $g(b,y)$,
        \begin{align}
            \Vert g(b,y) \Vert_2^2 &= \frac{1}{(s + 1)^2} \left\{(\mu + sy)^2 + s^2\right\} \\
            &\leq (\mu + sy)^2 + s^2 \\
            &\leq \Vert b \Vert_2^2 + \left(b^\mathrm{T} \begin{pmatrix} y \\ 1 \end{pmatrix}\right)^2 \\
            &\leq \Vert b \Vert_2^2(2 + \Vert y \Vert_2^2),
        \end{align}
        where we have used $(s + 1)^2 \geq 1$ and the Cauchy-Schwarz inequality. Thus,
        \begin{align}
            \Vert g(b,y) \Vert_2 &\leq \Vert b \Vert_2 \sqrt{2 + \Vert y \Vert_2^2} \\
            &\leq \sqrt{2} \Vert b \Vert_2 + \Vert b \Vert_2 \Vert y \Vert_2
        \end{align}
        Similarly, the bound on the Jacobian yields
        \begin{align}
            \left\Vert \frac{\partial}{\partial b} g(b,y) \right\Vert_2^2 &\leq \left\Vert \frac{\partial}{\partial b} g(b,y) \right\Vert_\mathrm{F}^2 \\
            &= \frac{1}{(s + 1)^4}\left\{(s + 1)^2 + y^2 + 1\right\} \\
            &= \frac{1}{(s + 1)^2} + \frac{1}{(s + 1)^4}(y^2 + 1) \\
            &\leq 2 + \Vert y \Vert_2^2.
        \end{align}
        Therefore, 
        \begin{align}
            \left\Vert \frac{\partial}{\partial b} g(b,y) \right\Vert_2 \leq \sqrt{2} + \Vert y \Vert_2
        \end{align}
        This shows that the jump function $g$ for the above univariate Gaussian model satisfies Assumption \ref{assum: function_g_bound}  with $(K_2, K_3, K_4, K_5) = (\sqrt{2}, \sqrt{2}, 1, 1)$ and $(L_1, L_2) = (1, 1).$
        \end{proof}
    \end{prop}
    \begin{prop}[Bounded Jump for Categorical Distribution]
        Let $b = (b_1, \dots, b_n)^\mathrm{T} \in \mathbb{R}^n$ be the $n$-dimensional belief state representing the categorical distribution over the underlying state $x \in \{1, \dots, n\}$. We choose the unnormalized form where the probability of $x = i$ is given by $b_i/\sum_{i=1}^n b_i$. Let the observation $y \in \{1, \dots, m\}$ be modeled by a conditional probability mass function $p(y \mid x) \in [0, 1]$. Then, the Bayesian update function $g$ for this belief system satisfies Assumption \ref{assum: function_g_bound}.
    \begin{proof}
        The Bayes rule gives $g(b,y) = \hat{b} \triangleq (\hat{b}_1, \dots, \hat{b}_n)$, where
        \begin{align}
            \begin{pmatrix} \hat{b}_1 \\ \hat{b}_2 \\ \vdots \\ \hat{b}_n \end{pmatrix}
            =
            \begin{pmatrix} p(y \mid 1) b_1 \\  p(y \mid 2) b_2 \\ \vdots \\ p(y \mid n) b_n
            \end{pmatrix}.
        \end{align}
        Therefore, we can easily bound the norm of the posterior belief $\hat{b}$ by
        \begin{align}
            \Vert g(b, y) \Vert_2 = \Vert \hat{b} \Vert_2 \leq \Vert b \Vert_2,
        \end{align}
        as $p(y \mid x) \leq 1$. The Jacobian is simply the diagonal matrix $\mathrm{diag}(p(y \mid 1), \dots, p(y \mid n))$, and hence
        \begin{align}
            \left\Vert \frac{\partial}{\partial b} g(b,y) \right\Vert_2 \leq 1.
        \end{align}
        This shows that the jump function $g$ for the categorical belief model above satisfies Assumption \ref{assum: function_g_bound} with $(K_2, K_3, K_4, K_5) = (1, 1, 0, 0)$ and $(L_1, L_2) = (1, 1).$
    \end{proof}
    \end{prop}
\end{remark}

\begin{assum}[Cost Model]
    \label{assum: cost}
    The instantaneous cost $c\colon \mathbb{R}^{n_x} \times \mathbb{R}^m \rightarrow \mathbb{R}$ is continuous. It is also continuously differentiable in $x$. The terminal cost $h\colon \mathbb{R}^{n_x} \rightarrow \mathbb{R}$ is differentiable. 
    Furthermore, we assume that there exist finite non-negative constants $K_6, K_7$ and a positive integer $L_3$ such that for all $x \in \mathbb{R}^{n_x}$ and $u \in B(0, \rho_{\max})$, the following relations hold:
    \begin{align}
        \left|c(x, u)\right| &\leq K_6 + K_7 \Vert x \Vert_2^{L_3} \\
        \left\Vert \frac{\partial}{\partial x} c(x, u) \right\Vert_2 &\leq K_6 + K_7 \Vert x \Vert_2^{L_3} \\
        \left| h(x) \right| &\leq K_6 + K_7 \Vert x \Vert_2^{L_3} \\
        \left\Vert \frac{\partial}{\partial x} h(x) \right\Vert_2 &\leq K_6 + K_7 \Vert x \Vert_2^{L_3}. 
    \end{align}
\end{assum}

\begin{remark}
    Assumption \ref{assum: cost} is to guarantee that the cost function is integrable with respect to stochastic observations, which are introduced in Assumption \ref{assum: observations}. Note that even though the above bound is not general enough to apply to all analytic functions, it does include all finite order polynomials of $\Vert x(t) \Vert_2$ and $\Vert u(t) \Vert_2$, for example, since $\Vert u(t) \Vert_2$ is bounded by Assumption \ref{assum: control}.
\end{remark}

\begin{assum} [Stochastic Observations]
    \label{assum: observations}
    Let $(\Omega, \mathcal{F}, \mathbb{P})$ be a probability space. Let $(Y_1, \dots, Y_T)$ be a sequence of random vectors in $\mathbb{R}^{n_y}$ defined on this space, representing the sequence of observations. Assume that for each $Y_i$ all the moments of the $\ell^2$ norm is finite. That is,
    \begin{align}
        \forall i \in \{1, \dots, T\}\; \forall k \in \mathbb{N}\;\; \mathbb{E}\left[\Vert Y_i \Vert_2^k\right] < \infty.
    \end{align}
\end{assum}

\begin{defn}[Perturbed Control]
    Let $u \in U$ be a control. For $\tau \in (0, 1)$ and $v \in B(0, \rho_{\max})$, define the perturbed control $u^{\epsilon}$ by
    \begin{align}
        \label{eq: perturbed_control_main_body}
        u^{\epsilon}(t) \triangleq \begin{cases}
            v & \text{if}\;\; t \in (\tau - \epsilon, \tau] \\
            u(t) & \text{otherwise},
        \end{cases}
    \end{align}
    where $\epsilon \in [0, \tau]$. By definition if $\epsilon = 0$ then $u^\epsilon$ is the same as $u$. We assume that the nominal control $u(t)$ is left continuous in $t$ at $t = \tau$. 
\end{defn}
\setcounter{defn}{0}

\subsection{Main Results}
The main result of the analysis is the following theorem.
\begin{theorem}[Mode Insertion Gradient]
    \label{theorem: mode_insertion_gradient}
     Suppose that Assumptions \ref{assum: control} -- \ref{assum: observations} are satisfied. For a given $(\tau, v)$, let $u^{\epsilon}$ denote the perturbed control of the form \eqref{eq: perturbed_control_main_body}. The perturbed control $u^{\epsilon}$ and the stochastic observations $(Y_1, \dots, Y_T)$ result in the stochastic perturbed state trajectory $x^{\epsilon}$. For such $u^{\epsilon}$ and $x^{\epsilon}$, let us define the mode insertion gradient of the expected total cost as
    \begin{align}
        \label{eq: stochastic_mode_insertion_gradient}
        \frac{\partial_+}{\partial\epsilon} \mathbb{E}\left[\int_{0}^{T} c(x^{\epsilon}(t), u^{\epsilon}(t)) dt + h(x^{\epsilon}(T))\right]\Bigg\vert_{\epsilon=0}.
    \end{align}
    Then, this right derivative exists and we have
    \begin{multline}
        \label{eq: stochastic_mode_insertion_gradient_equality}
         \frac{\partial_+}{\partial\epsilon} \mathbb{E}\left[\int_{0}^{T} c(x^{\epsilon}(t), u^{\epsilon}(t)) dt + h(x^{\epsilon}(T))\right]\Bigg\vert_{\epsilon=0} \\
          = c(x(\tau), v) - c(x(\tau), u(\tau)) \\ +
          \mathbb{E}\Bigg[\int_{\tau}^{T} \frac{\partial}{\partial x} c(x(t), u(t))^\mathrm{T} \Psi(t) dt   \\ + \frac{\partial}{\partial x} h(x(T))^\mathrm{T}\Psi(T)\Bigg],
    \end{multline}
    where $\Psi(t) = \frac{\partial_+}{\partial\epsilon} x^{\epsilon}(t) \big\vert_{\epsilon=0}$ is the state variation.
\end{theorem}
\setcounter{theorem}{0}

The proof of the theorem is deferred to Appendix \ref{sec: appendix}. One can see that the mode insertion gradient \eqref{eq: stochastic_mode_insertion_gradient} is a natural generalization of the ones discussed in \cite{egerstedt2006transition, Wardi2012switched, ansari2016sac} to stochastic hybrid systems. Furthermore, by comparing \eqref{eq: stochastic_mode_insertion_gradient_equality} with \eqref{eq: total_cost_variation_forward_integration} it is apparent that the right hand side of \eqref{eq: stochastic_mode_insertion_gradient_equality} is mathematically equivalent to $\mathbb{E}[\nu(t_f)]$, the quantity to be optimized with the SACBP algorithm in Section \ref{sec: SACBP}.

The fact that SACBP optimizes \eqref{eq: stochastic_mode_insertion_gradient} leads to a certain performance guarantee of the algorithm. In the open-loop nominal control case, the term $\mathbb{E}[\nu(t_f)]$ as in \eqref{eq: expected linear variational cost} or \eqref{eq: expected control-affine belief linear variational cost} becomes $0$ if the control perturbation $v$ is equal to the nominal control $u(\tau)$. Therefore, as long as $u(\tau)$ is a feasible solution to \eqref{eq: convex program} the optimal value is guaranteed to be less than or equal to zero. Furthermore, in expectation the actual value of $\mathbb{E}[\nu(t_f)]$ matches the one approximated with samples, since the Monte Carlo estimate is unbiased. In other words, the perturbation $(\tau^*, v^*)$ computed by the algorithm is expected to result in a non-positive mode insertion gradient. If the mode insertion gradient is negative, there always exists a sufficiently small $\epsilon>0$ such that the expected total cost is decreased by the control perturbation. In the corner case that the mode insertion gradient is zero, one can set $\epsilon=0$ to not perturb the control at all. Therefore, for an appropriate choice of $\epsilon$ the expected performance of the SACBP algorithm over the planning horizon is at least as good as that of the nominal control.

The same discussion holds for the case of closed-loop nominal control policies, when the expression for $\mathbb{E}[\nu(t_f)]$ is given by \eqref{eq: expected closed-loop linear variational cost} or \eqref{eq: expected belief closed-loop linear variational cost}. This is because Theorem 1 still holds if the nominal control $u(t)$ is a closed-loop policy as stated in Remark \ref{remark: closed_loop_policy} (Appendix \ref{sec: appendix}). Therefore, the expected worst-case performance of the algorithm is lower-bounded by that of the nominal policy. This implies that if a reasonable nominal policy is known, at run-time SACBP is expected to further improve it while synthesizing continuous-time control inputs efficiently.


\section{Simulation Results}
\label{sec: results}
We evaluated the performance of SACBP in the following two simulation studies: (i) active multi-target tracking with range-only observations; (ii) object manipulation under model uncertainty. SACBP as well as other baseline methods were implemented in Julia\footnote{The code for SACBP as well as the baseline methods is publicly available at \url{https://github.com/StanfordMSL/SACBP.jl}.}, except for T-LQG \cite{rafieisakhaei2017tlqg} whose NLP problems were modeled by CasADi \cite{andersson2019casadi} in Python and then solved by Ipopt \cite{wachter2006implementation}, a standard NLP solver based on interior-point methods.
All the computation was performed on a desktop computer with Intel\textregistered\ Core\texttrademark\ i7-8750H CPU and 32.1GB RAM. The Monte Carlo sampling of SACBP was parallelized on multiple cores of the CPU.
\subsection{Active Multi-Target Tracking with Range-only Observations}

\begin{figure*}[tpb]
\begin{center}
	\begin{tabular}{c}
		\begin{minipage}[t]{0.31\hsize}
			\centering
			\includegraphics[clip,width=1.0\columnwidth]{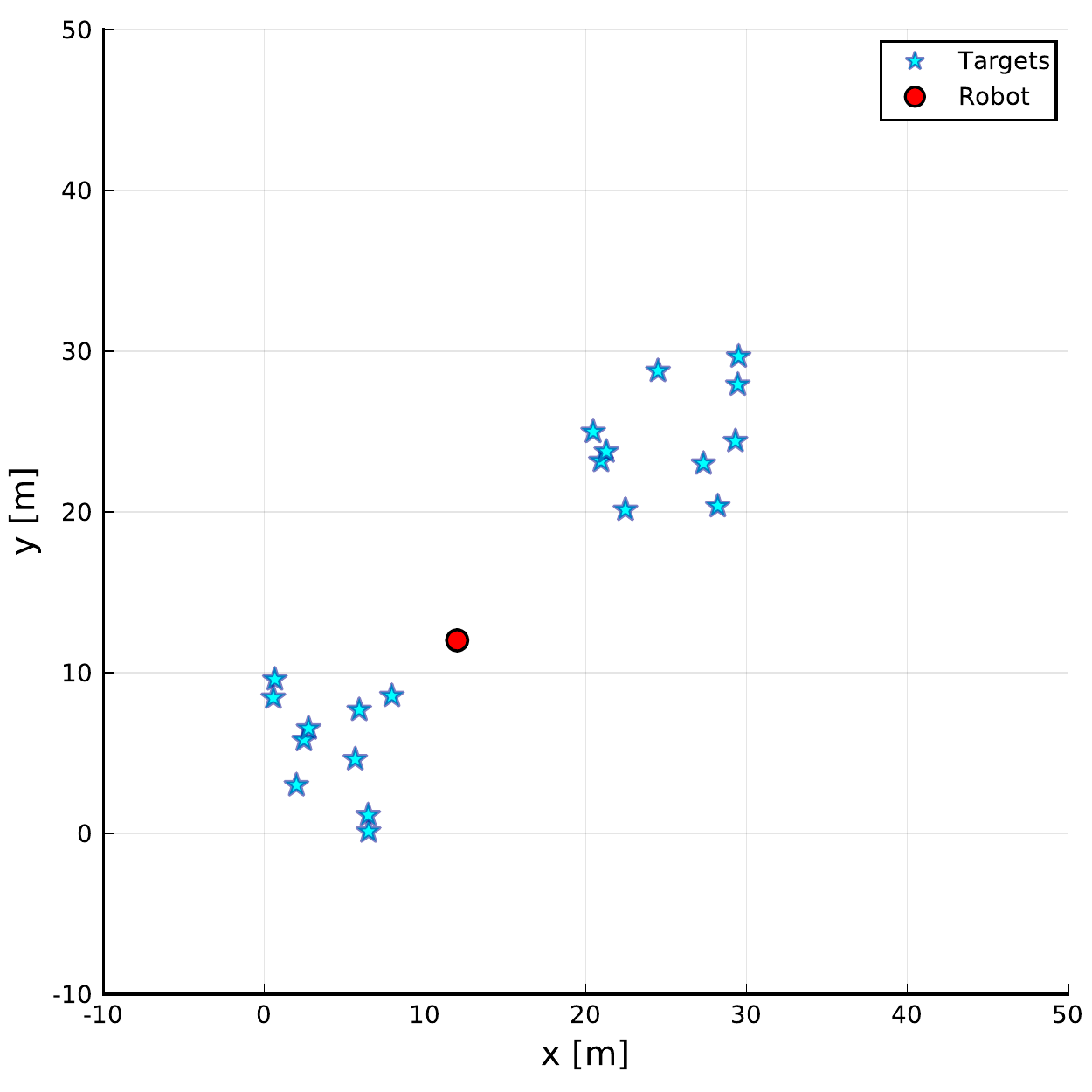}
		\end{minipage}
		\begin{minipage}[t]{0.31\hsize}
			\centering
			\includegraphics[clip,width=1.0\columnwidth]{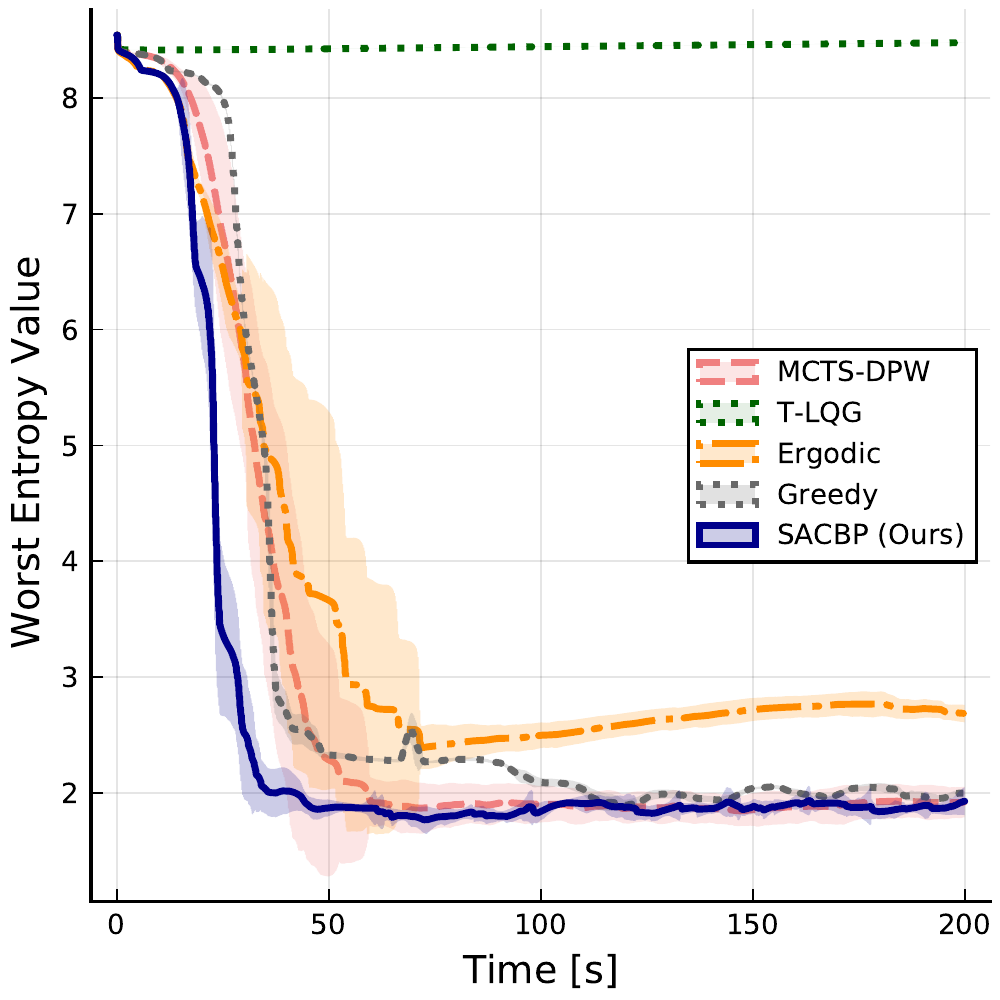}
		\end{minipage}
		\begin{minipage}[t]{0.37\hsize}
			\centering
			\includegraphics[clip,width=1.0\columnwidth]{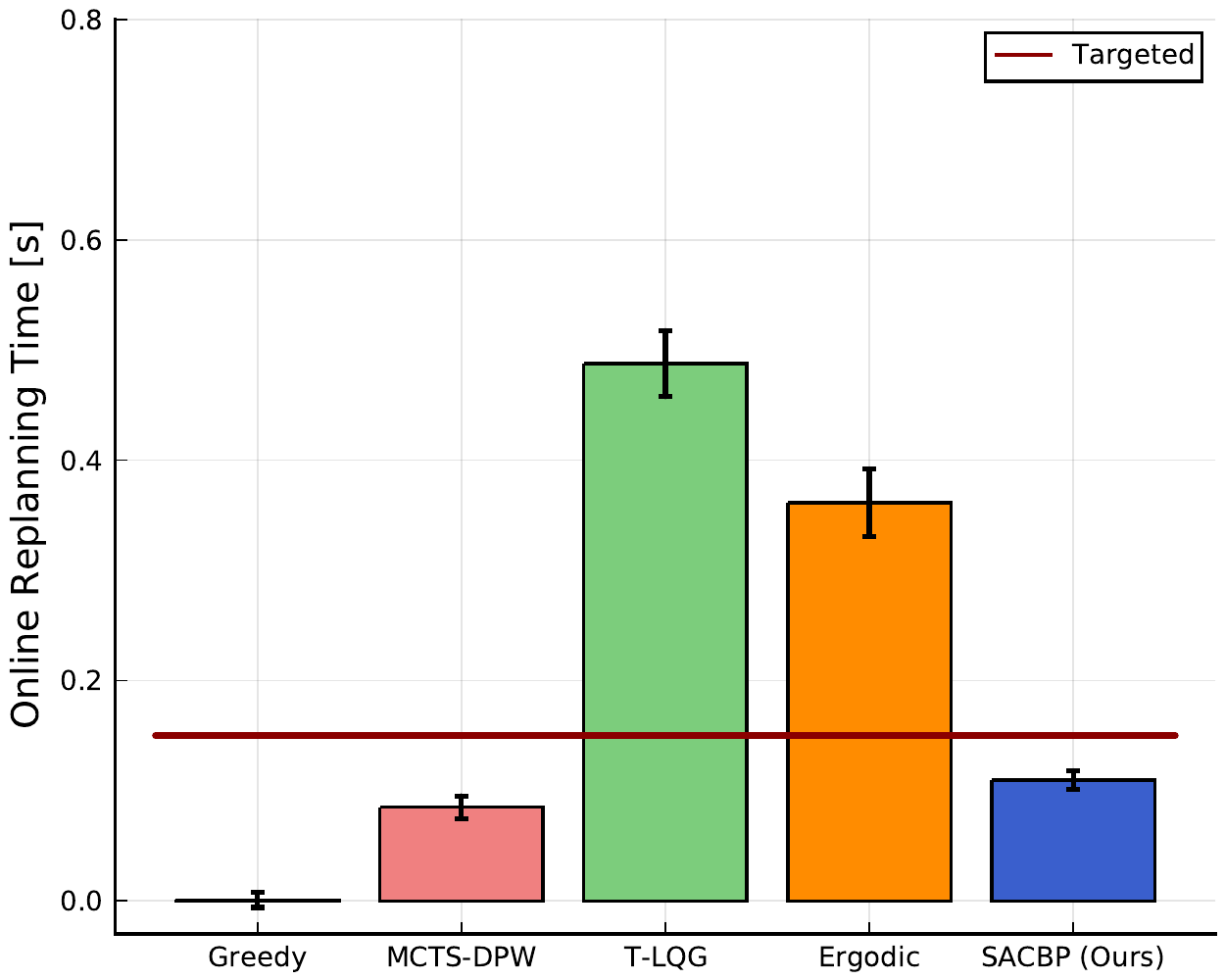}
		\end{minipage}
	\end{tabular}
	\caption{(Left) Simulation environment with 20 targets and a surveillance robot. (Middle) The history of the worst entropy value among the targets averaged over 20 random runs, plotted with the standard deviation. With the budget of 10 Monte Carlo samples, SACBP showed small variance for the performance curve and resulted in the fastest reduction of the worst entropy value compared to every other baseline. (Right) Computation times for Greedy, MCTS-DPW, and SACBP achieved real-time performance, taking less time than simulated $t_\text{calc} = 0.15[s]$.}
	\label{fig: tracking_1}
\end{center}
\end{figure*}
\begin{figure*}[ht]
	\begin{center}
		\begin{tabular}{c}
		    \vspace{1.0cm} \\
			\begin{minipage}[b]{0.32\hsize}
				\centering
				\includegraphics[clip,width=1.0\columnwidth]{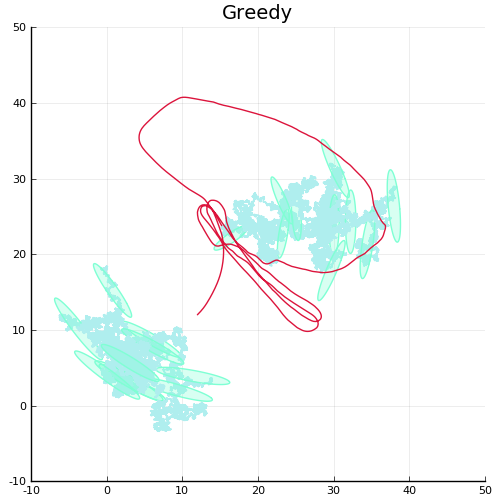}
			\end{minipage}
			\begin{minipage}[b]{0.32\hsize}
				\centering
				\includegraphics[clip,width=1.0\columnwidth]{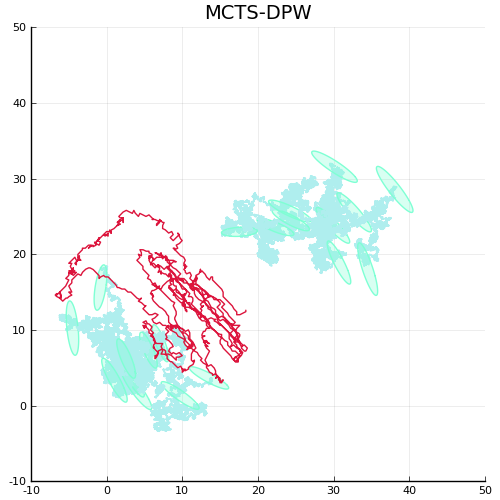}
			\end{minipage}
			\begin{minipage}[b]{0.32\hsize}
				\centering
				\includegraphics[clip,width=1.0\columnwidth]{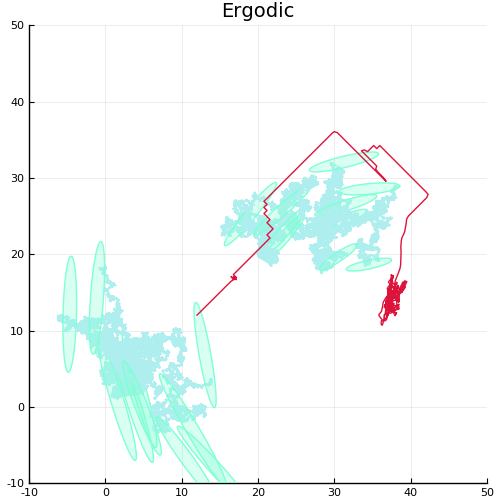}
			\end{minipage} \\
			\begin{minipage}[b]{0.32\hsize}
				\centering
				\includegraphics[clip,width=1.0\columnwidth]{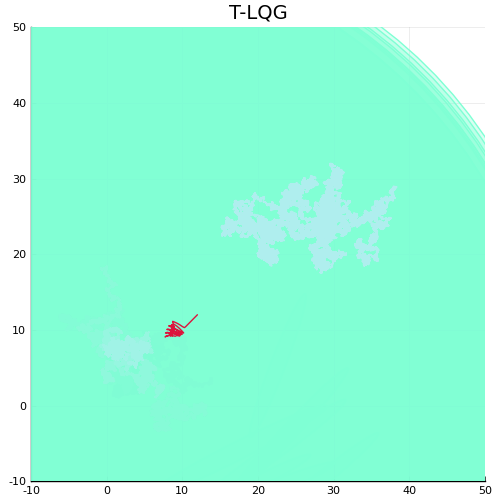}
			\end{minipage}
			\begin{minipage}[b]{0.32\hsize}
			    \includegraphics[clip,width=1.0\columnwidth]{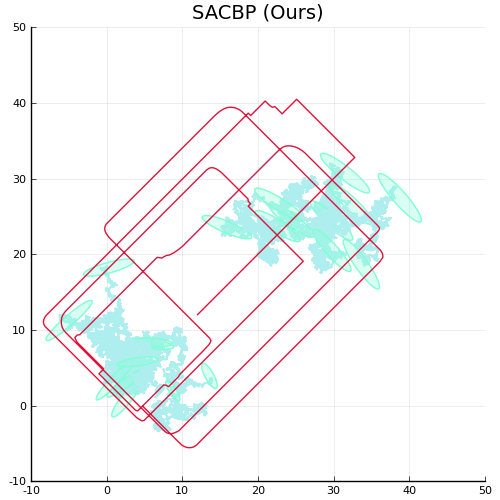}
			\end{minipage}
		\end{tabular}
	\end{center}
	\caption{Sample robot trajectories (depicted in red) generated by each algorithm. Greedy, MCTS-DPW, and Ergodic did not result in a trajectory that fully covers the two groups of the targets. T-LQG failed to reduce the estimation uncertainty even after 200[s], due to insufficient time to solve the NLP with high-dimensional joint states in an online manner. SACBP successfully explored the space and periodically revisited both of the two target groups. With SACBP, the robot traveled into one of the four diagonal directions for most of the time. This is due to the fact that SACBP optimizes a convex quadratic under a box saturation constraint, which tends to find optimal solutions at the corners. In all the figures, the blue lines represent the target trajectories and the shaded ellipses are 99\% error ellipses at $t = 200$[s].}
	\label{fig: tracking_2}
\end{figure*}

\begin{figure*}[t]
	\begin{center}
		\begin{tabular}{c}
			\begin{minipage}[t]{0.40\hsize}
				\centering
				\includegraphics[clip,width=1.0\columnwidth]{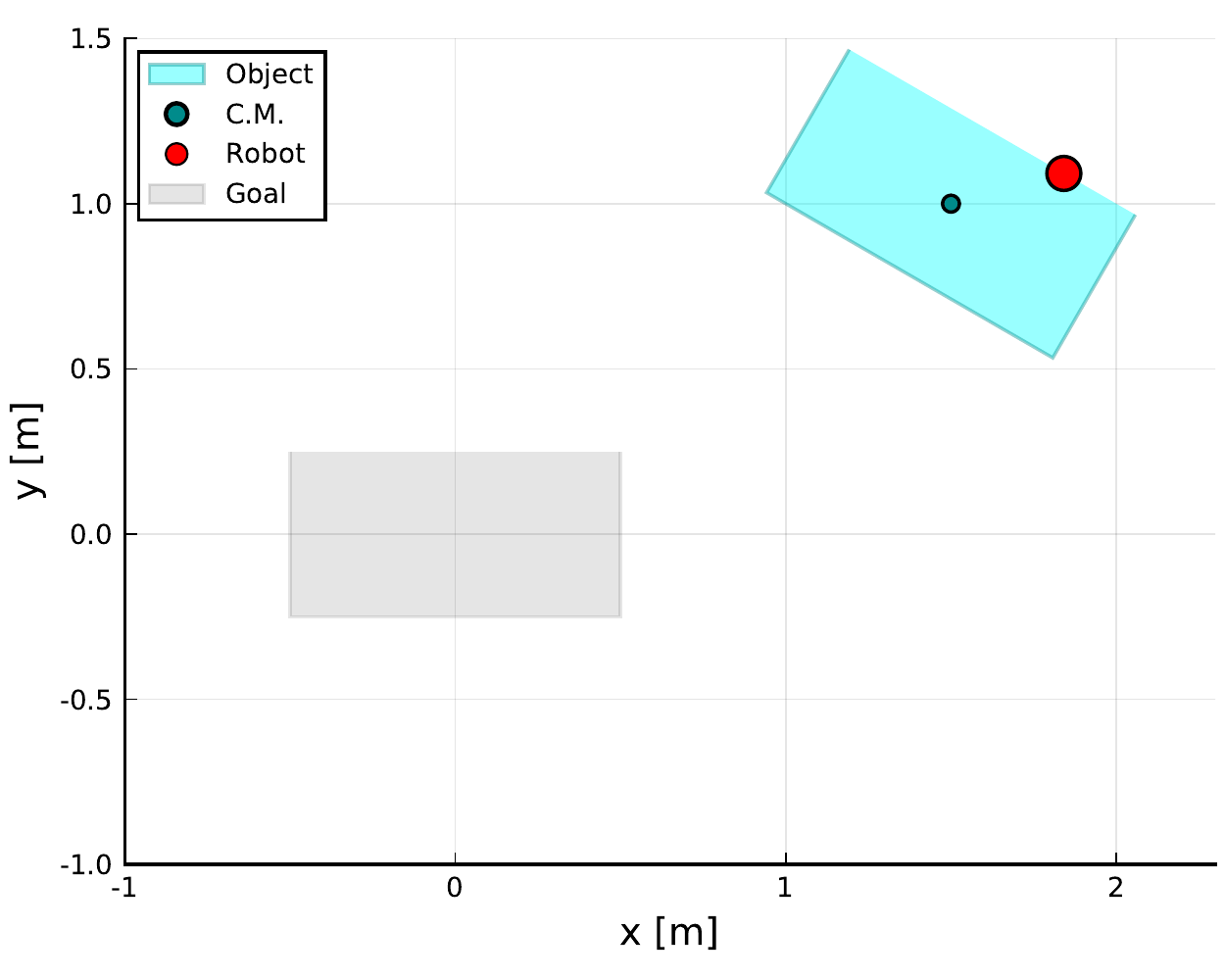}
			\end{minipage}
			\begin{minipage}[t]{0.33\hsize}
				\centering
				\includegraphics[clip,width=1.0\columnwidth]{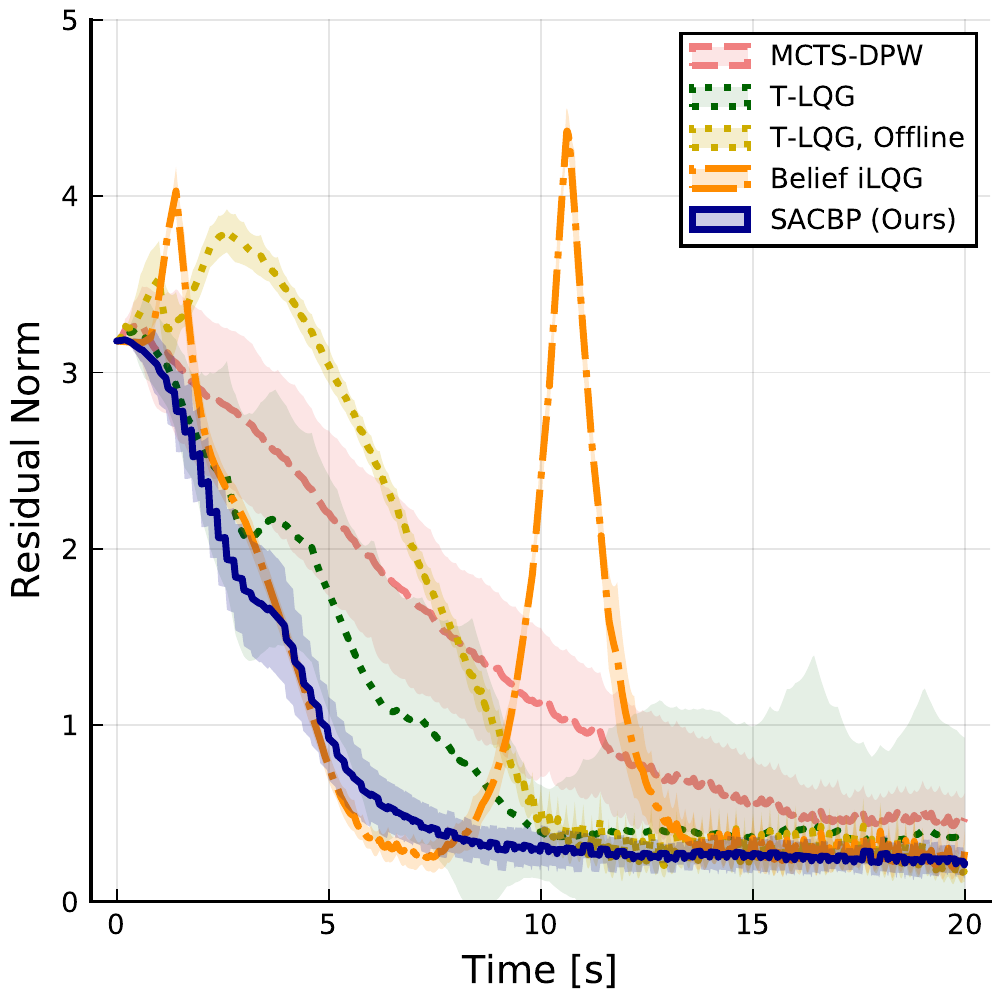}
			\end{minipage}
			\begin{minipage}[t]{0.24\hsize}
				\centering
				\includegraphics[clip,width=1.0\columnwidth]{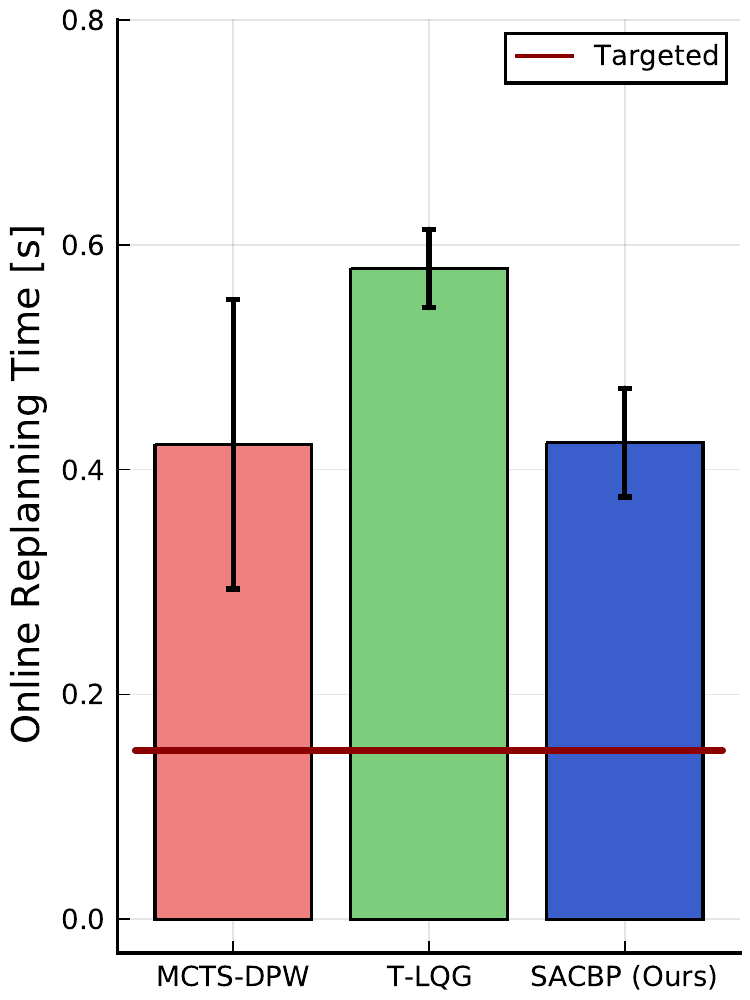}
			\end{minipage}
		\end{tabular}
		\caption{(Left) The robot is attached to the rectangular object. (Middle) The history of the $l_2$ norm of the residual between the goal state and the true object state averaged over 20 runs. SACBP with $N=10$ samples successfully brought the object close to the goal. The reduction of the residual norm was much slower for MCTS-DPW. T-LQG was not as successful either, regardless of whether the policy was derived offline (without re-planning) or online (with re-planning), although it eventually achieved similar residual norms to SACBP. Belief iLQG resulted in large overshoots at around 2[s] and 11[s]. (Right) Computation time of SACBP was increased from the multi-target tracking problem due to increased complexity related to the continuous-discrete belief dynamics, but still achieved a reasonable value. Note that the computation times of the offline algorithms were significantly longer and are not shown in this plot.}
		\label{fig: manipulation_1}
	\end{center}
\end{figure*}

This problem focuses on pure information gathering, namely identifying where the moving targets are in the environment. In doing so, the surveillance robot modeled as a single integrator can only use relative distance observations.
The robot's position $p$ is fully observable and the transitions are deterministic. Assuming perfect data association, the observation for target $i$ is $d_i = ||q_i - p + v_i||_2$, where $q_i$ is the true target position and $v_i$ is zero-mean Gaussian white noise with state-dependent covariance $R(p,q_i) = R_0 + ||q_i - p||_2R_1$. We use $0.01 I_{2\times2}$ for the nominal noise $R_0$. The range-dependent noise $R_1 = 0.001 I_{2\times2}$ degrades the observation quality as the robot gets farther from the target. The discrete-time UKF is employed for state estimation in tracking 20 independent targets. The target dynamics are modeled by a 2D Brownian motion with covariance $Q = 0.1I_{2\times2}$. Similarly to \cite{spinello2010statedepedent}, an approximated observation covariance $R(p,\mu_i)$ is used in the filter to obtain tractable estimation results, where $\mu_i$ is the most recent mean estimate of $q_i$.

The SACBP algorithm produces a continuous robot trajectory over 200[s] with planning horizon $t_f - t_0 = 2$[s], update interval $\Delta t_o = 0.2[s]$, perturbation duration $\epsilon = 0.16[s]$, and $N = 10$ Monte Carlo samples. The Euler scheme is used for integration with $\Delta t_c = 0.01[s]$. The Jacobians and the gradients are computed either analytically or using an automatic differentiation tool \cite{revels2016forwarddiff} to retain both speed and precision. In this simulation $t_\text{calc} = 0.15[s]$ is assumed no matter how long the actual control update takes. We use $c(p,b,u) = 0.05u^\mathrm{T}u$ for the running cost and $h(p,b)= \sum_{i=1}^{20}\exp(\text{entropy}(b_i))$ for the terminal cost, with an intention to reduce the worst-case uncertainty among the targets. This expression for $h(p,b)$ is equivalent to:
\begin{align}
    h(p,b) = \sum_{i=1}^{20} \sqrt{\text{det}(2\pi e \Sigma_i)},
\end{align}
where $\Sigma_i$ is the covariance for the $i$-th target.
The nominal control $u(t)$ is constantly zero.

We compared SACBP against four baseline methods: (i) a greedy algorithm based on the gradient descent of terminal cost $h$, similar to \cite{schwager2017info}; (ii) MCTS-DPW \cite{couetoux2011cuct, egorov2017pomdps} in the Gaussian belief space; (iii) projection-based trajectory optimization for ergodic exploration \cite{miller2013ergodic, miller2016ergodic, dressel2018ergoric}; (iv) T-LQG \cite{rafieisakhaei2017tlqg}.
We also attempted to implement the belief iLQG \cite{vandenberg2012beliefilqg} algorithm, but the policy did not converge for this problem. We suspect that the non-convex terminal cost $h$ contributed to this behavior, which in fact violates one of the underlying assumptions made in the paper \cite{vandenberg2012beliefilqg}.

MCTS-DPW uses the same planning horizon as SACBP, however it draws $N = 25$ samples from the belief tree so the computation time of the two algorithms matches approximately.

For T-LQG, an NLP is formulated so that the optimization objective is a discrete-time equivalent of the SACBP objective, with the discrete time interval of $\Delta t_o = 0.2[s]$.  
Unfortunately, the Ipopt solver takes a significantly long time to solve this NLP; on average the optimization over 10 timesteps (thus 2[s]) takes 43.9[s] to converge to a local minimum, with the worst case of over 250[s], due to the high dimensionality of the joint state space. This high computation cost makes it prohibitive to derive an offline policy over the entire 200[s] simulation horizon prior to the execution. Therefore, we only test T-LQG in an online manner for this problem, with the same planning horizon as SACBP. For a fair comparison with the other methods, we adjust the $\textrm{max\_cpu\_time}$ parameter of the solver so that the Ipopt iterations terminate within a computation time that is the same order of magnitude as SACBP, MCTS-DPW, and ergodic exploration.
Note also that T-LQG is a closed-loop planning method that comes with a local LQG controller to track the optimized nominal belief trajectory. Designed to require minimal online computation, re-planning for T-LQG only happens if the end of the planning horizon is reached or the symmetric KL-divergence between the nominal belief state and the actual belief state exceeds a certain threshold $d_{\text{th}}$.
We set $d_{\text{th}}$ to $1e6$ for this problem, since we noticed that a smaller value would result in re-planning after almost each observation update. With our choice, re-planning almost only happens at the end of the planning horizon, which allows for efficient execution of the policy. More details about the re-planning for T-LQG can be found in \cite{rafieisakhaei2016belief, rafieisakhaei2017tlqg}.

Ergodic exploration is not a belief space planning approach but has been used in the active sensing literature. Beginning with the nominal control of zero, it locally optimizes the ergodicity of the trajectory with respect to the spatial information distribution based on Fisher information. This optimization is open-loop since the effect of future observations is not considered. As a new observation becomes available, the distribution and the trajectory are recomputed. 

All the controllers were saturated at the same limit. The results presented in Figure \ref{fig: tracking_1} clearly indicates superior performance of SACBP while achieving real-time computation. More notably, SACBP generated a trajectory that periodically revisited the two groups whereas other methods failed to do so (Figure \ref{fig: tracking_2}). With SACBP the robot was moving into one of the four diagonal directions for most of the time. This is plausible, as SACBP solves the quadratic program with a box input constraint \eqref{eq: convex program}, which tends to find optimal solutions at the corners. MCTS-DPW resulted in a highly non-smooth trajectory and failed to fully explore the environment. The greedy approach improved the smoothness, but the robot eventually followed a cyclic trajectory in a small region of the environment. To our surprise, the ergodic method did not generate a trajectory that covers the two groups of the targets. This is likely due to the use of a projection-based trajectory optimization method, which has been recently found to perform rather poorly with rapid re-planning \cite{dressel2018efficient}. Among all the baselines implemented, T-LQG performed the worst, leaving large estimation covariances even after the full 200[s]. This is due to the insufficient time budget to solve the NLP, which indicates that T-LQG is not suited for high-dimensional belief space planning problems such as this active multi-target tracking task.

\subsection{Object Manipulation under Model Uncertainty}
This problem is identical to the model-based Bayesian reinforcement learning problem studied in \cite{slade2017manip}, therefore a detailed description of the nonlinear dynamics and the observation models are omitted. See Figure \ref{fig: manipulation_1} for the illustration of the environment. A 2D robot attached to a rigid body object applies forces and torques to move the object to the origin. The object's mass, moment of inertia, moment arm lengths, and linear friction coefficient are unknown. These parameters as well as the object's 2D state need to be estimated using EKF, with noisy sensors which measure the robot's position, velocity, and acceleration in the global frame. The same values for $t_f - t_0$, $\Delta t_o$, $\Delta t_c$, $t_\text{calc}$ as in the previous problem are assumed. Each simulation is run for 20[s]. SACBP uses $\epsilon = 0.04[s]$ and $N = 10$. The nominal control for SACBP is a position controller whose input is the mean x-y position and the rotation estimates of the object. The cost function is quadratic in the true state $x$ and control $u$, given by $\frac{1}{2}x^\mathrm{T}C_x x + \frac{1}{2} u^\mathrm{T}C_u u$. Taking expectations yields the equivalent cost in the Gaussian belief space $c(b,u) = \frac{1}{2} \mu^\mathrm{T}C_x \mu + \frac{1}{2}\text{tr}(C_x \Sigma) + \frac{1}{2}u^\mathrm{T}C_u u$, where $\Sigma$ is the covariance matrix. We let terminal cost $h$ be the same as $c$ except that it is without the control term.

We compared SACBP against (i) MCTS-DPW in the Gaussian belief space, (ii) belief iLQG, and (iii) T-LQG. MCTS-DPW uses the same planning horizon as SACBP, and is set to draw $N=240$ samples so that the computation time is approximately equal to that of SACBP. Furthermore, MCTS-DPW uses the position controller mentioned above as the rollout policy, which is suggested in \cite{slade2017manip}. 

The policy computation with belief iLQG is performed offline over the entire simulation horizon of 20[s]. The solver is initialized with a nominal trajectory generated by the same position controller. The entire computation of the policy takes 5945[s], a significant amount of time until convergence to a locally optimal affine feedback policy. Note however that the online policy execution can be performed instantaneously. 

For T-LQG, we test the algorithm in both the online and the offline modes. The online mode is equivalent to the implementation for the active multi-target tracking problem. Re-planning happens when the end of the planning horizon is reached or the symmetric KL-divergence surpasses $d_{\text{th}} = 25$, which was rarely exceeded during the simulation for this problem. Furthermore, the $\textrm{max\_cpu\_time}$ parameter of Ipopt is adjusted so the computation time is comparable to SACBP and MCTS-DPW. This is because the full optimization over 10 timesteps (i.e. 2[s]) takes 2.5[s] on average and 7[s] in the worst case, which is better than in the active multi-target tracking problem but still prohibitively slow for online control computation, taking about $15\times t_{\mathrm{calc}}$ to $45\times t_{\mathrm{calc}}$ with $t_{\mathrm{calc}} = 0.15$[s].
On the other hand, the offline mode computes the closed-loop policy once for the entire simulation horizon without online re-planning or limiting $\textrm{max\_cpu\_time}$. This setup is identical to belief iLQG. In this mode, it takes T-LQG 1065[s] to compute the policy, which is about 5 to 6 times faster than belief iLQG. This improved efficiency is congruous with the complexity analysis provided in \cite{rafieisakhaei2017tlqg}. Note also that we use the aforementioned position controller to initialize the NLP solver in both the online and the offline modes.

Overall, the results presented in Figure \ref{fig: manipulation_1} demonstrate that SACBP outperformed all the baselines in this task with only 10 Monte Carlo samples, bringing the object close to the goal within 10[s]. Although the computation time increased from the previous problem due to the continuous-discrete filtering, it still achieved near real-time performance with less than $3\times t_{\mathrm{calc}}$[s] on average. Compared to SACBP, reduction of the residual norm was slower for MCTS-DPW and online T-LQG. The two offline algorithms tested, belief iLQG and offline T-LQG, both had a large overshoot at around 2[s] and 3[s], respectively. We suppose that this was caused by the offline nature of the policies, as well as a mismatch between the discrete-time model used for planning and the continuous-time model employed for dynamics simulation. Another large overshoot for belief iLQG at 11[s] was likely due to a locally optimal behavior of the iLQG solver.

\section{Conclusions and Future Work}
\label{sec: conclusions}
In this paper we present SACBP, a novel belief space planning algorithm for continuous-time dynamical systems. We view the stochastic belief dynamics as a hybrid system with time-driven switching and derive the optimal control perturbation based on the perturbation theory of differential equations. The resulting algorithm extends the framework of SAC to stochastic belief dynamics and is highly parallelizable to run in near real-time. The rigorous mathematical analysis shows that the notion of mode insertion gradient can be generalized to stochastic hybrid systems, which leads to the property of SACBP that the algorithm is expected to perform at least as good as the nominal policy with an appropriate choice of the perturbation duration. Through an extensive simulation study, we have confirmed that SACBP outperforms other algorithms including a greedy algorithm and non-myopic closed-loop planners that are based on approximate dynamic programming and/or local trajectory optimization. In future work, we are interested to consider a distributed multi-robot version of SACBP as well as problems with hard state constraints. We also plan to provide additional case studies for more complex belief distributions with efficient implementation.

\begin{funding}
Toyota Research Institute (``TRI")  provided funds to assist the authors with their research but this article solely reflects the opinions and conclusions of its authors and not TRI or any other Toyota entity. This work was also supported in part by NSF grant CMMI1562335, ONR grant N00014-18-1-2830, a JASSO fellowship, and a Masason Foundation fellowship. The authors are grateful for this support.
\end{funding}

\bibliographystyle{SageH}
\bibliography{references.bib}

\begin{thebibliography}{54}
\providecommand{\natexlab}[1]{#1}
\providecommand{\url}[1]{\texttt{#1}}
\providecommand{\urlprefix}{URL }
\expandafter\ifx\csname urlstyle\endcsname\relax
  \providecommand{\doi}[1]{DOI:\discretionary{}{}{}#1}\else
  \providecommand{\doi}{DOI:\discretionary{}{}{}\begingroup
  \urlstyle{rm}\Url}\fi

\bibitem[{{Agha-mohammadi} et~al.(2018){Agha-mohammadi}, {Agarwal}, {Kim},
  {Chakravorty} and {Amato}}]{agha2018slap}
{Agha-mohammadi} A, {Agarwal} S, {Kim} S, {Chakravorty} S and {Amato} NM (2018)
  Slap: Simultaneous localization and planning under uncertainty via dynamic
  replanning in belief space.
\newblock \emph{IEEE Transactions on Robotics} 34(5): 1195--1214.

\bibitem[{Agha-Mohammadi et~al.(2014)Agha-Mohammadi, Chakravorty and
  Amato}]{agha2014firm}
Agha-Mohammadi AA, Chakravorty S and Amato NM (2014) Firm: Sampling-based
  feedback motion-planning under motion uncertainty and imperfect measurements.
\newblock \emph{The International Journal of Robotics Research} 33(2):
  268--304.

\bibitem[{Andersson et~al.(2019)Andersson, Gillis, Horn, Rawlings and
  Diehl}]{andersson2019casadi}
Andersson JAE, Gillis J, Horn G, Rawlings JB and Diehl M (2019) {CasADi} -- {A}
  software framework for nonlinear optimization and optimal control.
\newblock \emph{Mathematical Programming Computation} 11(1): 1--36.
\newblock \doi{10.1007/s12532-018-0139-4}.

\bibitem[{Ansari and Murphey(2016)}]{ansari2016sac}
Ansari AR and Murphey TD (2016) {Sequential Action Control: Closed-Form Optimal
  Control for Nonlinear and Nonsmooth Systems}.
\newblock \emph{IEEE Transactions on Robotics} 32(5): 1196--1214.

\bibitem[{Bajcsy(1988)}]{bajcsy1988activeperception}
Bajcsy R (1988) Active perception.
\newblock \emph{Proceedings of the IEEE} 76(8): 966--1005.

\bibitem[{Bangura and Mahony(2014)}]{bangura2014real}
Bangura M and Mahony R (2014) Real-time model predictive control for
  quadrotors.
\newblock \emph{IFAC Proceedings Volumes} 47(3): 11773--11780.

\bibitem[{Bourbaki and Spain(2004)}]{bourbaki2004theory}
Bourbaki N and Spain P (2004) \emph{Elements of Mathematics Functions of a Real
  Variable: Elementary Theory}.
\newblock Berlin, Heidelberg: Springer Berlin Heidelberg.

\bibitem[{Bourgault et~al.(2002)Bourgault, Makarenko, Williams, Grocholsky and
  Durrant-Whyte}]{bourgault2002information}
Bourgault F, Makarenko A, Williams S, Grocholsky B and Durrant-Whyte H (2002)
  {Information based adaptive robotic exploration}.
\newblock In: \emph{IEEE/RSJ International Conference on Intelligent Robots and
  System}, volume~1. IEEE, pp. 540--545.

\bibitem[{Browne et~al.(2012)Browne, Powley, Whitehouse, Lucas, Cowling,
  Rohlfshagen, Tavener, Perez, Samothrakis and Colton}]{browne2012mcts}
Browne CB, Powley E, Whitehouse D, Lucas SM, Cowling PI, Rohlfshagen P, Tavener
  S, Perez D, Samothrakis S and Colton S (2012) A survey of monte carlo tree
  search methods.
\newblock \emph{IEEE Transactions on Computational Intelligence and AI in
  Games} 4(1): 1--43.

\bibitem[{Chaudhari et~al.(2013)Chaudhari, Karaman, Hsu and
  Frazzoli}]{chaudhari2013sampling}
Chaudhari P, Karaman S, Hsu D and Frazzoli E (2013) Sampling-based algorithms
  for continuous-time pomdps.
\newblock In: \emph{2013 American Control Conference}. IEEE, pp. 4604--4610.

\bibitem[{Cou{\"{e}}toux et~al.(2011)Cou{\"{e}}toux, Hoock, Sokolovska, Teytaud
  and Bonnard}]{couetoux2011cuct}
Cou{\"{e}}toux A, Hoock JB, Sokolovska N, Teytaud O and Bonnard N (2011)
  {Continuous Upper Confidence Trees}.
\newblock In: \emph{2011 International Conference on Learning and Intelligent
  Optimization}. Springer, Berlin, Heidelberg, pp. 433--445.

\bibitem[{Diestel and Uhl(1977)}]{diestel1977vector}
Diestel J and Uhl J (1977) \emph{Vector Measures}.
\newblock American Mathematical Society.

\bibitem[{Dressel and Kochenderfer(2018)}]{dressel2018ergoric}
Dressel L and Kochenderfer MJ (2018) Tutorial on the generation of ergodic
  trajectories with projection-based gradient descent.
\newblock \emph{IET Cyber-Physical Systems: Theory \& Applications} .

\bibitem[{Dressel(2018)}]{dressel2018efficient}
Dressel LK (2018) \emph{Efficient and Low-cost Localization of Radio Sources
  with an Autonomous Drone}.
\newblock PhD Thesis, Stanford University.

\bibitem[{Egerstedt et~al.(2006)Egerstedt, Wardi and
  Axelsson}]{egerstedt2006transition}
Egerstedt M, Wardi Y and Axelsson H (2006) Transition-time optimization for
  switched-mode dynamical systems.
\newblock \emph{IEEE Transactions on Automatic Control} 51(1): 110--115.

\bibitem[{Egorov et~al.(2017)Egorov, Sunberg, Balaban, Wheeler, Gupta and
  Kochenderfer}]{egorov2017pomdps}
Egorov M, Sunberg ZN, Balaban E, Wheeler TA, Gupta JK and Kochenderfer MJ
  (2017) Pomdps. jl: A framework for sequential decision making under
  uncertainty.
\newblock \emph{Journal of Machine Learning Research} 18(26): 1--5.

\bibitem[{Elijah(1997)}]{elijah1997optimization}
Elijah P (1997) \emph{Optimization: Algorithms and Consistent Approximations}.
\newblock Springer Verlage Publications.

\bibitem[{Erez and Smart(2010)}]{erez2010continuouspomdps}
Erez T and Smart WD (2010) A scalable method for solving high-dimensional
  continuous pomdps using local approximation.
\newblock In: \emph{Proceedings of the Twenty-Sixth Conference on Uncertainty
  in Artificial Intelligence}, UAI'10. Arlington, Virginia, United States: AUAI
  Press, pp. 160--167.

\bibitem[{Gowrisankaran(1972)}]{gowrisankaran1972measurability}
Gowrisankaran K (1972) Measurability of functions in product spaces.
\newblock \emph{Proceedings of the American Mathematical Society} 31(2):
  485--488.

\bibitem[{Heemels et~al.(2009)Heemels, Lehmann, Lunze and
  Schutter}]{heemels2009hybrid}
Heemels W, Lehmann D, Lunze J and Schutter BD (2009) Introduction to hybrid
  systems.
\newblock In: Lunze J and Lamnabhi-Lagarrigue F (eds.) \emph{Handbook of Hybrid
  Systems Control -- Theory, Tools, Applications}, chapter~1. Cambridge
  University Press, pp. 3--30.

\bibitem[{Hollinger and Sukhatme(2014)}]{hollinger2014immp}
Hollinger GA and Sukhatme GS (2014) {Sampling-based robotic information
  gathering algorithms}.
\newblock \emph{The International Journal of Robotics Research} 33(9):
  1271--1287.

\bibitem[{Huber(2009)}]{huber2009probabilistic}
Huber M (2009) \emph{Probabilistic framework for sensor management}, volume~7.
\newblock KIT Scientific Publishing.

\bibitem[{Kaelbling et~al.(1998)Kaelbling, Littman and
  Cassandra}]{kaelbling1998planning}
Kaelbling LP, Littman ML and Cassandra AR (1998) Planning and acting in
  partially observable stochastic domains.
\newblock \emph{Artif. Intell.} 101(1-2): 99--134.

\bibitem[{Khalil and Grizzle(2002)}]{khalil2002nonlinear}
Khalil HK and Grizzle JW (2002) \emph{Nonlinear systems}, volume~3.
\newblock Prentice hall Upper Saddle River, NJ.

\bibitem[{Kochenderfer(2015)}]{kochenderfer2015decision}
Kochenderfer MJ (2015) \emph{Decision making under uncertainty: theory and
  application}.
\newblock MIT press.

\bibitem[{Kurniawati et~al.(2008)Kurniawati, Hsu and
  Lee}]{kurniawati2008sarsop}
Kurniawati H, Hsu D and Lee WS (2008) Sarsop: Efficient point-based pomdp
  planning by approximating optimally reachable belief spaces.
\newblock In: \emph{Robotics: Science and Systems (RSS)}. Zurich, Switzerland.

\bibitem[{{Le Ny} and Pappas(2009)}]{leny2009active}
{Le Ny} J and Pappas GJ (2009) {On trajectory optimization for active sensing
  in Gaussian process models}.
\newblock In: \emph{Proceedings of the 48h IEEE Conference on Decision and
  Control (CDC) held jointly with 2009 28th Chinese Control Conference}. IEEE,
  pp. 6286--6292.

\bibitem[{Liberzon(2011)}]{liberzon2012calculus}
Liberzon D (2011) \emph{Calculus of variations and optimal control theory: A
  concise introduction}.
\newblock Princeton University Press.

\bibitem[{Madani et~al.(1999)Madani, Hanks and
  Condon}]{madani1999undecidability}
Madani O, Hanks S and Condon A (1999) On the undecidability of probabilistic
  planning and infinite-horizon partially observable markov decision problems.
\newblock In: \emph{AAAI/IAAI}. pp. 541--548.

\bibitem[{Mavrommati et~al.(2018)Mavrommati, Tzorakoleftherakis, Abraham and
  Murphey}]{mavrommati2018ergodic}
Mavrommati A, Tzorakoleftherakis E, Abraham I and Murphey TD (2018) Real-time
  area coverage and target localization using receding-horizon ergodic
  exploration.
\newblock \emph{IEEE Transactions on Robotics} : 62--80.

\bibitem[{Mihaylova et~al.(2002)Mihaylova, Lefebvre, Bruyninckx, Gadeyne and
  {De Schutter}}]{mihaylova2002activesensing}
Mihaylova L, Lefebvre T, Bruyninckx H, Gadeyne K and {De Schutter} J (2002)
  {Active Sensing for Robotics - A Survey}.
\newblock In: \emph{Proceedings of 5th International Conference on Numerical
  Methods and Applications}. pp. 316--324.

\bibitem[{Miller and Murphey(2013)}]{miller2013ergodic}
Miller LM and Murphey TD (2013) Trajectory optimization for continuous ergodic
  exploration.
\newblock In: \emph{American Control Conference (ACC), 2013}. IEEE, pp.
  4196--4201.

\bibitem[{Miller et~al.(2016)Miller, Silverman, MacIver and
  Murphey}]{miller2016ergodic}
Miller LM, Silverman Y, MacIver MA and Murphey TD (2016) Ergodic exploration of
  distributed information.
\newblock \emph{IEEE Transactions on Robotics} 32(1): 36--52.

\bibitem[{Nishimura and Schwager(2018{\natexlab{a}})}]{nishimura2018activembc}
Nishimura H and Schwager M (2018{\natexlab{a}}) Active motion-based
  communication for robots with monocular vision.
\newblock In: \emph{2018 IEEE International Conference on Robotics and
  Automation (ICRA)}. pp. 2948--2955.

\bibitem[{Nishimura and Schwager(2018{\natexlab{b}})}]{nishimura2018sacbp}
Nishimura H and Schwager M (2018{\natexlab{b}}) Sacbp: Belief space planning
  for continuous-time dynamical systems via stochastic sequential action
  control.
\newblock In: \emph{The 13th International Workshop on the Algorithmic
  Foundations of Robotics (WAFR)}. M{\'e}rida, M{\'e}xico.

\bibitem[{Papadimitriou and Tsitsiklis(1987)}]{papadimitriou1987complexity}
Papadimitriou CH and Tsitsiklis JN (1987) The complexity of markov decision
  processes.
\newblock \emph{Mathematics of operations research} 12(3): 441--450.

\bibitem[{Patil et~al.(2014)Patil, Kahn, Laskey, Schulman, Goldberg and
  Abbeel}]{patil2014gaussianbsp}
Patil S, Kahn G, Laskey M, Schulman J, Goldberg K and Abbeel P (2014) Scaling
  up gaussian belief space planning through covariance-free trajectory
  optimization and automatic differentiation.
\newblock In: \emph{{WAFR}}, \emph{Springer Tracts in Advanced Robotics},
  volume 107. Springer, pp. 515--533.

\bibitem[{Platt(2013)}]{platt2013convex}
Platt R (2013) Convex receding horizon control in non-gaussian belief space.
\newblock In: \emph{Algorithmic Foundations of Robotics X}. Springer, pp.
  443--458.

\bibitem[{Platt et~al.(2010)Platt, Tedrake, Kaelbling and
  Lozano-Perez}]{platt2010bsp}
Platt R, Tedrake R, Kaelbling L and Lozano-Perez T (2010) Belief space planning
  assuming maximum likelihood observations.
\newblock In: \emph{Robotics: Science and Systems (RSS)}.

\bibitem[{Popović et~al.(2017)Popović, Hitz, Nieto, Sa, Siegwart and
  Galceran}]{popovic2017uav}
Popović M, Hitz G, Nieto J, Sa I, Siegwart R and Galceran E (2017) Online
  informative path planning for active classification using uavs.
\newblock In: \emph{2017 IEEE International Conference on Robotics and
  Automation (ICRA)}. pp. 5753--5758.

\bibitem[{Rafieisakhaei et~al.(2016)Rafieisakhaei, Chakravorty and
  Kumar}]{rafieisakhaei2016belief}
Rafieisakhaei M, Chakravorty S and Kumar P (2016) Belief space planning
  simplified: Trajectory-optimized lqg (t-lqg).
\newblock \emph{arXiv preprint arXiv:1608.03013} .

\bibitem[{Rafieisakhaei et~al.(2017)Rafieisakhaei, Chakravorty and
  Kumar}]{rafieisakhaei2017tlqg}
Rafieisakhaei M, Chakravorty S and Kumar PR (2017) T-lqg: Closed-loop belief
  space planning via trajectory-optimized lqg.
\newblock In: \emph{2017 IEEE International Conference on Robotics and
  Automation (ICRA)}. pp. 649--656.

\bibitem[{{Revels} et~al.(2016){Revels}, {Lubin} and
  {Papamarkou}}]{revels2016forwarddiff}
{Revels} J, {Lubin} M and {Papamarkou} T (2016) Forward-mode automatic
  differentiation in julia.
\newblock \emph{arXiv:1607.07892 [cs.MS]} .

\bibitem[{Schwager et~al.(2017)Schwager, Dames, Rus and
  Kumar}]{schwager2017info}
Schwager M, Dames P, Rus D and Kumar V (2017) A multi-robot control policy for
  information gathering in the presence of unknown hazards.
\newblock In: \emph{Robotics Research : The 15th International Symposium ISRR}.
  Springer International Publishing, pp. 455--472.

\bibitem[{Seekircher et~al.(2011)Seekircher, Laue and
  R{\"{o}}fer}]{seekircher2011entropy}
Seekircher A, Laue T and R{\"{o}}fer T (2011) Entropy-based active vision for a
  humanoid soccer robot.
\newblock In: \emph{RoboCup 2010: Robot Soccer World Cup XIV}, volume 6556
  LNCS. Springer Berlin Heidelberg, pp. 1--12.

\bibitem[{Slade et~al.(2017)Slade, Culbertson, Sunberg and
  Kochenderfer}]{slade2017manip}
Slade P, Culbertson P, Sunberg Z and Kochenderfer M (2017) Simultaneous active
  parameter estimation and control using sampling-based bayesian reinforcement
  learning.
\newblock In: \emph{2017 IEEE/RSJ International Conference on Intelligent
  Robots and Systems (IROS)}. pp. 804--810.

\bibitem[{Somani et~al.(2013)Somani, Ye, Hsu and Lee}]{somani2013despot}
Somani A, Ye N, Hsu D and Lee WS (2013) Despot: Online pomdp planning with
  regularization.
\newblock In: \emph{Proceedings of the 26th International Conference on Neural
  Information Processing Systems - Volume 2}, NIPS'13. USA: Curran Associates
  Inc., pp. 1772--1780.

\bibitem[{Spinello and Stilwell(2010)}]{spinello2010statedepedent}
Spinello D and Stilwell DJ (2010) Nonlinear estimation with state-dependent
  gaussian observation noise.
\newblock \emph{IEEE Transactions on Automatic Control} 55(6): 1358--1366.

\bibitem[{Sunberg and Kochenderfer(2017)}]{sunberg2017pomcpow}
Sunberg Z and Kochenderfer MJ (2017) {POMCPOW:} an online algorithm for pomdps
  with continuous state, action, and observation spaces.
\newblock \emph{CoRR} abs/1709.06196.

\bibitem[{van~den Berg et~al.(2012)van~den Berg, Patil and
  Alterovitz}]{vandenberg2012beliefilqg}
van~den Berg J, Patil S and Alterovitz R (2012) Motion planning under
  uncertainty using iterative local optimization in belief space.
\newblock \emph{The International Journal of Robotics Research} 31(11):
  1263--1278.

\bibitem[{W{\"a}chter and Biegler(2006)}]{wachter2006implementation}
W{\"a}chter A and Biegler LT (2006) On the implementation of an interior-point
  filter line-search algorithm for large-scale nonlinear programming.
\newblock \emph{Mathematical programming} 106(1): 25--57.

\bibitem[{Wardi and Egerstedt(2012)}]{Wardi2012switched}
Wardi Y and Egerstedt M (2012) Algorithm for optimal mode scheduling in
  switched systems.
\newblock In: \emph{2012 American Control Conference (ACC)}. pp. 4546--4551.

\bibitem[{Williams et~al.(2018)Williams, Drews, Goldfain, Rehg and
  Theodorou}]{williams2018information}
Williams G, Drews P, Goldfain B, Rehg JM and Theodorou EA (2018)
  Information-theoretic model predictive control: Theory and applications to
  autonomous driving.
\newblock \emph{IEEE Transactions on Robotics} 34(6): 1603--1622.

\bibitem[{Xie et~al.(2007)Xie, Popa and Lewis}]{xie2007estimation}
Xie L, Popa D and Lewis FL (2007) \emph{Optimal and robust estimation: with an
  introduction to stochastic control theory}.
\newblock CRC press.

\end{thebibliography}

\appendix
\section{Detailed Analysis of Mode Insertion Gradient for Stochastic Hybrid Systems}
In this appendix, we provide a thorough analysis of the stochastic hybrid systems with time-driven switching that satisfy Assumptions \ref{assum: control} -- \ref{assum: observations}. Our goal is to prove Theorem \ref{theorem: mode_insertion_gradient}.
\label{sec: appendix}
\subsection{Nominal Trajectory under Specific Observations}
\label{sec: analysis_nominal_specific}
First, we analyze the properties of the system $x(t)$ for $t \in [0, T]$ under a given initial condition $x_0$, control $u \in U$, and a specific sequence of observations $(y_1, \dots, y_T)$ sampled from $(Y_1,\dots,Y_T)$.
\begin{prop} [Existence and Uniqueness of Solutions]
    \label{prop: existence}
    Given a control $u \in U$ and a sequence of observations $(y_1, \dots, y_T)$, the system $x(t)$ starting at $x_0$ has a unique solution for $t \in [0, T]$.
\begin{proof}
    We will show that each $x_i$ for $i \in \{1, \dots, T\}$ has a unique solution, and thus $x$ is uniquely determined as a whole. First, by Assumption \ref{assum: function_f}, \ref{assum: function_f_bound} and the Picard Lemma (Lemma 5.6.3 in \cite{elijah1997optimization}), the differential equation
    \begin{align}
        \dot{x}_1(t) = f(x_1(t), u(t))
    \end{align}
    with initial condition $x_1(0) = x_0$ has a solution for $t \in [0, 1]$. Furthermore, Proposition 5.6.5 in \cite{elijah1997optimization} assures that the solution $x_1$ is unique under Assumption \ref{assum: function_f} and \ref{assum: function_f_bound}. This guarantees that the initial condition for $x_2$ defined by $x_2(1) = g(x_1(1), y_1)$ is unique. Therefore, proceeding by induction each $x_1, \dots, x_T$ has a unique solution, which completes the proof.
\end{proof}
\end{prop}

\begin{cor}[Right Continuity]
    \label{cor: right_continuity}
    Given a control $u \in U$ and a sequence of observations $(y_1, \dots, y_T)$, the system $x(t)$ starting at $x_0$ is right continuous in $t$ on $[0, T]$.
\begin{proof}
    By Proposition \ref{prop: existence} each $x_i$ has a unique solution that follows $\dot{x}_i = f(x, u)$. Clearly each $x_i$ is continuous on $[i-1, i]$, which proves that $x(t) \triangleq x_i(t)\;  \forall t \in [i-1, i)\; \forall i \in \{1,2,\dots, T\}$ with $x(T) \triangleq g(x_T(T), y_T)$ is right continuous on $[0, T]$.
\end{proof}
\end{cor}

\begin{lemma}
    \label{lemma: bound_1}
    Let $\xi_i$ denote the initial condition for $x_i$. Then, there exists a constant $K_8 < \infty$ such that for all $i \in \{1, \dots, T\}$, 
    \begin{align}
        \forall t \in [i-1, i]\;\; \Vert x_i(t) \Vert_2 \leq \left(1 + \Vert \xi_i \Vert_2\right) e^{K_8}
    \end{align}
\begin{proof}
    Using Assumption \ref{assum: function_f} and \ref{assum: function_f_bound}, the claim follows directly from Proposition 5.6.5 in \cite{elijah1997optimization}.
\end{proof}
\end{lemma}

\begin{prop}[Lipschitz Continuity]
    \label{prop: lipschitz_continuity}
    For each $i \in \{1, \dots, T\}$, let $\xi_i^\prime$ and $\xi_i^{\prime\prime}$ be two distinct initial conditions for $x_i$. Furthermore, let $u^\prime$ and $u^{\prime\prime}$ be two controls from $U$. The pairs $(\xi^{\prime}, u^{\prime})$ and $(\xi^{\prime\prime}, u^{\prime\prime}$) respectively define two solutions $x_i^\prime$ and $x_i^{\prime\prime}$ to ODE $\dot{x}_i = f(x_i, u)$ over $[i-1, i]$. Then, there exists an $L < \infty$, independent of $\xi_i^\prime, \xi_i^{\prime\prime}, u^\prime$ and $u^{\prime\prime}$, such that
    \begin{multline}
        \forall i \in \{1, \dots, T\}\; \forall t \in [i-1, i]\;\;
        \Vert x_i^\prime(t) - x_i^{\prime\prime}(t) \Vert_2 \leq \\
        L\left(\Vert \xi_i^\prime - \xi_i^{\prime\prime} \Vert_2 + \int_{i-1}^{i} \Vert u^\prime(t) - u^{\prime\prime}(t) \Vert_2 dt \right).
    \end{multline}
\begin{proof}
    The proof is similar to that of Lemma 5.6.7 in \cite{elijah1997optimization}. Making use of the Picard Lemma (Lemma 5.6.3 in \cite{elijah1997optimization}) and Assumption \ref{assum: function_f_bound}, we obtain
    \begin{multline}
        \Vert x_i^\prime(t) - x_i^{\prime\prime}(t) \Vert_2 \leq \\
        e^{K_1}\left (\Vert\xi_i^\prime - \xi_i^{\prime\prime} \Vert_2 
        + K_1 \int_{i-1}^i \Vert u^\prime(t) - u^{\prime\prime}(t) \Vert_2 dt\right).
    \end{multline}
    As $K_1 \geq 1$,
    \begin{multline}
        \Vert x_i^\prime(t) - x_i^{\prime\prime}(t) \Vert_2 \leq\\
        K_1 e^{K_1}\left (\Vert\xi_i^\prime - \xi_i^{\prime\prime} \Vert_2 
        + \int_{i-1}^i \Vert u^\prime(t) - u^{\prime\prime}(t) \Vert_2 dt\right).
    \end{multline}
    Defining $L \triangleq K_1 e^{K_1} < \infty$ completes the proof.
\end{proof}
\end{prop}

\begin{cor} [Uniform Continuity in Initial Conditions]
    \label{cor: uniform_continuity}
    Let $u \in U$ be a given control. For each $i \in \{1, \dots, T\}$, let $\xi_i^\prime$ and $\xi_i^{\prime\prime}$ be two distinct initial conditions for $x_i$. The pairs $(\xi^{\prime}, u)$ and $(\xi^{\prime\prime}, u)$ respectively define two solutions $x_i^\prime$ and $x_i^{\prime\prime}$ to ODE $\dot{x}_i = f(x_i, u)$ over $[i-1, i]$. Note that they share the same control but have different initial conditions, unlike Proposition \ref{prop: lipschitz_continuity}. Then, for any $\epsilon > 0$ there exists $\delta > 0$ such that
    \begin{multline}
        \forall i \in \{1, \dots, T\}\; \forall t \in [i-1, i] \;\; \Vert \xi_i^\prime - \xi_i^{\prime\prime} \Vert_2 < \delta \\
        \Rightarrow \Vert x_i^\prime(t) - x_i^{\prime\prime}(t) \Vert_2 < \epsilon.
    \end{multline}
\begin{proof}
    By Proposition \ref{prop: lipschitz_continuity}, we find
    \begin{align}
        \Vert x_i^\prime(t) - x_i^{\prime\prime}(t) \Vert_2 \leq L \Vert \xi_i^\prime - \xi_i^{\prime\prime} \Vert_2
    \end{align}
    for all $i \in \{1, \dots, T\}$ and $t \in [i-1, i]$, where $L < \infty$. Take any $\delta < \frac{\epsilon}{L}$ to prove the claim.
\end{proof}
\end{cor}

\begin{prop} [Continuity in Observations]
    \label{prop: continuity_in_observations}
    Given a control $u \in U$, the map $(y_1,\dots,y_T) \mapsto x(t)$ is continuous for all $t \in [0, T]$, where $x$ represents the solution to the system under Assumption \ref{assum: dynamics} starting at $x(0) = x_0$.
\begin{proof}
    We will show the continuity of $(y_1, \dots, y_T) \mapsto x_i(t)$ for each $i \in \{1,\dots,T\}$ by mathematical induction. First, by Assumption \ref{assum: dynamics} the value of $x_1(t)$ is solely determined by $x_0$ and $u$. Therefore, for any $t \in [0, 1]$ the function $(y_1, \dots, y_T) \mapsto x_1(t)$ is a constant map, which is continuous. Next, suppose that $\forall t \in [i-1, i]\;\; (y_1, \dots, y_T) \mapsto x_i(t)$ is continuous for some $i \in \{1,\dots, T-1\}$. Now consider $x_{i+1}$. Let $F_{i+1}(\xi_{i+1}, t)$ be the map from an initial condition $x_{i+1}(i) = \xi_{i+1}$ to the solution $x_{i+1}$ at $t \in [i, i+1]$ under the given $u$. In other words, $F_{i+1}(\xi_{i+1}, t)$ is equivalent to the integral equation
    \begin{align}
        F_{i+1}(\xi_{i+1}, t) \triangleq \xi_{i+1} + \int_{i}^t f(x_{i+1}(a), u(a)) da
    \end{align}
    and takes initial condition $\xi_{i+1}$ as well as time $t$ as its arguments. 
    Substituting $\xi_{i+1} = g(x_i(i), y_i)$, $F_{i+1}\left(g(x_i(i), y_i), t\right)$ gives the actual value of $x_{i+1}(t)$. We will prove the continuity of $F_{i+1}(g(x_i(i), y_i), t)$ as follows. First, note that the map from $(y_1,\dots,y_T)$ to $F_{i+1}\left(g(x_i(i), y_i), t\right)$ is the result of the composition:
    \begin{align}
        \begin{pmatrix} y_1 \\ \vdots \\ y_T\end{pmatrix} \mapsto \begin{pmatrix} x_i(i) \\  y_i\end{pmatrix} \mapsto g(x_i(i), y_i) \mapsto F_{i+1}(g(x_i(i), y_i), t).
    \end{align}
    The first map is continuous since $x_i(i)$ is continuous in $(y_1, \dots, y_T)$ by the induction hypothesis. The second map is also continuous by Assumption \ref{assum: function_g}. Lastly, Corollary \ref{cor: uniform_continuity} shows that $F_{i+1}(\xi_{i+1}, t)$ is (uniformly) continuous in $\xi_{i+1}$ for $t \in [i, i+1]$. Therefore, $(y_1,\dots,y_T) \mapsto x_{i+1}(t)$ is continuous for all $t \in [i, i+1]$.
\end{proof}
\end{prop}

\begin{prop} [Bounded State Trajectory]
    \label{prop: bound_1}
    Given a control $u \in U$ and a sequence of observations $(y_1, \dots, y_T)$, the system $x(t)$ starting at $x(0) = x_1(0) = x_0$ has the following bound:
    \begin{multline}
        \forall i \in \{2,\dots, T\}\; \forall t \in [i-1, i] \;\; \\
        \Vert x_i(t) \Vert_2 \leq \sum_{(j_1, \dots, j_{i-1}) \in \mathcal{K}_i} \alpha_i^{(j_1,\dots,j_{i-1})}(x_0) \prod_{m=1}^{i-1}\Vert y_m \Vert_2^{j_m},
    \end{multline}
    where $\mathcal{K}_i$ is a finite set of sequences of non-negative integers of length $i-1$, and $\alpha_i^{(j_1,\dots,j_{i-1})}(x_0)$ is a finite positive constant that depends on $x_0$ and $(j_1,\dots, j_{i-1})$ but not on any of the observations or the control. 
    
    For $i = 1$ the bound is given by $\forall t \in [0, 1]\;\; \Vert x_1(t) \Vert_2 \leq \alpha_1 (x_0)$ for some finite positive constant $\alpha_1(x_0)$.
\begin{proof}
    For $i = 1$, Lemma \ref{lemma: bound_1} gives $\forall t \in [0, 1]\;\;$ $\Vert x_1(t) \Vert_2 \leq \left(1 + \Vert x_0 \Vert_2\right) e^{K_8} \triangleq \alpha_1(x_0)$. For $i = 2$, by Assumption \ref{assum: function_g_bound} and the case for $i = 1$, we have
    \begin{align}
        \Vert x_2(1) \Vert_2 &= \Vert g(x_1(1), y_1) \Vert_2 \\ 
        &\begin{multlined}
            \leq K_2 + K_3 \Vert x_1(1) \Vert_2^{L_1} + K_4 \Vert y_1 \Vert_2^{L_2} \\ + K_5 \Vert x_1(1) \Vert_2^{L_1} \Vert y_1 \Vert_2^{L_2}
        \end{multlined} \\
        &\begin{multlined}
            \leq K_2 + K_3\alpha_1(x_0)^{L_1} + K_4\Vert y_1\Vert_2^{L_2} \\ 
            + K_5 \alpha_1(x_0)^{L_1} \Vert y_1\Vert_2^{L_2}.
        \end{multlined}
    \end{align}
    Then, by Lemma \ref{lemma: bound_1}, $\forall t \in [1, 2]\;$
    \begin{align}
        \Vert x_2(t) \Vert_2 &\leq (1 + \Vert x_2(1) \Vert_2)e^{K_8} \\
        &\begin{multlined}
            \leq e^{K_8}\bigg(1 + K_2 + K_3\alpha_1(x_0)^{L_1} + K_4 \Vert y_1 \Vert_2^{L_2} \\ 
            + K_5 \alpha_1(x_0)^{L_1} \Vert y_1 \Vert_2^{L_2}\bigg)
        \end{multlined} \\
        &\triangleq  \sum_{(j_1) \in \mathcal{K}_2} \alpha_2^{(j_1)}(x_0) \Vert y_1 \Vert_2^{j_1},
    \end{align}
    where $\mathcal{K}_2  = \{(0), (L_2)\}$, and
    \begin{align}
        \alpha_2^{(0)}(x_0) &= e^{K_8}\left( 1 + K_2 + K_3 \alpha_1(x_0)^{L_1}\right) \\
        \alpha_2^{(L_2)}(x_0) &= e^{K_8}\left(K_4 + K_5\alpha_1(x_0)^{L_1}\right)
    \end{align}
    are both finite positive constants that depend on $x_0$ but not on any of the observations or the control.
    
    Next, suppose that the claim holds for some $i \geq 2$. That is,
    \begin{align}
        \label{eq: ih_1}
        &\forall t \in [i-1, i] \;\; \nonumber \\
        &\Vert x_i(t) \Vert_2 \leq \sum_{(j_1, \dots, j_{i-1}) \in \mathcal{K}_i} \alpha_i^{(j_1,\dots,j_{i-1})}(x_0) \prod_{m=1}^{i-1}\Vert y_m \Vert_2^{j_m},
    \end{align}
    where $\mathcal{K}_i$ and $\alpha_i^{(j_1,\dots,j_{i-1})}(x_0)$ are as defined in the statement of the proposition.
    Making use of this assumption, Assumption \ref{assum: function_g_bound} and Lemma \ref{lemma: bound_1}, we find that for all $t \in [i, i+1]$,
    \begin{align}
        \Vert x_{i+1}(t) \Vert_2 &\leq (1 + \Vert g(x_i(i), y_i)\Vert_2)e^{K_8} \\
        &\begin{multlined}
            \label{eq: bound_1_terms}
            \leq e^{K_8}(1 + K_2) + e^{K_8}K_4 \Vert y_i\Vert_2^{L_2} \\
            + e^{K_8} \Vert x_i(i)\Vert_2^{L_1}(K_3 + K_5\Vert y_i\Vert_2^{L_2}).
        \end{multlined}
    \end{align}
    The first two terms in the above sum can be rewritten as
    \begin{align}
        e^{K_8}(1 + K_2) &= e^{K_8}(1 + K_2) \Vert y_1\Vert_2^0 \times\cdots\times \Vert y_i\Vert_2^0 \label{eq: bound_1_term_1}\\
        e^{K_8} K_4 \Vert y_i \Vert_2^{L_2} &= e^{K_8} K_4 \Vert y_1 \Vert_2^0 \times\cdots\times \Vert y_{i-1}\Vert_2^0 \times \Vert y_i \Vert_2^{L_2} \label{eq: bound_1_term_2}
    \end{align}
    For the last term, we can use \eqref{eq: ih_1} and the multinomial theorem to write
    \begin{align}
        &\begin{multlined}
            \Vert x_i(i) \Vert_2^{L_1} \\
            \leq \left(\sum_{(j_1, \dots, j_{i-1}) \in \mathcal{K}_i} \alpha_i^{(j_1,\dots,j_{i-1})}(x_0) \prod_{m=1}^{i-1}\Vert y_m \Vert_2^{j_m}\right)^{L_1}
        \end{multlined}\\
        &\begin{multlined}
            = \sum_{k_1 + \cdots + k_{|\mathcal{K}_i|} = L_1} \begin{pmatrix} L_1 \\ k_1,\dots,k_{|\mathcal{K}_i|} \end{pmatrix} \\
            \times \left\{\prod_{l=1}^{|\mathcal{K}_i|} \alpha_i^{(j_1^{(l)},\dots, j_{i-1}^{(l)})}(x_0)^{k_l}\right\}\times \prod_{m=1}^{i-1}\Vert y_m \Vert_2^{\sum_{l=1}^{|\mathcal{K}_i|} k_l j_m^{(l)}}, \label{eq: bound_1_term_3}
        \end{multlined}
    \end{align}
    where $(j_1^{(l)},\dots,j_{i-1}^{(l)})$ is the $l$-th element in $\mathcal{K}_i$.
    Note that $k_l$ is non-negative for all $l \in \{1,\dots, |\mathcal{K}_i|\}$. By the induction hypothesis \eqref{eq: ih_1}, exponent $\sum_{l=1}^{|\mathcal{K}_i|} k_l j_m^{(l)}$ is also non-negative for all $m \in \{1,\dots, i-1\}$.
    
    Thus, substituting \eqref{eq: bound_1_term_1}, \eqref{eq: bound_1_term_2}, and \eqref{eq: bound_1_term_3} into \eqref{eq: bound_1_terms}, rearranging the sums and re-labeling the sequences of integer exponents, we can write
    \begin{align}
        \Vert x_{i+1}(t) \Vert_2 \leq \sum_{(j_1, \dots, j_i) \in \mathcal{K}_{i+1}} \alpha_{i+1}^{(j_1,\dots,j_i)}(x_0) \prod_{m=1}^i \Vert y_m \Vert_2^{j_m}
    \end{align}
    for all $t \in [i, i+1]$, where $\mathcal{K}_{i+1}$ is a set of sequences of non-negative integers of length $i$, and each $\alpha_{i+1}^{(j_1,\dots,j_i)}(x_0)$ does not depend on any of the observations or the control. Here the cardinality of $\mathcal{K}_{i+1}$ is at most finite, since
    \begin{align}
        | \mathcal{K}_{i+1} | \leq 2 + 2 \begin{pmatrix} L_1 + |\mathcal{K}_i| - 1 \\ |\mathcal{K}_i| - 1 \end{pmatrix} 
    \end{align}
    by \eqref{eq: bound_1_terms} and the multinomial theorem.
    
    Finally, proceeding by mathematical induction over $i \in \{2, \dots, T\}$ completes the proof.
\end{proof}
\end{prop}

\begin{prop} [Bounded Cost Functions]
    \label{prop: bounded_cost}
    Given a control $u \in U$ and a sequence of observations $(y_1, \dots, y_T)$, the instantaneous cost $c(x(t), u(t))$ induced by the state trajectory $x(t)$ has the following bound.
    \begin{multline}
        \forall i \in \{2,\dots,T\}\; \forall t \in [i-1, i]\;\; |c(x_i(t), u(t))| \\
        \leq \sum_{(j_1,\dots,j_{i-1}) \in \mathcal{K}^{\prime}_i} \alpha_i^{\prime (j_1,\dots,j_{i-1})}(x_0) \prod_{m=1}^{i-1} \Vert y_m \Vert_2^{j_m},
    \end{multline}
    where $\mathcal{K}^{\prime}_i$ is a finite set of sequences of non-negative integers of length $i-1$, and $\alpha_i^{\prime (j_1,\dots,j_{i-1})}(x_0)$ is a finite positive constant that depends on $x_0$ and $(j_1,\dots, j_{i-1})$ but not on any of the observations or the control. 
    
    For $i = 1$ the bound is given by
    \begin{align}
        \forall t \in [0, 1]\;\; |c(x_1(t), u(t))| \leq \alpha^{\prime}_1 (x_0)
    \end{align}
    for some finite positive constant $\alpha^{\prime}_1(x_0)$.
    
    Similarly, the terminal cost $h(x(T))$ is bounded by
    \begin{align}
        |h(x(T))| \leq \sum_{(j_1,\dots,j_T) \in \mathcal{K}^{\prime}_{T+1}} \alpha^{\prime(j_1,\dots,j_T)}_{T+1}(x_0) \prod_{m=1}^T \Vert y_m \Vert_2^{j_m}.
    \end{align}
\begin{proof}
    For the instantaneous cost function $c(x_i(t), u(t))$, Assumption \ref{assum: cost} along with Proposition \ref{prop: bound_1} yields the following bound:
    \begin{multline}
        \forall i \in \{2,\dots,T\}\; \forall t \in [i-1, i]\;\; |c(x_i(t), u(t))| \leq K_6 \\
        + K_7 \left(\sum_{(j_1, \dots, j_{i-1}) \in \mathcal{K}_i} \alpha_i^{(j_1,\dots,j_{i-1})}(x_0) \prod_{m=1}^{i-1}\Vert y_m \Vert_2^{j_m}\right)^{L_3}.
    \end{multline}
    Making use of the multinomial expansion formula in the same manner as in the proof of Proposition \ref{prop: bound_1}, we conclude that
    \begin{align}
        \forall i \in \{2,\dots,T\}\; \forall t \in [i-1, i]\;\; |c(x_i(t), u(t))| \nonumber \\
        \leq \sum_{(j_1,\dots,j_{i-1}) \in \mathcal{K}^{\prime}_i} \alpha_i^{\prime (j_1,\dots,j_{i-1})}(x_0) \prod_{m=1}^{i-1} \Vert y_m \Vert_2^{j_m}
    \end{align}
    for some finite set $\mathcal{K}^{\prime}_{i}$ of sequences of non-negative integers and finite positive constants $\alpha_i^{\prime(j_1,\dots,j_{i-1})}(x_0)$.
    Similarly, for $i=1$ we obtain
    \begin{align}
        |c(x_1(t), u(t))| \leq K_6 + K_7\alpha_1(x_0)^{L_3} \triangleq \alpha_1^{\prime}(x_0)
    \end{align}
    for all $t \in [0, 1]$.
    
    To bound the terminal cost, note that
    \begin{align}
        |h(x(T))| &\leq K_6 + K_7 \Vert x(T) \Vert_2^{L_3} \\
        &= K_6 + K_7 \Vert g(x_T(T), y_T) \Vert_2^{L_3} \\
        &\begin{multlined}
        \label{eq: bound_2}
        \leq K_6 + K_7\bigg\{K_2 + K_4 \Vert y_T \Vert_2^{L_2} \\
        + \Vert x_T(T) \Vert_2^{L_1}\left(K_3 + K_5 \Vert y_T \Vert_2^{L_2}\right)\bigg\}^{L_3}
        \end{multlined}
    \end{align}
    by Assumptions \ref{assum: function_g_bound} and \ref{assum: cost}. Since $L_3 < \infty$, \eqref{eq: bound_2} yields a polynomial of $\Vert x_T(T) \Vert_2$ and $\Vert y_T \Vert_2$ of finite terms, for each of which we can use Proposition \ref{prop: bound_1} and apply the multinomial expansion formula to show
    \begin{align}
        |h(x(T))| \leq \sum_{(j_1,\dots,j_T) \in \mathcal{K}^{\prime}_{T+1}} \alpha^{\prime(j_1,\dots,j_T)}_{T+1}(x_0) \prod_{m=1}^T \Vert y_m \Vert_2^{j_m}
    \end{align}
    for some finite set $\mathcal{K}^{\prime}_{T+1}$ of sequence of non-negative integers and finite positive constants $\alpha_{T+1}^{\prime(j_1,\dots,j_T)}(x_0)$.
\end{proof}
\end{prop}

\subsection{Perturbed Trajectory under Specific Observations}
Next, we will perturb the nominal state $x$ of the system while assuming the same initial condition $x_0$ and the specific observations $(y_1,\dots, y_T)$ as in Section \ref{sec: analysis_nominal_specific}. The perturbed control is an open-loop perturbation as defined below.
\begin{defn}[Perturbed Control]
    \label{def: perturbed_control}
    Let $u \in U$ be a control. For $\tau \in (0, 1)$ and $v \in B(0, \rho_{\max})$, define the perturbed control $u^{\epsilon}$ by
    \begin{align}
        \label{eq: perturbed_control}
        u^{\epsilon}(t) \triangleq \begin{cases}
            v & \text{if}\;\; t \in (\tau - \epsilon, \tau] \\
            u(t) & \text{otherwise},
        \end{cases}
    \end{align}
    where $\epsilon \in [0, \tau]$. By definition if $\epsilon = 0$ then $u^\epsilon$ is the same as $u$. We assume that the nominal control $u(t)$ is left continuous in $t$ at $t = \tau$. 
\end{defn}
\begin{remark}
    It is obvious that $u^{\epsilon}(t)$ is piecewise continuous on $[0, T]$. Therefore, for $v \in B(0,\rho_{\max})$ we have $u^\epsilon \in U$, i.e. $u^\epsilon$ is an admissible control. Thus, Proposition \ref{prop: existence} assures that there exists a unique solution $x^{\epsilon}$ for the trajectory of the system under the control perturbation. In the remainder of the analysis, we assume that $(\tau, v)$ is given and fixed.
\end{remark}

\begin{lemma}
    \label{lemma: lipschitz_continuity}
    Let $\epsilon, \epsilon^{\prime} \in [0, \tau]$, and let $u^{\epsilon}, u^{\epsilon^{\prime}} \in U$ be two perturbed controls of the form \eqref{eq: perturbed_control}. Let $x_1^{\epsilon}, x_1^{\epsilon^{\prime}}$ be the solutions of $x_1$ for $t \in [0, 1]$ by applying $u^{\epsilon}$ and $u^{\epsilon^{\prime}}$ respectively to the initial condition $x_0$. Then, there exists an $L^{\prime} < \infty$, independent of $\epsilon, \epsilon^{\prime}, x_1^{\epsilon}$ and $x_1^{\epsilon^{\prime}}$, such that
    \begin{align}
        \forall \epsilon, \epsilon^{\prime} \in [0, \tau]\; \forall t \in [0, 1]\;\; \Vert x_1^{\epsilon}(t) - x_1^{\epsilon^{\prime}}(t) \Vert_2 \leq L^{\prime} | \epsilon - \epsilon^{\prime}|.
    \end{align}
\begin{proof}
    By Proposition \ref{prop: lipschitz_continuity}, we find that
    \begin{align}
        \label{eq: lipschitz_perturbation_1}
        \forall t \in [0, 1]\;\; \Vert x_1^{\epsilon}(t) - x_1^{\epsilon^{\prime}}(t) \Vert_2 &\leq L \int_{0}^{1} \Vert u^{\epsilon}(t) - u^{\epsilon^{\prime}}(t) \Vert_2 dt
    \end{align}
    for some $L < \infty$. Let us derive an upper-bound on the integral on the right hand side. If $\epsilon \geq \epsilon^{\prime}$,
    \begin{align}
        \int_{0}^{1} \Vert u^{\epsilon}(t) - u^{\epsilon^{\prime}}(t) \Vert_2 dt &= \int_{\tau - \epsilon}^{\tau - \epsilon^{\prime}} \Vert v - u(t) \Vert_2 dt,
    \end{align}
    where $u(t)$ is the nominal control that both $u^{\epsilon}$ and $u^{\epsilon^{\prime}}$ are based on. Since $u(t) \in B(0, \rho_{\max})$, we obtain
    \begin{align}
        \int_{0}^{1} \Vert u^{\epsilon}(t) - u^{\epsilon^{\prime}}(t) \Vert_2 dt \leq \sup_{u \in B(0, \rho_{\max})} \Vert v - u \Vert_2 (\epsilon - \epsilon^{\prime}).
    \end{align}
    Similarly, if $\epsilon < \epsilon^{\prime}$ we have
    \begin{align}
        \int_{0}^{1} \Vert u^{\epsilon}(t) - u^{\epsilon^{\prime}}(t) \Vert_2 dt \leq \sup_{u \in B(0, \rho_{\max})} \Vert v - u \Vert_2 (\epsilon^\prime - \epsilon).
    \end{align}
    Put these two cases together and substitute into \eqref{eq: lipschitz_perturbation_1} to get
    \begin{align}
        \forall t \in [0, 1]\;\; \Vert x_1^{\epsilon}(t) - x_1^{\epsilon^{\prime}}(t) \Vert_2 &\leq L^{\prime} | \epsilon - \epsilon^{\prime}|,
    \end{align}
    where $L^{\prime} \triangleq L \sup_{u \in B(0, \rho_{\max})} \Vert v - u \Vert_2 \leq 2L\rho_{\max}< \infty$.
\end{proof}
\end{lemma}

\begin{lemma}
    \label{lemma: state_variation_1}
    Let $u^{\epsilon}$ and $x_1^{\epsilon}$ be as in Lemma \ref{lemma: lipschitz_continuity}. Then,
    \begin{multline}
        \lim_{\epsilon \rightarrow 0^{+}} \frac{1}{\epsilon} \int_{\tau - \epsilon}^{\tau} \left\{ f(x_1^{\epsilon}(t), u^{\epsilon}(t)) - f(x_1(t), u(t))\right\} dt \\
        = f(x_1(\tau), v) - f(x_1(\tau), u(\tau)),
    \end{multline}
    where $x_1$ denotes the solution under the nominal control $u \in U$.
\begin{proof}
    We will show that the difference norm
    \begin{multline}
        \label{eq: difference_norm_variation}
        \bigg\Vert \frac{1}{\epsilon} \int_{\tau - \epsilon}^{\tau} \left\{ f(x_1^{\epsilon}(t), u^{\epsilon}(t)) - f(x_1(t), u(t))\right\} dt \\
        - f(x_1(\tau), v) + f(x_1(\tau), u(\tau)) \bigg\Vert_2
    \end{multline}
    converges to $0$ as $\epsilon \rightarrow 0^{+}$. Indeed, \eqref{eq: difference_norm_variation} becomes
    \begin{align}
        &\begin{multlined}
            \bigg\Vert \frac{1}{\epsilon} \int_{\tau - \epsilon}^{\tau} \{ f(x_1^{\epsilon}(t), u^{\epsilon}(t)) - f(x_1(t), u(t)) \\
            - f(x_1(\tau), v) + f(x_1(\tau), u(\tau))\} dt \bigg\Vert_2
        \end{multlined} \\
        &\begin{multlined}
            \leq \frac{1}{\epsilon} \int_{\tau - \epsilon}^{\tau} \Vert f(x_1^{\epsilon}(t), u^{\epsilon}(t)) - f(x_1(t), u(t)) \\
            - f(x_1(\tau), v) + f(x_1(\tau), u(\tau)) \Vert_2 dt 
        \end{multlined} \\
        &\begin{multlined}
            \label{eq: difference_norm_variation_2}
            \leq \frac{1}{\epsilon} \int_{\tau - \epsilon}^{\tau} \big\{\Vert f(x_1^{\epsilon}(t), u^{\epsilon}(t)) - f(x_1(\tau), v) \Vert_2 \\
            + \Vert f(x_1(t), u(t)) - f(x_1(\tau), u(\tau)) \Vert_2\big\} dt.
        \end{multlined}
    \end{align}
    We used the triangle inequality in \eqref{eq: difference_norm_variation_2}. Now, making use of the fact that $\forall t \in (\tau - \epsilon, \tau]\;\; u^\epsilon(t) = v $, Assumption \ref{assum: function_f_bound} yields
    \begin{align}
        \Vert f(&x_1^{\epsilon}(t), u^{\epsilon}(t)) - f(x_1(\tau), v) \Vert_2 \nonumber \\
        &\leq K_1 \Vert x_1^{\epsilon}(t) - x_1(\tau) \Vert_2 \\
        &= K_1 \Vert x_1^{\epsilon}(t) - x_1(t) + x_1(t) - x_1(\tau) \Vert_2 \\
        &\leq K_1 \Vert x_1^{\epsilon}(t) - x_1(t) \Vert_2 + K_1\Vert x_1(t) - x_1(\tau) \Vert_2 \\
        &\leq K_1L^{\prime} \epsilon + K_1\Vert x_1(t) - x_1(\tau) \Vert_2 \label{eq: difference_norm_lemma_11_applied}
    \end{align}
    for any $t \in (\tau - \epsilon, \tau]$, where we applied Lemma \ref{lemma: lipschitz_continuity} in \eqref{eq: difference_norm_lemma_11_applied} with $\epsilon^{\prime} = 0$.
    Similarly,
    \begin{multline}
        \Vert f(x_1(t), u(t)) - f(x_1(\tau), u(\tau)) \Vert_2 \\
        \leq K_1 \Vert x_1(t) - x_1(\tau) \Vert_2 + K_1 \Vert u(t) - u(\tau) \Vert_2.
    \end{multline}
    Therefore, \eqref{eq: difference_norm_variation} is upper-bounded by
    \begin{align}
        &\begin{multlined}
            \frac{1}{\epsilon} \int_{\tau-\epsilon}^{\tau} \{K_1 L^{\prime}\epsilon + 2K_1  \Vert x_1(t) - x_1(\tau) \Vert_2 \\
            + K_1 \Vert u(t) - u(\tau) \Vert_2 \} dt
        \end{multlined} \\
        &\begin{multlined}
            \leq K_1 L^{\prime}\epsilon + 2K_1 \sup_{t \in [\tau-\epsilon, \tau]} \Vert x_1(t) - x_1(\tau)\Vert_2 \\
            + K_1 \sup_{t \in [\tau-\epsilon, \tau]} \Vert u(t) - u(\tau)\Vert_2,
        \end{multlined}
    \end{align}
    which converges to $0$ as $\epsilon \rightarrow 0^{+}$, since
    \begin{align}
        0 \leq \sup_{t \in [\tau-\epsilon, \tau]} \Vert x_1(t) - x_1(\tau)\Vert_2 \longrightarrow 0
    \end{align}
    and
    \begin{align}
        0 \leq \sup_{t \in [\tau-\epsilon, \tau]} \Vert u(t) - u(\tau)\Vert_2 \longrightarrow 0
    \end{align}
    as $\epsilon \rightarrow 0^{+}$.
\end{proof}
\end{lemma}

\begin{lemma}
    \label{lemma: state_variation_2}
    Let $x_1^{\epsilon}$ and $x_1$ be as in Lemma \ref{lemma: state_variation_1}. Let $\Psi_1(t)$ be the right derivative of $x_1^{\epsilon}(t)$ with respect to $\epsilon$ evaluated at $\epsilon = 0$. That is,
    \begin{align}
        \Psi_1(\tau) = \frac{\partial_+}{\partial \epsilon} x_1^{\epsilon}(\tau) \bigg\vert_{\epsilon=0} \triangleq \lim_{\epsilon \rightarrow 0^+} \frac{x_1^{\epsilon}(\tau) - x_1(\tau)}{\epsilon}.
    \end{align}
    Then, we have
    \begin{align}
         \Psi_1(\tau) = f(x_1(\tau), v) - f(x_1(\tau), u(\tau)).
    \end{align}
\begin{proof}
    Let us express both $x_1^{\epsilon}(\tau)$ and $x_1(\tau)$ in the integral form:
    \begin{align}
        x_1^{\epsilon}(\tau) &= x_0 + \int_{0}^\tau f(x_1^{\epsilon}(t), u^\epsilon(t))dt \\
        x_1(\tau) &= x_0 + \int_{0}^\tau f(x_1(t), u(t))dt.
    \end{align}
    Note that $u^{\epsilon}(t) = u(t)$ and $x_1^{\epsilon}(t) = x_1(t)$ for $t \in [0, \tau - \epsilon]$, since no perturbation is applied to the system until $t > \tau - \epsilon$.
    Therefore,
    \begin{multline}
        x_1^{\epsilon}(\tau) - x_1(\tau) \\
        = \int_{\tau - \epsilon}^\tau \left\{ f(x_1^{\epsilon}(t), u^{\epsilon}(t)) - f(x_1(t), u(t))\right\} dt.
    \end{multline}
    Making use of Lemma \ref{lemma: state_variation_1}, we conclude that
    \begin{align}
        \lim_{\epsilon \rightarrow 0^+} \frac{x_1^{\epsilon}(\tau) - x_1(\tau)}{\epsilon} = f(x_1(\tau), v) - f(x_1(\tau), u(\tau)).
    \end{align}
\end{proof}
\end{lemma}

\begin{lemma}
    \label{lemma: bounded_dynamics_jacobian}
    Suppose that Assumption \ref{assum: function_f_bound} is satisfied. Then, we have the following bound on the (matrix) norm of the Jacobian of function $f$:
    \begin{align}
        \left\Vert \frac{\partial}{\partial x} f(x^{\prime}, u^{\prime}) \right\Vert_2 \leq K_1,
    \end{align}
    for any $x^{\prime} \in \mathbb{R}^{n_x}$ and $u^{\prime} \in B(0, \rho_{\max} )$.
\begin{proof}
    Let $u^{\prime\prime}$ be equal to $u^{\prime}$ in Assumption \ref{assum: function_f_bound}. Then, we have
    \begin{align}
        \label{eq: lipschitz_bound_f}
        \Vert f(x^{\prime}, u^{\prime}) - f(x^{\prime\prime}, u^{\prime})\Vert_2 \leq K_1 \Vert x^{\prime} - x^{\prime\prime}\Vert_2
    \end{align}
    for any $x^{\prime}, x^{\prime\prime} \in \mathbb{R}^{n_x}$ and any $u^{\prime} \in B(0, \rho_{\max})$.
    Now, let us define some non-zero scalar $t$ and some unit vector $v \in \mathbb{R}^{n_x}$ that satisfies $\Vert v \Vert_2 = 1$, and let
    \begin{align}
        \label{eq: reparameterization_x_double_prime}
        x^{\prime\prime} \triangleq x^{\prime} + t v.
    \end{align}
    Substituting \eqref{eq: reparameterization_x_double_prime} into \eqref{eq: lipschitz_bound_f}, we get
    \begin{align}
       \frac{\Vert f(x^{\prime} + tv, u^{\prime}) - f(x^{\prime}, u^{\prime})\Vert_2}{|t|} \leq K_1.
    \end{align}
    Note that this holds for any $t$ and $v$ as defined above.
    From multivariate calculus, on the other hand, the directional derivative of $f$ with respect to $v$ is given by
   \begin{align}
       \lim_{t\rightarrow0} \frac{f(x^{\prime} + tv, u^{\prime}) - f(x^{\prime}, u^{\prime})}{t} = \left(\frac{\partial}{\partial x} f(x^{\prime}, u^{\prime})\right) v.
   \end{align}
   Therefore, we have
   \begin{align}
       \left\Vert \left(\frac{\partial}{\partial x} f(x^{\prime}, u^{\prime})\right) v \right\Vert_2 \leq K_1,
   \end{align}
    which implies
    \begin{align}
        \left\Vert \frac{\partial}{\partial x} f(x^{\prime}, u^{\prime}) \right\Vert_2 \triangleq \sup_{\Vert v\Vert_2 = 1} \left\Vert \left(\frac{\partial}{\partial x} f(x^{\prime}, u^{\prime})\right) v \right\Vert_2 \leq K_1.
    \end{align}
\end{proof}
\end{lemma}
   
\begin{lemma}
    \label{lemma: state_variation_3}
     Let $x_1^{\epsilon}$ and $x_1$ be as in Lemma \ref{lemma: state_variation_1}. Then, $\Psi_1(t) = \frac{\partial_+}{\partial\epsilon} x_1^{\epsilon}(t)\big\vert_{\epsilon=0}$ uniquely exists for $t \in [\tau, 1]$ and follows the ODE:
    \begin{align}
        \dot{\Psi}_1(t) = \frac{\partial}{\partial x_1} f(x_1(t), u(t)) \Psi_1(t),
    \end{align}
    with the initial condition $\Psi_1(\tau)$ given by Lemma \ref{lemma: state_variation_2}.
\begin{proof}
    Taking some $a \in (\tau, 1]$, let us express $x_1^{\epsilon}(a)$ and $x_1(a)$ in the integral form:
    \begin{align}
        x_1^{\epsilon}(a) &= x_1^{\epsilon}(\tau) + \int_{\tau}^a f(x_1^{\epsilon}(t), u^\epsilon(t))dt \\
        x_1(a) &= x_1(\tau) + \int_{\tau}^a f(x_1(t), u(t))dt.
    \end{align}
    Thus, we have
    \begin{multline}
        \label{eq: variational_equation_1}
        \Psi_1(a) = \Psi_1(\tau) \\
        + \lim_{\epsilon \rightarrow 0{+}} \int_{\tau}^a \frac{1}{\epsilon}\{f(x_1^{\epsilon}(t), u(t)) - f(x_1(t), u(t))\} dt,
    \end{multline}
    where we used $\forall t \in [\tau, a]\;\; u^{\epsilon}(t) = u(t)$.
    We will take a measure-theoretic approach to prove that the order of the limit and the integration can be switched in \eqref{eq: variational_equation_1}. First, think of the integral as a Lebesgue integral on the measure space $([\tau, a], \mathcal{B}([\tau, a]), \lambda)$, where $\mathcal{B}([\tau, a])$ is the Borel $\sigma$-algebra on $[\tau, a]$ and $\lambda$ is the Lebesgue measure. Furthermore, consider the integrand as a function from $[\tau, a]$ into the Banach space $(\mathbb{R}^{n_x}, \Vert\cdot\Vert_2)$, i.e. the Euclidean space endowed with the $\ell^2$ norm. By the piecewise continuity of $u^\epsilon, u$ and the continuity of $x_1^\epsilon$, $x_1$, and $f$, the integrand is a piecewise continuous function with respect to $t$, which is $\mathcal{\lambda}$-measurable. In fact it is also Bochner-integrable, since for $t \in [\tau, a]$ we have the constant bound:
    \begin{align}
        \frac{1}{\epsilon} \Vert f(x_1^{\epsilon}(t), u(t)) &- f(x_1(t), u(t)) \Vert_2 \nonumber \\ &\leq \frac{1}{\epsilon} K_1 \Vert x_1^\epsilon(t) - x_1(t) \Vert_2 \\
        &\leq K_1 L^{\prime} 
    \end{align}
    by Assumption \ref{assum: function_f_bound} and Lemma \ref{lemma: lipschitz_continuity}.
    Furthermore, the chain rule gives
    \begin{align}
        \lim_{\epsilon \rightarrow 0^{+}} \frac{1}{\epsilon}(f(x_1^{\epsilon}(t), u(t)) - f(x_1(t), u(t))) \nonumber \\
        =\frac{\partial}{\partial x_1} f(x_1(t), u(t)) \Psi_1(t),
    \end{align}
    assuming that $\Psi_1(t)$ exists for $t \in [\tau, a]$.
    Therefore, by the Bochner-integral version of the dominated convergence theorem (Theorem 3 in \cite{diestel1977vector}, Chapter II), we obtain
    \begin{align}
        \Psi_1(a) = \Psi_1(\tau) + \int_{\tau}^a \frac{\partial}{\partial x_1} f(x_1(t), u(t)) \Psi_1(t) dt.
    \end{align}
    This is equivalent to the ordinary differential equation:
    \begin{align}
        \label{eq: variational_equation_2}
        \dot{\Psi}_1(t) = \frac{\partial}{\partial x_1} f(x_1(t), u(t)) \Psi_1(t).
    \end{align}
    It remains to show that the solution $\Psi_1(t)$ that satisfies \eqref{eq: variational_equation_2} does exist and is unique. First, let $\Psi_1^{\prime}$ and $\Psi_1^{\prime\prime}$ be two systems that follow \eqref{eq: variational_equation_2} and share the same initial condition $\Psi_1(\tau)$. Then, by Lemma \ref{lemma: bounded_dynamics_jacobian} we have
    \begin{align}
        \Vert \dot{\Psi}_1^{\prime}(t) &- \dot{\Psi}_1^{\prime\prime}(t) \Vert_2 \nonumber \\
        &= \left\Vert \frac{\partial}{\partial x_1} f(x_1(t), u(t)) (\Psi_1^{\prime}(t) - \Psi_1^{\prime\prime}(t)) \right\Vert_2 \\
        &\leq \left\Vert \frac{\partial}{\partial x_1} f(x_1(t), u(t))\right\Vert_2 \cdot \left\Vert \Psi_1^{\prime}(t) - \Psi_1^{\prime\prime}(t) \right\Vert_2 \\
        &\leq K_1 \left\Vert \Psi_1^{\prime}(t) - \Psi_1^{\prime\prime}(t) \right\Vert_2.
    \end{align}
    Existence follows from this inequality in conjunction with the Picard Lemma (Lemma 5.6.3 in \cite{elijah1997optimization}). To show the uniqueness, apply the Bellman-Gronwall Lemma (Lemma 5.6.4 in \cite{elijah1997optimization}) to the following integral inequality:
    \begin{align}
        &\forall a \in [\tau, 1]\;\; \nonumber  \\
        &\Vert \Psi_1^{\prime}(a) - \Psi_1^{\prime\prime}(a) \Vert_2 \leq K_1 \int_{\tau}^{a} \Vert \Psi_1^{\prime}(t) - \Psi_1^{\prime\prime}(t) \Vert_2 dt.
    \end{align}
\end{proof} 
\end{lemma}

\begin{prop}[Variational Equation]
    \label{prop: state_variation}
    Let $u \in U$ and $x_1$ be the nominal control and the resulting state trajectory. Let $x_1^{\epsilon}$ be the perturbed state induced by the perturbed control $u^\epsilon$ of the form \eqref{eq: perturbed_control}. Propagating $x_1^\epsilon(t)$ through the hybrid dynamics, we get a series of modes $x^{\epsilon}_2, \dots, x^{\epsilon}_T$ that constitutes the entire trajectory $x^{\epsilon}(t)$ for $t \in [\tau, T]$. Define the state variation $\Psi(t)$ for $t \in [\tau, T]$ by
    \begin{align}
        \Psi(t) = \frac{\partial_+}{\partial\epsilon} x^{\epsilon}(t)\bigg|_{\epsilon=0} \triangleq \lim_{\epsilon\rightarrow 0^{+}} \frac{x^{\epsilon}(t) - x(t)}{\epsilon}.
    \end{align}
    Then, $\Psi(t)$ exists for $t \in [\tau, T]$ and follows the hybrid system with time-driven switching:
    \begin{align}
        \label{eq: ih_prop_14_1}
        \Psi(t) = \begin{cases}
            \Psi_1(t) & \forall t \in [\tau, 1) \\
            \Psi_i(t) & \forall t \in [i-1, i)\; \forall i \in \{2,\dots,T \}
            \end{cases}
    \end{align}
    with $\Psi(T) = \frac{\partial}{\partial x_T} g(x_T(T), y_T) \Psi_{T}(T)$,
    where $\Psi_1$ is defined on $[\tau, 1]$ as in Lemma \ref{lemma: state_variation_3}, and $\Psi_i$ for $i \geq 2$ is defined on $[i-1, i]$ with
    \begin{align}
        \Psi_i(i-1) &= \frac{\partial}{\partial x_{i-1}}g(x_{i-1}(i-1), y_{i-1}) \Psi_{i-1}(i-1) \label{eq: ih_prop_14_2}\\
        \dot{\Psi}_i(t) &= \frac{\partial}{\partial x_i} f(x_i(t), u(t)) \Psi_i(t) \;\;\forall t \in [i-1, i]. \label{eq: ih_prop_14_3}
    \end{align}
    
\begin{proof}
    The case for $t \in [\tau, 1)$ follows from Lemma \ref{lemma: state_variation_3}, since in this case we have
    \begin{align}
        \Psi(t) &\triangleq \lim_{\epsilon\rightarrow 0^{+}} \frac{x^{\epsilon}(t) - x(t)}{\epsilon} \\
        &= \lim_{\epsilon\rightarrow 0^{+}} \frac{x^{\epsilon}_1(t) - x_1(t)}{\epsilon} \\
        &= \Psi_1(t).
    \end{align}
    At $t = 1$, we obtain
    \begin{align}  
        \Psi(1) &\triangleq \lim_{\epsilon\rightarrow 0^{+}} \frac{x^{\epsilon}(1) - x(1)}{\epsilon} \\
        &= \lim_{\epsilon\rightarrow 0^{+}} \frac{x^{\epsilon}_2(1) - x_2(1)}{\epsilon} \\
        &=  \lim_{\epsilon\rightarrow 0^{+}} \frac{g(x_1^{\epsilon}(1), y_1) - g(x_1(1), y_1)}{\epsilon} \\
        &= \frac{\partial}{\partial x} g(x_1(1), y_1) \Psi_1(1)
    \end{align}
    by \eqref{eq: x_jump_definition} and the chain rule. Let us define $\Psi_2$ by $\Psi_2(t) \triangleq \lim_{\epsilon\rightarrow 0^{+}} (x^{\epsilon}_2(t) - x_2(t))/\epsilon$. Similarly to the proof of Lemma \ref{lemma: state_variation_3}, one can show that $\Psi_2$ follows the integral equation:
    \begin{align}
        \Psi_2(a) = \Psi_2(1) + \int_{1}^a \frac{\partial}{\partial x_2} f(x_2(t), u(t)) \Psi_2(t) dt 
    \end{align}
    with $\Psi_2(1) = \Psi(1)$, and that $\Psi_2(t)$ that satisfies \eqref{eq: ih_prop_14_2} and \eqref{eq: ih_prop_14_3} uniquely exists. This proves the case for $t \in [1, 2)$. Proceeding by induction completes the proof.
\end{proof}
\end{prop}

So far we have focused entirely on the right derivative $\frac{\partial_+}{\partial \epsilon}x^\epsilon(t)$ evaluated at $\epsilon = 0$. The next proposition shows that $x_1^{\epsilon}$ is in fact right differentiable with respect to $\epsilon$ at all $\epsilon \in [0, \tau)$.
\begin{prop}[Right Differentiability of State Perturbation]
    \label{prop: right_differentiability_state}
    Let $x^{\epsilon}(t)$ be the perturbed state trajectory under the perturbed control $u^{\epsilon}$ defined by \eqref{eq: perturbed_control}. Let $\Psi^{\epsilon}(t)$ denote the right derivative $\frac{\partial_+}{\partial \epsilon}x^{\epsilon}(t)$ evaluated at a particular $\epsilon \in [0,\tau)$. (Note that when $\epsilon = 0$ we have $\Psi^\epsilon$ = $\Psi$.) Then, $\Psi^{\epsilon}(t)$ exists for $t \in [\tau - \epsilon, T]$ and follows the hybrid system with time-driven switching:
    \begin{align}
        \Psi^{\epsilon}(t) = 
        \begin{cases}
            \Psi_1^{\epsilon}(t) & \forall t \in [\tau - \epsilon, 1) \\
            \Psi_i^{\epsilon}(t) & \forall t \in [i-1, i)\; \forall i \in \{2,\dots,T \}
        \end{cases}
    \end{align}
    with $\Psi^{\epsilon}(T) = \frac{\partial}{\partial x^{\epsilon}_T} g(x^{\epsilon}_T(T), y_T) \Psi^{\epsilon}_{T}(T)$.
    $\Psi_1^{\epsilon}$ is defined on $[\tau - \epsilon, 1]$, where
    \begin{align}
        \Psi_1^{\epsilon}(\tau-\epsilon) = f(x_1^{\epsilon}(\tau-\epsilon), v) - f(x_1^{\epsilon}(\tau - \epsilon), u^{\epsilon} (\tau - \epsilon))
    \end{align}
    and
    \begin{align}
        \dot{\Psi}_1^{\epsilon}(t) = \frac{\partial}{\partial x_1^{\epsilon}} f(x_1^{\epsilon}(t), u^{\epsilon}(t)) \Psi_1^{\epsilon}(t) \;\;\forall t \in [\tau-\epsilon, 1].
    \end{align}
    $\Psi_i^{\epsilon}$ for $i \geq 2$ is defined on $[i-1, i]$ with
    \begin{align}
        \Psi_i^{\epsilon}(i-1) &= \frac{\partial}{\partial x_{i-1}^{\epsilon}}g(x_{i-1}^{\epsilon}(i-1), y_{i-1}) \Psi_{i-1}^{\epsilon}(i-1) \\
        \dot{\Psi}_i^{\epsilon}(t) &= \frac{\partial}{\partial x_i^{\epsilon}} f(x_i^{\epsilon}(t), u^{\epsilon}(t)) \Psi_i^{\epsilon}(t) \;\;\forall t \in [i-1, i].
    \end{align}
    
\begin{proof}
    The proposition can be proven by considering the perturbed control $u^{\epsilon}$ as the new nominal control, and defining a further perturbation based on this nominal control. Formally, define $\tilde{u}^{\epsilon^{\prime}}$ by
    \begin{align}
        \label{eq: perturbed_control_further}
        \tilde{u}^{\epsilon^\prime}(t) \triangleq \begin{cases}
            v & \text{if}\;\; t \in (\tau - \epsilon - \epsilon^{\prime}, \tau - \epsilon] \\
            u^{\epsilon}(t) & \text{otherwise}
        \end{cases}
    \end{align}
    with $\epsilon^{\prime} \in [0, \tau - \epsilon)$, where $(\tau, v)$ is the same pair of values as for $u^{\epsilon}(t)$. Since $u^{\epsilon}(t)$ is left continuous in $t$ at $t = \tau - \epsilon$, $\tilde{u}^{\epsilon^{\prime}}$ is a valid perturbed control of the form \eqref{eq: perturbed_control} with $u^{\epsilon}$ being the nominal control. Here $\epsilon$ is considered as fixed and the parameters defining this new perturbation are $v$, $\tau-\epsilon$, and $\epsilon^{\prime}$. This perturbation yields the new perturbed state trajectory $\tilde{x}^{\epsilon^{\prime}}$ based on the nominal trajectory $x^{\epsilon}$. Note that when $\epsilon^{\prime} = 0$ we have $\tilde{u}^{\epsilon^{\prime}} = u^{\epsilon}$ and $\tilde{x}^{\epsilon^{\prime}} = x^{\epsilon}$. We can define the new state variation:
    \begin{align}
        \widetilde{\Psi}(t) = \frac{\partial_+}{\partial\epsilon^{\prime}}\tilde{x}^{\epsilon^{\prime}}(t) \bigg|_{\epsilon^{\prime}=0} \triangleq \lim_{\epsilon^{\prime}\rightarrow 0^{+}} \frac{\tilde{x}^{\epsilon^{\prime}}(t) - x^{\epsilon}(t)}{\epsilon^{\prime}}.
    \end{align}
    Applying Proposition \ref{prop: state_variation} to this new setting, we find that $\widetilde{\Psi}(t)$ exists for $t \in [\tau - \epsilon, T]$ and follows the hybrid system with time-driven switching:
    \begin{align}
        \widetilde{\Psi}(t) = 
        \begin{cases}
            \widetilde{\Psi}_1(t) & \forall t \in [\tau - \epsilon, 1) \\
            \widetilde{\Psi}_i(t) & \forall t \in [i-1, i)\; \forall i \in \{2,\dots,T \}
        \end{cases}
    \end{align}
    with $\widetilde{\Psi}(T) = \frac{\partial}{\partial x^{\epsilon}_T} g(x^{\epsilon}_T(T), y_T) \widetilde{\Psi}_{T}(T)$.
    $\widetilde{\Psi}_1$ is defined on $[\tau - \epsilon, 1]$, where
    \begin{align}
        \widetilde{\Psi}_1(\tau-\epsilon) = f(x_1^{\epsilon}(\tau-\epsilon), v) - f(x_1^{\epsilon}(\tau - \epsilon), u^{\epsilon} (\tau - \epsilon))
    \end{align}
    and
    \begin{align}
        \dot{\widetilde{\Psi}}_1(t) = \frac{\partial}{\partial x_1^{\epsilon}} f(x_1^{\epsilon}(t), u^{\epsilon}(t)) \widetilde{\Psi}_1(t) \;\;\forall t \in [\tau-\epsilon, 1].
    \end{align}
    $\widetilde{\Psi}_i$ for $i \geq 2$ is defined on $[i-1, i]$ with
    \begin{align}
        \widetilde{\Psi}_i(i-1) &= \frac{\partial}{\partial x_{i-1}^{\epsilon}}g(x_{i-1}^{\epsilon}(i-1), y_{i-1}) \widetilde{\Psi}_{i-1}(i-1) \\
        \dot{\widetilde{\Psi}}_i(t) &= \frac{\partial}{\partial x_i^{\epsilon}} f(x_i^{\epsilon}(t), u^{\epsilon}(t)) \widetilde{\Psi}_i^{\epsilon}(t) \;\;\forall t \in [i-1, i].
    \end{align}
    
    On the other hand, notice that the new perturbed control $\tilde{u}^{\epsilon^{\prime}}$ is actually equivalent to the perturbed control $u^{\epsilon + \epsilon^{\prime}}$ that is based on the original nominal control $u$. Namely,
    \begin{align}
        \tilde{u}^{\epsilon^{\prime}}(t) = u^{\epsilon + \epsilon^{\prime}}(t) \triangleq \begin{cases}
        v & \text{if}\;\; t \in (\tau - (\epsilon + \epsilon^{\prime}), \tau] \\
        u(t) & \text{otherwise}.
        \end{cases}
    \end{align}
    Consequently, the new perturbed state $\tilde{x}^{\epsilon^{\prime}}$ is equal to $x^{\epsilon + \epsilon^{\prime}}$, and thus
    \begin{align}
         \widetilde{\Psi}(t) &= \lim_{\epsilon^{\prime}\rightarrow 0^{+}} \frac{\tilde{x}^{\epsilon^{\prime}}(t) - x^{\epsilon}(t)}{\epsilon^{\prime}} \\
         &= \lim_{\epsilon^{\prime}\rightarrow 0^{+}} \frac{x^{\epsilon + \epsilon^{\prime}}(t) - x^{\epsilon}(t)}{\epsilon^{\prime}} \\
         &= \Psi^{\epsilon}(t).
    \end{align}
    This completes the proof.
\end{proof}
\end{prop}

\begin{lemma}
     \label{lemma: bound_2}
     Let $\Psi_1^{\epsilon},\dots, \Psi_T^{\epsilon}$ be as given by Proposition \ref{prop: right_differentiability_state}. Then, for $\Psi_1^{\epsilon}$ the following holds.
    \begin{align}
        \Vert \Psi_1^{\epsilon}(t) \Vert_2 &\leq 2 K_1 \rho_{\max} e^{K_1} &\forall t \in [\tau - \epsilon,1]
    \end{align}
    Similarly, for all $i \in \{2,\dots,T\}$ we have
    \begin{align}
        \Vert \Psi^{\epsilon}_i(t) \Vert_2 &\leq \Vert \Psi_i^{\epsilon}(i-1) \Vert_2 e^{K_1} &\forall t \in [i-1, i].
    \end{align}
\begin{proof}
    We begin with the integral equation:
    \begin{multline}
        \Psi_1^{\epsilon}(a) = \Psi_1^{\epsilon}(\tau - \epsilon) \\ 
        + \int_{\tau-\epsilon}^{a} \frac{\partial}{\partial x_1^{\epsilon}} f(x_1^{\epsilon}(t), u^{\epsilon}(t))\Psi_1^{\epsilon}(t) dt.
    \end{multline}
    Therefore,
    \begin{multline}
        \Vert \Psi_1^{\epsilon}(a)\Vert_2 \leq \Vert \Psi_1^{\epsilon}(\tau - \epsilon) \Vert_2 \\ 
        + \int_{\tau - \epsilon}^{a} \left\Vert \frac{\partial}{\partial x_1^{\epsilon}}f(x_1^{\epsilon}(t), u^{\epsilon}(t)) \right\Vert_2 \cdot \Vert \Psi_1^{\epsilon}(t) \Vert_2 dt.
    \end{multline}
    Using 
    \begin{align}
        \Psi_1^{\epsilon}(\tau - \epsilon) = f(x_1^{\epsilon}(\tau - \epsilon), v) - f(x_1^{\epsilon}(\tau - \epsilon), u^{\epsilon}(\tau - \epsilon)),
    \end{align}
    Assumption \ref{assum: function_f_bound}, and Lemma \ref{lemma: bounded_dynamics_jacobian}, we get
    \begin{align}
        \Vert \Psi_1^{\epsilon}(a)\Vert_2 &\leq K_1 \Vert v - u^{\epsilon}(\tau - \epsilon)\Vert_2
        + K_1 \int_{\tau - \epsilon}^{a} \Vert \Psi_1^{\epsilon}(t) \Vert_2 dt \\
        &\begin{multlined}
            \leq K_1 \sup_{u \in B(0, \rho_{\max})} \Vert v - u \Vert_2 \\
            + K_1 \int_{\tau - \epsilon}^{a} \Vert \Psi_1^{\epsilon}(t) \Vert_2 dt
        \end{multlined} \\
        &\leq 2K_1 \rho_{\max} + K_1 \int_{\tau - \epsilon}^{a} \Vert \Psi_1^{\epsilon}(t) \Vert_2 dt
    \end{align}
    for all $a \in [\tau - \epsilon, 1]$. Thus, by the Bellman-Gronwall Lemma (Lemma 5.6.4 in \cite{elijah1997optimization}) it follows that
    \begin{align}
        \Vert \Psi_1^{\epsilon}(t) \Vert_2 &\leq 2 K_1 \rho_{\max} e^{K_1} &\forall t \in [\tau - \epsilon,1].
    \end{align}
    For general $i \geq 2$, apply the Bellman-Gronwall Lemma to the similar integral inequality:
    \begin{align}
        &\forall a \in [i-1, i] \nonumber \\
        &\Vert \Psi_i^{\epsilon}(a) \Vert_2 \leq \Vert \Psi_i^{\epsilon}(i-1) \Vert_2 + K_1 \int_{i-1}^{a} \Vert \Psi_i^{\epsilon}(t) \Vert_2 dt
    \end{align}
    to get the result.
\end{proof}
\end{lemma}

\begin{prop}[Bounded State Variation]
    \label{prop: bound_2}
    Given $u^{\epsilon}$ and $(y_1, \dots, y_T)$, $\Psi^{\epsilon}$ defined in Proposition \ref{prop: right_differentiability_state} has the following bound:
    \begin{align}
        &\forall \epsilon \in [0, \tau)\; \forall i \in \{2, \dots, T\}\; \forall t \in [i - 1, i]\;\; \nonumber \\
        &\Vert \Psi_i^{\epsilon}(t) \Vert_2 \leq \sum_{(j_1,\dots,j_{i-1}) \in \mathcal{L}_i} \beta_i^{(j_1,\dots,j_{i-1})} (x_0) \prod_{m=1}^{i-1} \Vert y_m \Vert_2^{j_m},
    \end{align}
    where $\mathcal{L}_i$ is a finite set of sequences of non-negative integers of length $i-1$, and $\beta_i^{(j_1,\dots,j_{i-1})} (x_0)$ is a finite positive constant that depends on $x_0$ and $(j_1, \dots, j_{i-1})$ but not on $\epsilon$, $u^{\epsilon}$, or $(y_1,\dots,y_T)$.

\begin{proof}
    The proof of this proposition is similar to that of Proposition \ref{prop: bound_1}. Take any $\epsilon \in [0, \tau)$. For $i = 2$, we have $\forall t \in [1, 2]$
    \begin{align}
        \Vert \Psi_2^{\epsilon}(t) \Vert_2 &\leq \Vert \Psi_2^{\epsilon}(1) \Vert_2 e^{K_1} \\
        &\leq \left\Vert \frac{\partial}{\partial x_1^{\epsilon}} g(x_1^{\epsilon}(1), y_1) \Psi_1^{\epsilon}(1)\right\Vert_2 e^{K_1} \\
        &\leq \left\Vert \frac{\partial}{\partial x_1^{\epsilon}} g(x_1^{\epsilon}(1), y_1)\right\Vert_2 \cdot \left\Vert \Psi_1^{\epsilon}(1)\right\Vert_2 e^{K_1}\\
        &\begin{multlined}
            \label{eq: bound_3}
            \leq 2 K_1 \rho_{\max} e^{2K_1} \bigg\{K_2 + K_4 \Vert y_1 \Vert_2^{L_2} + \\ \left(K_3 + K_5 \Vert y_1 \Vert_2^{L_2}\right) \Vert x_1^{\epsilon}(1) \Vert_2^{L_1} \bigg\}
        \end{multlined}
    \end{align}
    by Assumption \ref{assum: function_g_bound}, Proposition \ref{prop: right_differentiability_state}, and Lemma \ref{lemma: bound_2}. Using Proposition \ref{prop: bound_1}, we can bound $x_1^{\epsilon}(1)$ by
    \begin{align}
        \label{eq: bound_4}
        \Vert x_1^{\epsilon}(1) \Vert_2 \leq \sum_{(j_1)\in\mathcal{K}_2} \alpha_2^{(j_1)}(x_0) \Vert y_1 \Vert_2^{j_1}.
    \end{align}
    Substituting \eqref{eq: bound_4} into \eqref{eq: bound_3} and using the multinomial theorem, one can verify that
    \begin{align}
        \forall t \in [1, 2]\;\; \Vert \Psi_2^{\epsilon}(t) \Vert_2 \leq \sum_{(j_1) \in \mathcal{L}_1} \beta_1^{(j_1)}(x_0) \Vert y_1 \Vert_2^{j_1}
    \end{align}
    for some finite set $\mathcal{L}_1$ and finite $\beta_1^{(j_1)} (x_0)$.
    
    Next, suppose that the claim holds for some $i \leq 2$. That is,
    \begin{align}
        &\forall t \in [i-1, i] \nonumber \\
        &\Vert \Psi_i^{\epsilon}(t) \Vert_2 \leq \sum_{(j_1,\dots,j_{i-1}) \in \mathcal{L}_i} \beta_i^{(j_1,\dots,j_{i-1})} (x_0) \prod_{m=1}^{i-1} \Vert y_m \Vert_2^{j_m},
    \end{align}
    where $\mathcal{L}_i$ and $\beta_i^{(j_1,\dots,j_{i-1})}(x_0)$ are as defined in the statement of the proposition. Similar to the case for $i = 2$, we have $\forall t \in [i, i+1]$ 
    \begin{multline}
        \label{eq: bound_5}
        \Vert \Psi_{i+1}^{\epsilon}(t) \Vert_2 \\ \leq e^{K_1} \left(\sum_{(j_1,\dots,j_{i-1}) \in \mathcal{L}_i} \beta_i^{(j_1,\dots,j_{i-1})} (x_0) \prod_{m=1}^{i-1} \Vert y_m \Vert_2^{j_m}\right) \\
        \times \bigg\{K_2 + K_4 \Vert y_i \Vert_2^{L_2} + \left(K_3 + K_5 \Vert y_i \Vert_2^{L_2}\right) \Vert x_i^{\epsilon}(i) \Vert_2^{L_1} \bigg\}
    \end{multline}
    Proposition \ref{prop: bound_1} gives the following bound:
    \begin{align}
        \label{eq: bound_6}
        \Vert x_i^{\epsilon}(i) \Vert_2 \leq \sum_{(j_1, \dots, j_{i-1}) \in \mathcal{K}_i} \alpha_i^{(j_1,\dots,j_{i-1})}(x_0) \prod_{m=1}^{i-1}\Vert y_m \Vert_2^{j_m}.
    \end{align}
    Substituting \eqref{eq: bound_6} into \eqref{eq: bound_5} and using the multinomial theorem, we conclude that
    \begin{align}
        &\forall t \in [i, i+1]\nonumber \\
        &\Vert \Psi_{i+1}^{\epsilon}(t) \Vert_2 \leq \sum_{(j_1,\dots,j_{i}) \in \mathcal{L}_{i+1}} \beta_{i+1}^{(j_1,\dots,j_{i})} (x_0) \prod_{m=1}^{i} \Vert y_m \Vert_2^{j_m}
    \end{align}
    for some finite set $\mathcal{L}_{i+1}$ and finite $\beta_{i+1}^{(j_1,\dots,j_i)}(x_0)$.
    
    Finally, proceeding by mathematical induction over $i \in \{2, \dots, T\}$ completes the proof.
\end{proof}
\end{prop}

\begin{lemma}
    \label{lemma: cost_variation_1}
    Let $u^{\epsilon}$ and $x_1^{\epsilon}$ be as in Lemma \ref{lemma: lipschitz_continuity}. Then, similarly to Lemma \ref{lemma: state_variation_1} we have
    \begin{multline}
        \lim_{\epsilon\rightarrow0^{+}} \frac{1}{\epsilon}\int_{\tau-\epsilon}^{\tau} \left\{c(x_1^{\epsilon}(t), u^{\epsilon}(t)) - c(x_1(t), u(t))\right\}dt \\
        = c(x_1(\tau), v) - c(x_1(\tau), u(\tau)).
    \end{multline}
\begin{proof}
    As $u^{\epsilon}(t) = v\;\; \forall t \in (\tau-\epsilon, \tau]$, we will show that
    \begin{multline}
        \lim_{\epsilon\rightarrow0^{+}} \frac{1}{\epsilon}\int_{\tau-\epsilon}^{\tau} \left\{c(x_1^{\epsilon}(t), v) - c(x_1(t), u(t))\right\}dt \\
        = c(x_1(\tau), v) - c(x_1(\tau), u(\tau)).
    \end{multline}
    By Assumption \ref{assum: cost} and the continuity of $x_1^{\epsilon}(t)$, $c(x_1^{\epsilon}(t), v)$ is continuous with respect to $t$ on $[\tau-\epsilon, \tau]$. Thus, the mean value theorem yields
    \begin{align}
       \frac{1}{\epsilon}\int_{\tau-\epsilon}^{\tau} c(x_1^{\epsilon}(t), v)dt = c(x_1^{\epsilon}(\tilde{t}), v)
    \end{align}
    for some $\tilde{t} \in [\tau-\epsilon, \tau]$.
    From the triangle inequality and Lemma \ref{lemma: lipschitz_continuity} it follows that
    \begin{align}
        \Vert x_1^{\epsilon}(\tilde{t}) - x_1(\tau) \Vert_2 &\leq \Vert x_1^{\epsilon}(\tilde{t}) - x_1(\tilde{t})\Vert_2 + \Vert x_1(\tilde{t}) - x_1(\tau) \Vert_2 \\
        &\leq L^{\prime}\epsilon + \Vert x_1(\tilde{t}) - x_1(\tau) \Vert_2.
    \end{align}
    Therefore, $\lim_{\epsilon\rightarrow 0^{+}} x_1^{\epsilon}(\tilde{t}) = x_1(\tau)$ and 
    \begin{align}
        \label{eq: instant_cost_limit_1}
        \lim_{\epsilon\rightarrow 0^{+}} \frac{1}{\epsilon} \int_{\tau-\epsilon}^{\tau} c(x_1^{\epsilon}(t), v)dt = c(x_1(\tau), v).
    \end{align}
    
    On the other hand, $u(t)$ is continuous on $[\tau - \epsilon, \tau]$ for all sufficiently small $\epsilon$, since $u$ is left continuous at $\tau$ by Definition \ref{def: perturbed_control}. Therefore, $c(x_1(t), u(t))$ is continuous with respect to $t$ on $[\tau-\epsilon, \tau]$ for $\epsilon$ small. The mean value theorem gives
    \begin{align}
         \frac{1}{\epsilon}\int_{\tau-\epsilon}^{\tau} c(x_1(t), u(t))dt = c(x_1(\tilde{t}), u(\tilde{t}))
    \end{align}
    for some $\tilde{t} \in [\tau-\epsilon, \tau]$.
    Taking the limit $\epsilon\rightarrow0^{+}$, the right hand side converges to $c(x_1(\tau), u(\tau))$. Combining this result with \eqref{eq: instant_cost_limit_1}, we conclude that
    \begin{multline}
        \lim_{\epsilon\rightarrow0^{+}} \frac{1}{\epsilon}\int_{\tau-\epsilon}^{\tau} \left\{c(x_1^{\epsilon}(t), v) - c(x_1(t), u(t))\right\}dt \\
        = c(x_1(\tau), v) - c(x_1(\tau), u(\tau)).
    \end{multline}
\end{proof}

\end{lemma}

\begin{lemma}
    \label{lemma: cost_variation_2}
    Given $\epsilon \in [0, \tau)$, $u^{\epsilon}$ and $(y_1,\dots,y_T)$, the right derivative of the instantaneous cost function with respect to $\epsilon$ is given by
    \begin{align}
        &\forall t \in [i-1, i] \nonumber \\
        &\frac{\partial_+}{\partial\epsilon} c(x_i^{\epsilon}(t), u^{\epsilon}(t)) = \frac{\partial}{\partial x_i^{\epsilon}} c(x_i^{\epsilon}(t), u(t))^\mathrm{T} \Psi_i^{\epsilon}(t).
    \end{align}
    for each $i \in \{2,\dots, T\}$. 
    
    For $i = 1$ we have
    \begin{align}
        &\forall t \in (\tau, 1] \nonumber \\
        &\frac{\partial_+}{\partial\epsilon} c(x_1^{\epsilon}(t), u^{\epsilon}(t)) = \frac{\partial}{\partial x_1^{\epsilon}} c(x_1^{\epsilon}(t), u(t))^\mathrm{T} \Psi_1^{\epsilon}(t),
    \end{align}
    
    Similarly, the right derivative of the terminal cost function with respect to $\epsilon$ is given by
    \begin{align}
        \frac{\partial_+}{\partial\epsilon} h(x^{\epsilon}(T)) = \frac{\partial}{\partial x^{\epsilon}}h(x^{\epsilon}(T))^\mathrm{T} \Psi^{\epsilon}(T).
    \end{align}
\begin{proof}
    To prove the claim for the instantaneous cost, note that $u^{\epsilon}(t) = u(t)$ for all $t \in (\tau, T]$ and use the chain rule. The case for the terminal cost also follows from the chain rule.
\end{proof}
\end{lemma}

\begin{prop}[Bounded Cost Variations]
    \label{prop: bounded_cost_variations}
    Given $\epsilon \in [0, \tau)$, $u^{\epsilon}$ and $(y_1,\dots,y_T)$, the right derivative of the instantaneous cost function with respect to $\epsilon$ has the following uniform bound:
    \begin{align}
        &\forall t \in [i-1, i] \nonumber \\
        &\begin{multlined}
            \left\vert\frac{\partial_+}{\partial\epsilon} c(x_i^{\epsilon}(t), u^{\epsilon}(t))\right\vert \\
            \leq \sum_{(j_1,\dots,j_{i-1}) \in \mathcal{L}_i^{\prime}} \beta_i^{\prime(j_1,\dots,j_{i-1})} (x_0) \prod_{m=1}^{i-1} \Vert y_m \Vert_2^{j_m}
        \end{multlined}
    \end{align}
    for each $i \in \{2,\dots, T\}$, where $\mathcal{L}^{\prime}_i$ is a finite set of sequences of nen-negative integers of length $i-1$, and $\beta_i^{\prime(j_1,\dots,j_{i-1})}(x_0)$ is a finite positive constant that depends on $x_0$ and $(j_1,\dots,j_{i-1})$ but not on $\epsilon$, $u^{\epsilon}$, or $(y_1,\dots, y_T)$. 
    
    For $i = 1$ the bound is given by
    \begin{align}
        \forall t \in (\tau, 1]\;\; \left\vert\frac{\partial_+}{\partial\epsilon} c(x_1^{\epsilon}(t),  u^{\epsilon}(t))\right\vert \leq \beta_1^{\prime}(x_0)
    \end{align}
    for some finite positive constant $\beta_1^{\prime}(x_0)$.
    
    Similarly, the right derivative of the terminal cost function with respect to $\epsilon$ has the following bound:
    \begin{multline}
        \left\vert \frac{\partial_+}{\partial\epsilon} h(x^{\epsilon}(T))\right\vert \\
        \leq \sum_{(j_1,\dots,j_{T}) \in \mathcal{L}_{T+1}^{\prime}} \beta_{T+1}^{\prime(j_1,\dots,j_{T})} (x_0) \prod_{m=1}^{T} \Vert y_m \Vert_2^{j_m}
    \end{multline}
    for some finite set $\mathcal{L}^{\prime}_{T+1}$ of sequence of non-negative integers and finite positive constants $\beta_{T+1}^{\prime(j_1,\dots,j_T)}(x_0)$.
\begin{proof}
    The proof of this proposition is similar to that of Proposition \ref{prop: bounded_cost}. Take any $\epsilon \in [0, \tau)$. For $i \in \{2,\dots,T\}$, Assumption \ref{assum: cost} along with Lemma \ref{lemma: cost_variation_2} yields
    \begin{align}
        \left\vert\frac{\partial_+}{\partial\epsilon} c(x_i^{\epsilon}(t), u^{\epsilon}(t))\right\vert &= \left\vert\frac{\partial}{\partial x_i^{\epsilon}} c(x_i^{\epsilon}(t), u(t))^\mathrm{T} \Psi_i^{\epsilon}(t)\right\vert \\
        &\leq \left\Vert\frac{\partial}{\partial x_i^{\epsilon}} c(x_i^{\epsilon}(t), u(t))\right\Vert_2 \cdot \left\Vert\Psi_i^{\epsilon}(t)\right\Vert_2 \\
        &\leq (K_6 + K_7 \Vert x_i^{\epsilon}(t)\Vert_2^{L_3}) \left\Vert \label{eq: bound_7} \Psi_i^{\epsilon}(t)\right\Vert_2.
    \end{align}
    By Propositions \ref{prop: bound_1} and \ref{prop: bound_2}, we have
    \begin{align}
        \Vert x_i^{\epsilon}(t) \Vert_2 &\leq \sum_{(j_1, \dots, j_{i-1}) \in \mathcal{K}_i} \alpha_i^{(j_1,\dots,j_{i-1})}(x_0) \prod_{m=1}^{i-1}\Vert y_m \Vert_2^{j_m} \label{eq: bound_8}\\
        \Vert \Psi_i^{\epsilon}(t) \Vert_2 &\leq \sum_{(j_1,\dots,j_{i-1}) \in \mathcal{L}_i} \beta_i^{(j_1,\dots,j_{i-1})} (x_0) \prod_{m=1}^{i-1} \Vert y_m \Vert_2^{j_m} \label{eq: bound_9}
    \end{align}
    for all $t \in [i-1, i]$. Substituting these into \eqref{eq: bound_7} and using the multinomial expansion formula, we conclude that
    \begin{align}
        \forall i \in \{2,\dots,T\}\; \forall t \in [i-1, i] \;\; \left\vert\frac{\partial_+}{\partial\epsilon} c(x_i^{\epsilon}(t), u^{\epsilon}(t))\right\vert \nonumber \\
        \leq \sum_{(j_1,\dots,j_{i-1}) \in \mathcal{L}_i^{\prime}} \beta_i^{\prime(j_1,\dots,j_{i-1})} (x_0) \prod_{m=1}^{i-1} \Vert y_m \Vert_2^{j_m}
    \end{align}
    for some finite set $\mathcal{L}_i^{\prime}$ of sequences of non-negative integers and finite positive constants $\beta_i^{\prime(j_1,\dots,j_{i-1})}(x_0)$. Similarly, for $i=1$ we have $\forall t \in (\tau, 1]$
    \begin{align}
        \left\vert\frac{\partial_+}{\partial\epsilon} c(x_1^{\epsilon}(t),  u^{\epsilon}(t))\right\vert &\leq (K_6 + K_7\Vert x_1^{\epsilon}(t) \Vert_2^{L_3}) \Vert \Psi_1^{\epsilon}(t)\Vert_2 \\
        &\leq 2(K_6 + K_7 \alpha_1(x_0)^{L_3})K_1\rho_{\max}e^{K_1} \\
        &\triangleq \beta_1^{\prime}(x_0),
    \end{align}
    by Proposition \ref{prop: bound_1} and Lemma \ref{lemma: bound_2}.

    To bound the right derivative of the terminal cost, note that
    \begin{align}
        \bigg\vert\frac{\partial_+}{\partial\epsilon} &h(x^{\epsilon}(T))\bigg\vert = \left\vert\frac{\partial}{\partial x^{\epsilon}}h(x^{\epsilon}(T))^\mathrm{T} \Psi^{\epsilon}(T)\right\vert \\
        &\leq (K_6 + K_7 \Vert x^{\epsilon}(T) \Vert_2^{L_3})\Vert\Psi^{\epsilon}(T)\Vert_2 \\
        &\begin{multlined}
        = \left(K_6 + K_7 \Vert g(x_T^{\epsilon}(T), y_T) \Vert_2^{L_3}\right) \Vert\Psi^{\epsilon}(T)\Vert_2
        \end{multlined} \\
        &\begin{multlined}
        \label{eq: bound_10}
        \leq \left(K_6 + K_7 \Vert g(x_T^{\epsilon}(T), y_T) \Vert_2^{L_3}\right) \\ \times \left\Vert\frac{\partial}{\partial x_T^{\epsilon}} g(x_T^{\epsilon}(T), y_T) \right\Vert_2 \cdot \Vert \Psi_{T}^{\epsilon}(T)\Vert_2
        \end{multlined}
    \end{align}
    by Assumption \ref{assum: cost}, Proposition \ref{prop: right_differentiability_state}, and Lemma \ref{lemma: cost_variation_2}. One can apply Assumption \ref{assum: function_g_bound} to bound the norms of $g$ and its Jacobian in terms of $\Vert x_T^{\epsilon}(T) \Vert_2$ and $\Vert y_T \Vert$. Then, \eqref{eq: bound_10} becomes a polynomial of $\Vert x_T^{\epsilon}(T) \Vert_2$ and $\Vert y_T \Vert$, multiplied by $ \Vert \Psi_{T}^{\epsilon}(T)\Vert_2$. Finally, using \eqref{eq: bound_8} and \eqref{eq: bound_9} with $i = T$ to replace $\Vert x_T^{\epsilon}(T) \Vert_2$ and $ \Vert \Psi_{T}^{\epsilon}(T)\Vert_2$, one can verify that
    \begin{multline}
        \left\vert \frac{\partial_+}{\partial\epsilon} h(x^{\epsilon}(T))\right\vert \\
        \leq \sum_{(j_1,\dots,j_{T}) \in \mathcal{L}_{T+1}^{\prime}} \beta_{T+1}^{\prime(j_1,\dots,j_{T})} (x_0) \prod_{m=1}^{T} \Vert y_m \Vert_2^{j_m}
    \end{multline}
    for some finite set $\mathcal{L}^{\prime}_{T+1}$ of sequence of non-negative integers and finite positive constants $\beta_{T+1}^{\prime(j_1,\dots,j_T)}(x_0)$.
\end{proof}
\end{prop}

\begin{lemma}
    \label{lemma: cost_continuity_1}
    Let $x^{\epsilon}$ be the perturbed state induced by the perturbed control $u^{\epsilon}$, and let $(y_1,\dots,y_T)$ be the given observations. Then, the function $\epsilon\mapsto c(x^{\epsilon}(t), u^{\epsilon}(t))$ is continuous with respect to $\epsilon\in[0,\tau]$ for all $t \in (\tau, T]$.
\begin{proof}
    Note that for $t \in (\tau, T]$ we have $u^{\epsilon}(t) = u(t)$. Thus, $c(x^{\epsilon}(t), u^{\epsilon}(t)) = c(x^\epsilon(t), u(t))$. The continuity of $x_1^{\epsilon}(t)$ with respect to $\epsilon\in[0, \tau]$ follows from Lemma \ref{lemma: lipschitz_continuity}. In particular, $x_{1}^{\epsilon}(1)$ is continuous with respect to $\epsilon$. Next, suppose that $\epsilon\mapsto x_{i}^{\epsilon}(i)$ is continuous for some $i \in \{1,\dots,T\}$. Then, by Assumption \ref{assum: function_g} and Corollary \ref{cor: uniform_continuity} it follows that $\epsilon\mapsto x_{i+1}^{\epsilon}(t)$ is continuous for all $t \in [i, i+1]$. Proceeding by mathematical induction, we conclude that $x^{\epsilon}(t)$ is continuous with respect to $\epsilon \in[0,\tau]$ for all $t \in (\tau, T]$. Therefore, $\epsilon\mapsto c(x^{\epsilon}(t), u(t))$ is continuous by Assumption \ref{assum: cost}.
\end{proof}
\end{lemma}

\begin{prop}
    \label{prop: cost_mean_value_1}
    Let $u \in U$ be a control, which yields the nominal state $x$. Let $x^{\epsilon}$ be the perturbed state induced by the perturbed control $u^{\epsilon}$, and let $(y_1,\dots,y_T)$ be the given observations. Then, the following bounds hold for all $\epsilon \in [0,\tau]$:
    \begin{align}
        &\forall t \in (\tau, 1] \nonumber \\
        &\left\vert c(x_1^{\epsilon}(t), u^{\epsilon}(t)) - c(x_1(t), u(t))\right\vert \leq \epsilon \beta_1^{\prime}(x_0)
    \end{align}
    and
    \begin{align}
        &\forall i \in \{2,\dots,T\}\;\; \forall t \in [i-1, i]\nonumber \\
        &\begin{multlined}
            \label{eq: cost_difference_bound_1}
            \left\vert c(x_i^{\epsilon}(t), u^{\epsilon}(t)) - c(x_i(t), u(t))\right\vert \\
            \leq \epsilon \sum_{(j_1,\dots,j_{i-1}) \in \mathcal{L}_i^{\prime}} \beta_i^{\prime(j_1,\dots,j_{i-1})} (x_0) \prod_{m=1}^{i-1} \Vert y_m \Vert_2^{j_m},
        \end{multlined}
    \end{align}
    where $\beta_1^{\prime}(x_0)$, $\beta_i^{\prime (j_1,\dots,j_{i-1})}(x_0)$ and $\mathcal{L}^{\prime}$ are as defined in Proposition \ref{prop: bounded_cost_variations}.
\begin{proof}
    For $\epsilon \in [0,\tau]$ and $t \in (\tau, T]$, the function $\epsilon\mapsto c(x^{\epsilon}(t), u^{\epsilon}(t))$ is continuous by Lemma \ref{lemma: cost_continuity_1} and finite by Proposition \ref{prop: bounded_cost}. It is also right differentiable with respect to $\epsilon$ for $\epsilon\in[0, \tau)$ and $t \in (\tau, T]$ by Lemma \ref{lemma: cost_variation_2}. Therefore, the mean value theorem (Corollary in \cite{bourbaki2004theory}, p.15) along with Proposition \ref{prop: bounded_cost_variations} proves the claim.
\end{proof}
\end{prop}

\begin{lemma}
    \label{lemma: cost_continuity_2}
    Let $x^{\epsilon}$ be the perturbed state induced by the perturbed control $u^{\epsilon}$, and let $(y_1,\dots,y_T)$ be the given observations. Then, the function $\epsilon\mapsto h(x^{\epsilon}(T))$ is continuous with respect to $\epsilon\in[0,\tau]$.
\begin{proof}
    By the proof of Lemma \ref{lemma: cost_continuity_1} it follows that the function $\epsilon\mapsto x^{\epsilon}(T)$ is continuous with respect to $\epsilon\in[0,\tau]$. The continuity of $h$ by Assumption \ref{assum: cost} completes the proof.
\end{proof}
\end{lemma}

\begin{prop}
    \label{prop: cost_mean_value_2}
    Let $u \in U$ be a control, which yields the nominal state $x$. Let $x^{\epsilon}$ be the perturbed state induced by the perturbed control $u^{\epsilon}$, and let $(y_1,\dots,y_T)$ be the given observations. Then, the following bound holds for all $\epsilon \in [0,\tau]$:
    \begin{multline}
        \left\vert h(x^{\epsilon}(T)) - h(x(T)) \right\vert \\ \leq \epsilon\sum_{(j_1,\dots,j_{T}) \in \mathcal{L}_{T+1}^{\prime}} \beta_{T+1}^{\prime(j_1,\dots,j_{T})} (x_0) \prod_{m=1}^{T} \Vert y_m \Vert_2^{j_m},
    \end{multline}
    where $\beta_{T+1}^{\prime(j_1,\dots,j_{T})} (x_0)$ and $\mathcal{L}_{T+1}^{\prime}$ are as defined in Proposition \ref{prop: bounded_cost_variations}.
\begin{proof}
    The proof is very similar to that of Proposition \ref{prop: cost_mean_value_1}. Use Proposition \ref{prop: bounded_cost}, Lemma \ref{lemma: cost_continuity_2}, and Lemma \ref{lemma: cost_variation_2} to show finiteness, continuity, and right differentiability of $\epsilon\mapsto h(x^{\epsilon}(T))$. Then use the same mean value theorem with Proposition \ref{prop: bounded_cost_variations} to prove the claim.
\end{proof}
\end{prop}

\subsection{Expected Total Cost under Stochastic Observations}
In this last part of the analysis, we finally let the observations $(y_1,\dots,y_T)$ take random values; more formally, we treat them as a sequence of random variables $(Y_1(\omega), \dots, Y_T(\omega))$ for $\omega \in \Omega$, where $(\Omega, \mathcal{F}, \mathbb{P})$ is the probability space and each $Y_i$ satisfies Assumption \ref{assum: observations}. With $(\tau, v)$ given and fixed, $(\omega, t)$ and $\epsilon$ uniquely determine the perturbed control $u^{\epsilon}$ and the observations, hence the resulting state trajectory $x^{\epsilon}$ and the costs $c(x^{\epsilon}(t), u^{\epsilon}(t)), h(x^{\epsilon}(T))$. 

\begin{lemma}
    \label{lemma: measurable_cost_1}
    Let $\left([\tau, T], \mathcal{B}\left([\tau, T]\right), \lambda\right)$ be a measure space, where $\mathcal{B}\left([\tau, T]\right)$ is the Borel $\sigma$-algebra on $[\tau, T]$ and $\lambda$ is the Lebesgue measure. Let $\mu \triangleq \lambda\times \mathbb{P}$ be the product measure defined on the product space $\left(\Omega\times[\tau,T], \mathcal{F}\otimes\mathcal{B}\left([\tau, T]\right)\right)$, where $\mathcal{F}\otimes\mathcal{B}\left([\tau, T]\right)$ is the product $\sigma$-algebra. Then, the function $(\omega, t)\mapsto c(x^{\epsilon}(t), u^{\epsilon}(t))$ is $\mathcal{F}\otimes\mathcal{B}\left([\tau, T]\right)$-measurable for every $\epsilon \in [0, \tau]$.
\begin{proof}
    Take any $\epsilon \in [0, \tau]$. Then, $u^{\epsilon}$ is in $U$ and thus the function $(Y_1(\omega),\dots,Y_T(\omega)) \mapsto x^{\epsilon}(t)$ is continuous for every $t \in [\tau, T]$ by Proposition \ref{prop: continuity_in_observations}. Therefore, the map $(\omega, t) \mapsto x^{\epsilon}(t)$ as a function of $\omega$ is $\mathcal{F}$-measurable  for every $t \in [\tau, T]$. By Corollary \ref{cor: right_continuity}, $x^{\epsilon}(t)$ is also right continuous with respect to $t$ for every $\omega \in \Omega$. Therefore, from Theorem 3 in \cite{gowrisankaran1972measurability} it follows that $(\omega, t)\mapsto x^{\epsilon}(t)$ is measurable with respect to the product $\sigma$-algebra $\mathcal{F}\otimes\mathcal{B}\left([\tau, T]\right)$.
    
    On the other hand, $u^{\epsilon}(t)$ is piecewise continuous in $t$ and is constant with respect to $(Y_1(\omega), \dots, Y_T(\omega))$. Therefore, $(\omega, t)\mapsto u^{\epsilon}(t)$ is also measurable with respect to $\mathcal{F}\otimes\mathcal{B}\left([\tau, T]\right)$.
    
    Finally, the continuity of the instantaneous cost $c$ by Assumption \ref{assum: cost} proves the claim.
\end{proof}
\end{lemma}

\begin{prop}
   \label{prop: diff_under_integral_1}
   For the perturbed control $u^{\epsilon}$ and the perturbed state $x^{\epsilon}$, we have
   \begin{multline}
       \label{eq: diff_under_integral_1}
       \frac{\partial_+}{\partial\epsilon} \mathbb{E}\left[\int_{\tau}^{T} c(x^{\epsilon}(t), u^{\epsilon}(t)) dt\right]\Bigg\vert_{\epsilon=0} \\
       = \mathbb{E}\left[\int_{\tau}^{T} \frac{\partial}{\partial x} c(x(t), u(t))^\mathrm{T} \Psi(t) dt\right],
   \end{multline}
   where $\Psi$ is the state variation defined in Proposition \ref{prop: state_variation}.
\begin{proof}
    By definition, the left hand side of \eqref{eq: diff_under_integral_1} is
    \begin{multline}
        \frac{\partial_+}{\partial\epsilon} \mathbb{E}\left[\int_{\tau}^{T} c(x^{\epsilon}(t), u^{\epsilon}(t)) dt\right]\Bigg\vert_{\epsilon=0} \\
        =\lim_{\epsilon\rightarrow 0^{+}} \mathbb{E}\left[\int_{\tau}^{T} \frac{1}{\epsilon}\left\{c(x^{\epsilon}(t), u^{\epsilon}(t)) - c(x(t), u(t))\right\} dt\right].
    \end{multline}
    Consider the expected value above as the equivalent Lebesgue integral:
    \begin{multline}
        \int_{\Omega}\left(\int_{[\tau, T]} \frac{1}{\epsilon}\left\{c(x^{\epsilon}(t), u^{\epsilon}(t)) - c(x(t), u(t))\right\} d\lambda(t)\right)\\ \times d\mathbb{P}(\omega).
    \end{multline}
    By Lemma \ref{lemma: measurable_cost_1} the integrand is measurable in the product space. In addition, Proposition \ref{prop: cost_mean_value_1} shows that the absolute value:
    \begin{align}
        \frac{1}{\epsilon}\left\vert c(x^{\epsilon}(t), u^{\epsilon}(t)) - c(x(t), u(t))\right\vert
    \end{align}
    is bounded by an integrable function for $\mu$-a.e. $(\omega, t)$. Indeed, if we let $\hat{c}(\omega, t)$ to be a function defined by
    \begin{align}
        \hat{c}(\omega, t) = \begin{cases} 
            \beta_1^{\prime}(x_0) \;\; \forall t \in [\tau, 1) \\
            \sum_{(j_1,\dots,j_{i-1}) \in \mathcal{L}_i^{\prime}} \beta_i^{\prime(j_1,\dots,j_{i-1})} (x_0) \prod_{m=1}^{i-1} \Vert Y_m \Vert_2^{j_m} \\ \;\;\;\;\;\;\;\;\;\;\;\;\forall t \in [i-1, i)\;\forall i \in \{1,\dots, T\},
        \end{cases}
    \end{align}
    then we have
    \begin{align}
        \frac{1}{\epsilon}\left\vert c(x^{\epsilon}(t), u^{\epsilon}(t)) - c(x(t), u(t))\right\vert \leq \hat{c}(\omega, t)
    \end{align}
    for every non-zero $\epsilon$ and $\mu$-a.e. $(\omega, t)$, and
    \begin{multline}
        \int_{\Omega}\left(\int_{[\tau, T]} \hat{c}(\omega, t) d\lambda(t)\right) d\mathbb{P}(\omega) = \beta_1^{\prime}(x_0)(1 - \tau) \\ + \sum_{i=2}^{T}\sum_{(j_1,\dots,j_{i-1}) \in \mathcal{L}_i^{\prime}} \beta_i^{\prime(j_1,\dots,j_{i-1})} (x_0) \\ \times \mathbb{E}\left[\prod_{m=1}^{i-1} \Vert Y_m \Vert_2^{j_m}\right],
    \end{multline}
    where
    \begin{align}
        \mathbb{E}\left[\prod_{m=1}^{i-1} \Vert Y_m \Vert_2^{j_m}\right] &\leq \sqrt{\mathbb{E}\left[\Vert Y_1 \Vert_2^{j_1}\right]} \sqrt{\mathbb{E}\left[\prod_{m=2}^{i-1} \Vert Y_m \Vert_2^{j_m}\right]} \\
        &\leq \nonumber \\
        &\vdots \nonumber \\
        &\leq\prod_{m=1}^{i-1} \left(\mathbb{E}\left[\Vert Y_m \Vert_2^{j_m}\right]\right)^{\frac{1}{2^m}} < \infty
    \end{align}
    by the Cauchy-Schwarz inequality and Assumption \ref{assum: observations}.
    
    Furthermore, Lemma \ref{lemma: cost_variation_2} proves that
    \begin{multline}
        \lim_{\epsilon\rightarrow 0^{+}} \frac{1}{\epsilon}\left\{c(x^{\epsilon}(t), u^{\epsilon}(t)) - c(x(t), u(t))\right\} \\
         = \frac{\partial}{\partial x} c(x(t), u(t))^\mathrm{T} \Psi(t)
    \end{multline}
    for $\mu$-a.e. $(\omega, t)$. Therefore, the dominated convergence theorem yields
    \begin{multline}
        \frac{\partial_+}{\partial\epsilon} \int_{\Omega}\left(\int_{\tau}^{T} c(x^{\epsilon}(t), u^{\epsilon}(t)) d\lambda(t) \right) d\mathbb{P}(\omega)  \Bigg\vert_{\epsilon=0} \\
        =  \int_{\Omega}\left(\int_{\tau}^{T}  \frac{\partial}{\partial x} c(x(t), u(t))^\mathrm{T} \Psi(t) d\lambda(t) \right) d\mathbb{P}(\omega).
    \end{multline}
    
\end{proof}
\end{prop}

\begin{prop}
    \label{prop: diff_under_integral_2}
    For the perturbed control $u^{\epsilon}$ and the perturbed state $x^{\epsilon}$, we have
    \begin{align}
        \label{eq: diff_under_integral_2}
        \frac{\partial_+}{\partial\epsilon} \mathbb{E}\left[h(x^{\epsilon}(T)\right]\bigg\vert_{\epsilon=0} = \mathbb{E}\left[\frac{\partial}{\partial x}h(x(T))^\mathrm{T} \Psi(T)\right],
    \end{align}
    where $\Psi$ is the state variation defined in Proposition \ref{prop: state_variation}.
\begin{proof}
    The proof is similar to that of Proposition \ref{prop: diff_under_integral_1}. By Assumption \ref{assum: cost} and Proposition \ref{prop: continuity_in_observations}, $(Y_1(\omega), \dots, Y_T(\omega))\mapsto h(x^\epsilon(T))$ is a continuous map for every $\epsilon \in [0, \tau]$. Therefore, the function $\omega \mapsto h(x^\epsilon(T))$ is $\mathcal{F}$-measurable.
    
    By definition, the left hand side of \eqref{eq: diff_under_integral_2} is
    \begin{multline}
        \frac{\partial_+}{\partial\epsilon} \mathbb{E}\left[h(x^{\epsilon}(T)\right]\bigg\vert_{\epsilon=0} \\
        = \lim_{\epsilon\rightarrow 0^{+}} \int_{\Omega} \frac{1}{\epsilon}\left\{h(x^{\epsilon}(T)) - h(x(T))\right\} d\mathbb{P}(\omega)
    \end{multline}
    The analysis above implies that the integrand is $\mathcal{F}$-measurable. In addition, Proposition \ref{prop: cost_mean_value_2} shows that the absolute value:
    \begin{align}
        \frac{1}{\epsilon}\left\vert h(x^{\epsilon}(T)) - h(x(T))\right\vert
    \end{align}
    is bounded by an integrable function for every non-zero $\epsilon$ and every $\omega \in \Omega$. Furthermore, Lemma \ref{lemma: cost_variation_2} proves that
    \begin{align}
        \lim_{\epsilon\rightarrow 0^{+}} \frac{1}{\epsilon}\left\{h(x^{\epsilon}(T)) - h(x(T))\right\} = \frac{\partial}{\partial x} h(x(T))^\mathrm{T} \Psi(T)
    \end{align}
    for every $\omega \in \Omega$. Therefore, the dominated convergence theorem yields
    \begin{align}
        \frac{\partial_+}{\partial\epsilon} \int_{\Omega} h(x^{\epsilon}(T)) d\mathbb{P}(\omega)\bigg\vert_{\epsilon=0} = \int_{\Omega}\frac{\partial}{\partial x} h(x(T))^\mathrm{T} \Psi(T) d\mathbb{P}(\omega).
    \end{align}
\end{proof}
\end{prop}

\begin{theorem}[Mode Insertion Gradient]
    Suppose that Assumptions \ref{assum: control} -- \ref{assum: observations} are satisfied. For a given $(\tau, v)$, let $u^{\epsilon}$ denote the perturbed control of the form \eqref{eq: perturbed_control}. The perturbed control $u^{\epsilon}$ and the stochastic observations $(Y_1, \dots, Y_T)$ result in the stochastic perturbed state trajectory $x^{\epsilon}$. For such $u^{\epsilon}$ and $x^{\epsilon}$, let us define the mode insertion gradient of the expected total cost as
    \begin{align}
        \frac{\partial_+}{\partial\epsilon} \mathbb{E}\left[\int_{0}^{T} c(x^{\epsilon}(t), u^{\epsilon}(t)) dt + h(x^{\epsilon}(T))\right]\Bigg\vert_{\epsilon=0}.
    \end{align}
    Then, this right derivative exists and we have
    \begin{multline}
         \frac{\partial_+}{\partial\epsilon} \mathbb{E}\left[\int_{0}^{T} c(x^{\epsilon}(t), u^{\epsilon}(t)) dt + h(x^{\epsilon}(T))\right]\Bigg\vert_{\epsilon=0} \\
          = c(x(\tau), v) - c(x(\tau), u(\tau)) \\ +
          \mathbb{E}\Bigg[\int_{\tau}^{T} \frac{\partial}{\partial x} c(x(t), u(t))^\mathrm{T} \Psi(t) dt   \\ + \frac{\partial}{\partial x} h(x(T))^\mathrm{T}\Psi(T)\Bigg],
    \end{multline}
    where $\Psi$ is the state variation defined in Proposition \ref{prop: state_variation}.
\begin{proof}
    We first consider the instantaneous cost $c$. Split the integration interval to get
    \begin{multline}
        \label{eq: split_expected_cost}
        \mathbb{E}\left[\int_0^{T} c(x^{\epsilon}(t), u^{\epsilon}(t)) dt\right] \\=
        \mathbb{E}\left[\int_0^{\tau-\epsilon} c(x^{\epsilon}(t), u^{\epsilon}(t)) dt \right] \\ +
        \mathbb{E}\left[\int_{\tau-\epsilon}^{\tau} c(x^{\epsilon}(t), u^{\epsilon}(t)) dt\right] \\+ \mathbb{E}\left[\int_{\tau}^{T} c(x^{\epsilon}(t), u^{\epsilon}(t)) dt\right]
    \end{multline}
    For the first two terms in the sum, recall that the evolution of the state $x^{\epsilon}(t)$ is not affected by any observations for all $t \in [0, \tau]$. Thus,
    \begin{align}
        \mathbb{E}\left[\int_0^{\tau-\epsilon} c(x^{\epsilon}(t), u^{\epsilon}(t)) dt\right] &= \int_0^{\tau-\epsilon} c(x^{\epsilon}(t), u^{\epsilon}(t)) dt \label{eq: integrated_instantaneous_cost_constant}\\
        \mathbb{E}\left[\int_{\tau-\epsilon}^{\tau} c(x^{\epsilon}(t), u^{\epsilon}(t)) dt\right] &= \int_{\tau-\epsilon}^{\tau} c(x^{\epsilon}(t), u^{\epsilon}(t)) dt.
        \label{eq: integrated_instantaneous_cost_deterministice}
    \end{align}
    Note that \eqref{eq: integrated_instantaneous_cost_constant} is constant with respect to $\epsilon$, since for all $t \in [0, \tau-\epsilon]$ we have $u^{\epsilon}(t) = u(t)$ and $x^{\epsilon}(t) = x(t)$. On the other hand, for \eqref{eq: integrated_instantaneous_cost_deterministice} we can apply Lemma \ref{lemma: cost_variation_1} to obtain
    \begin{multline}
        \frac{\partial_+}{\partial \epsilon}\mathbb{E}\left[\int_{\tau-\epsilon}^{\tau} c(x^{\epsilon}(t), u^{\epsilon}(t)) dt\right]\Bigg\vert_{\epsilon=0} \\= c(x(\tau), v) - c(x(\tau), u(\tau))
    \end{multline}
    
    For the last term, Proposition \ref{prop: diff_under_integral_1} gives
    \begin{multline}
        \frac{\partial_+}{\partial\epsilon} \mathbb{E}\left[\int_{\tau}^{T} c(x^{\epsilon}(t), u^{\epsilon}(t)) dt\right]\Bigg\vert_{\epsilon=0} \\
       = \mathbb{E}\left[\int_{\tau}^{T} \frac{\partial}{\partial x} c(x(t), u(t))^\mathrm{T} \Psi(t) dt\right].
    \end{multline}
    
    Finally, for the terminal cost $h$ we have 
    \begin{align}
        \frac{\partial_+}{\partial\epsilon} \mathbb{E}\left[h(x^{\epsilon}(T)\right]\bigg\vert_{\epsilon=0} = \mathbb{E}\left[\frac{\partial}{\partial x}h(x(T))^\mathrm{T} \Psi(T)\right]
    \end{align}
    by Proposition \ref{prop: diff_under_integral_2}.
\end{proof}
\end{theorem}

\begin{remark}[Closed-loop Nominal Policy]
    \label{remark: closed_loop_policy}
    As far as the control is concerned, the analysis above only requires that the nominal control $u$ is in $U$ (as in Assumption \ref{assum: control}) and that the perturbed control $u^{\epsilon}$ is measurable with respect to $\mathcal{F} \otimes \mathcal{B}([\tau, T])$ (as in Lemma \ref{lemma: measurable_cost_1}). To satisfy these requirements with a closed-loop nominal policy $\pi: \mathbb{R}^{n_x} \rightarrow \mathbb{R}^m$, it is sufficient that $\pi$ is a measurable map and that the induced nominal control trajectory $u(t) = \pi(x(t))$ for $t \in [0, T]$ belongs to $U$ for any given observations $(y_1,\dots,y_T)$. 
    
    Notice that even with state feedback, both the nominal control trajectory $u \in U$ and the nominal state trajectory $x$ are uniquely determined under a fixed sequence of observations. Note also that the type of the control perturbation considered in this sensitivity analysis is still open-loop:
    \begin{align}
        u^{\epsilon}(t) = \begin{cases}
            v & \text{if}\;\; t \in (\tau - \epsilon, \tau] \\
            \pi(x(t)) & \text{otherwise},
        \end{cases}
    \end{align}
    because we still follow Definition \ref{def: perturbed_control} for the control perturbation model. That is, the nominal state trajectory $x$ is used in the control feedback. This is not to be confused with the closed-loop perturbation:
    \begin{align}
        u^{\epsilon}_\text{closed}(t) = \begin{cases}
            v & \text{if}\;\; t \in (\tau - \epsilon, \tau] \\
            \pi(x^{\epsilon}(t)) & \text{otherwise},
        \end{cases}
    \end{align}
    where the perturbed state trajectory $x^{\epsilon}$ is fed back to the perturbed control.
\end{remark}

\end{document}